\documentclass[12pt]{article}
\usepackage{amsmath,amsthm,amssymb,enumerate,mathtools}
\usepackage{microtype}
\usepackage{color}
\usepackage{pxfonts}
\newtheorem{thm}{Theorem}[section]

\newtheorem{lem}[thm]{Lemma}
\newtheorem{prop}[thm]{Proposition}

\newtheorem{rmk}[thm]{Remark}

\newcommand{\cD}{{\mathcal D}}
\newcommand{\cX}{{\mathcal X}}

\newcommand{\eps}{\varepsilon}

\newcommand{\R}{\mathbb{R}}
\newcommand{\dd}{\,\mathrm{d}}
\newcommand{\calE}{\mathcal{E}}

\newcommand{\abs}[1]{\left\lvert #1\right\rvert}

\def\ni{\noindent}
\begin{document}
\title{A matching construction of self-similar profiles for the fast diffusion equation}
\author{Kin Ming Hui\\
Institute of Mathematics, Academia Sinica\\
Taipei, Taiwan, R. O. C.\\
e-mail: kmhui@as.edu.tw}
\date{Aug 31, 2026}
\smallbreak \maketitle
\begin{abstract}
Let $n\ge 3$, $0<m<\frac{n-2}{n}$, $\eta>0$, $\eta_0>0$, $\rho_1>0$, $\beta_-(\rho_1)=-\frac{\rho_1}{2}$, $\beta_+(\rho_1)=\frac{m\rho_1}{n-2-nm}$ and $\alpha=\frac{2\beta+\rho_1}{1-m}$. For any $\beta_-(\rho_1)\le\beta\le\beta_+(\rho_1)$, we construct the unique maximal positive radial branch of
\[
 \Delta(f^m/m)+\alpha f+\beta x\cdot\nabla f=0
\]
 issuing from prescribed  origin data  $f(0)=\eta_0$ and $f_r(0)=0$. For any $\beta\le\beta_+(\rho_1)$, we construct
the unique maximal positive radial branch at infinity satisfying
\[
\lim_{|x|\to\infty}|x|^{\frac{n-2}{m}}f(x)=\eta.
\]
  The unmatched branches may reach zero at a finite radius. 
We formulate the shooting construction through the slope--amplitude equations of the increasing and decreasing half-branches before their first turning points.   As a consequence we obtain a new proof of the existence result of Peletier and Zhang \cite{PeZ}: there exists $\beta\in (\beta_-(\rho_1),\beta_+(\rho_1))$ for which the equation $\Delta(f^m/m)+\alpha f+\beta x\cdot\nabla f=0$, $f>0$, in $\mathbb{R}^n$ has a positive radial solution $f$ satisfying
\[ f(0)=\eta _0,\qquad f_r(0)=0,
 \qquad
 \lim_{r\to\infty}r^{\frac{n-2}{m}}f(r)
 =C_*\rho _1^{-\frac{n-2}{2m}}
   \eta _0^{-\frac{n-2-nm}{2m}}.
\]
for some constant $C_*>0$ depending on  $n$, $m$, and is independent of $\rho_1$ and $\eta_0$. 
For every selected matching value of $\beta$, this solution is unique among positive radial solutions with the prescribed origin data.  When $\rho _1=1$, the function $V(x,t)=(T-t)^\alpha f((T-t)^\beta x)$ is a backward self-similar solution of $u_t=\Delta(u^m/m)$ in $\mathbb{R}^n\times (-\infty,T)$.
\end{abstract}

\vskip 0.2truein



\vskip 0.2truein
\setcounter{equation}{0}
\setcounter{section}{0}

\section{Introduction}
\setcounter{equation}{0}
\setcounter{thm}{0}

There has been considerable study of the equation
\begin{equation}\label{fde}
u_t=\Delta (u^m/m).
\end{equation}
When $m=1$, \eqref{fde} is the heat equation.  When $m>1$, it is the porous medium equation, which models diffusion through porous media \cite{A} and also arises in the large-time analysis of compressible Euler equations with damping \cite{LZ}.  When $n\ge 3$ and $m=\frac{n-2}{n+2}$, \eqref{fde} is related to the Yamabe flow \cite{DKS}, \cite{DS1}, \cite{DS2}, \cite{PS}.  When $0<m<1$, it is the fast diffusion equation, which also appears in models of impurity diffusion in silicon \cite{K1}, \cite{K2}.  We refer to the books of P.~Daskalopoulos and C.E.~Kenig \cite{DK} and J.L.~V\'azquez \cite{V} for background and further results.

As observed by P.~Daskalopoulos and N.~Sesum \cite{DS1}, M.~del Pino and M.~S\'aez \cite{PS}, and J.L.~V\'azquez \cite{V}, when $n\ge 3$ and $0<m<\frac{n-2}{n}$, there are solutions of \eqref{fde} in $\mathbb{R}^n\times(0,T)$ that vanish at a finite time $T>0$.  Under appropriate hypotheses on the initial data, such solutions approach backward self-similar profiles near the extinction time.  It is therefore natural to study backward self-similar solutions of \eqref{fde}.

Note that a function
\begin{equation*}
V(x,t)=(T-t)^{\alpha}f((T-t)^{\beta} x)\quad\forall x\in\mathbb{R}^n, t<T
\end{equation*} 
is a self-similar solution of \eqref{fde} in $\mathbb{R}^n\times (-\infty,T)$ if and only if
$f$ satisfies
\begin{equation}\label{f-elliptic-eqn}
\Delta (f^m/m)+\alpha f+\beta x\cdot\nabla f=0,\quad f>0,
\end{equation}
 in $\mathbb{R}^n$ with
\begin{equation}\label{alpha-beta-relation}
\alpha=\frac{2\beta+1}{1-m}.
\end{equation}
When $m=\frac{n-2}{n+2}$ and $\alpha$, $\beta$, satisfy \eqref{alpha-beta-relation},  \eqref{f-elliptic-eqn} also appears in the study of locally conformally flat shrinking Yamabe solitons \cite{DS2}. In \cite{Hs1} S.Y.~Hsu proved that
for any 
\begin{equation}\label{n-m-range}
n\ge 3,\quad 0<m<\frac{n-2}{n},
\end{equation}
and
\begin{equation}\label{alpha-beta-relation2}
 \rho_1>0,\quad\beta\ge\frac{m\rho_1}{n-2-nm}\quad\mbox{ and }\quad\alpha=\frac{2\beta+\rho_1}{1-m},
\end{equation}
 there exists a unique radially symmetric solution $f$ of \eqref{f-elliptic-eqn}.
In \cite{Hs2} S.Y.~Hsu proved that when
\begin{equation*}
\rho_1>0,\quad\beta>\frac{m\rho_1}{n-2-nm}\quad\mbox{ and }\quad\alpha=\frac{2\beta+\rho_1}{1-m},
\end{equation*}
holds,
such radially symmetric solution $f$ satisfies
\begin{equation}\label{f-infty-behaviour1}
\lim_{r\to\infty}r^2f(r)^{1-m}=\frac{2(n-2-nm)}{(1-m)(\alpha(1-m)-2\beta)}.
\end{equation}
By \eqref{n-m-range} and \eqref{f-infty-behaviour1}, such solution $f$ satisfies $f\not\in L^1(\mathbb{R}^n)$. On the other hand in the paper \cite{PeZ} M.A.~Peletier and H.~Zhang used phase plane analysis to prove that when \eqref{n-m-range} holds, then  there exists 
$\alpha=\frac{2\beta+1}{1-m}$ and $\beta=\beta (n,m)$ satisfying
\begin{equation}\label{alpha-beta-relation3}
-\frac{1}{2}<\beta<\frac{m}{n-2-nm}
\end{equation}
such that \eqref{f-elliptic-eqn} has a positive radially symmetric solution $f$ which satisfies 
\begin{equation}\label{r-fr-f-ratio=limit1}
\frac{rf_r(r)}{f(r)}\to-\frac{n-2}{m}\quad\mbox{ as }r\to\infty
\end{equation}
and
\begin{equation}\label{r-fr-f-ratio=limit2}
\frac{rf_r(r)}{f(r)}\to 0\quad\mbox{ as }r\to 0.
\end{equation}
Note that by \eqref{r-fr-f-ratio=limit1} and \eqref{r-fr-f-ratio=limit2}, $f\in L^1(\mathbb{R}^n)$. For any $\rho_1>0$, let
\begin{equation*}
\beta_-(\rho_1)=-\frac{\rho_1}{2}\quad\mbox{ and }\quad\beta_+(\rho_1)=\frac{m\rho_1}{n-2-nm}.
\end{equation*}
In this paper we will extend the results of \cite{Hs1} and \cite{PeZ} and prove that when 
\eqref{n-m-range} holds, then for any  $\eta_0>0$, $\rho_1>0$, $\alpha=\frac{2\beta+\rho_1}{1-m}$ and $\beta\in [\beta_-(\rho_1),\beta_+(\rho_1)]$, there exists a unique maximal positive radial branch of
\eqref{f-elliptic-eqn}
in $\mathbb{R}^n$ issuing from prescribed positive data at the origin satisfying 
\begin{equation}\label{f-fr-origin-value}
f(0)=\eta_0,\quad f_r(0)=0,
\end{equation} 
 and for any $\beta\le\beta_+(\rho_1)$ we construct
the unique maximal positive branch at infinity satisfying 
\begin{equation}\label{f-infty-behaviour2}
\lim_{r\to\infty}r^{\frac{n-2}{m}}f(r)=\eta.
\end{equation}
 The unmatched branches may reach zero at a finite radius.

A natural question is whether there exists a constant $\beta$  such that \eqref{f-elliptic-eqn} with $\alpha=\frac{2\beta+\rho_1}{1-m}$ and $\rho_1>0$ has a radially symmetric solution in $\mathbb{R}^n$ which satisfies 
\eqref{f-fr-origin-value} and \eqref{f-infty-behaviour2} for some constant $\eta$ that may depend on $\eta_0$. In the paper \cite{PeZ}  Peletier and Zhang gave an affirmative answer
to this question. They used a phase-plane method to construct such solution. Their proof uses a different phase-plane formulation. 

In this paper, using a matching argument, we give a new proof of the result on the existence of $\beta\in (\beta_-(\rho_1),\beta_+(\rho_1))$ such that the equation \eqref{f-elliptic-eqn} in $\mathbb{R}^n$ with $\alpha=\frac{2\beta+\rho_1}{1-m}$, $\rho_1>0$,
has a unique radially symmetric solution which satisfies 
\begin{equation}\label{f-infty-limit}
\lim_{|x|\to\infty}|x|^{\frac{n-2}{m}}f(x)=C_*\rho _1^{-\frac{n-2}{2m}}
   \eta _0^{-\frac{n-2-nm}{2m}}
\end{equation}
for some constant $C_*>0$ depending on  $n$, $m$, and is independent of $\rho_1$ and $\eta_0$. 
Our proof does not use the Peletier--Zhang \((t,z)\)-transformation, their
phase-plane half-orbits, or any result from that paper.  

Instead we formulate the shooting construction in terms of the corresponding slope--amplitude equations of the increasing and decreasing half-branches $f^o_{\beta}$ and $f^{\infty}_{\beta}$ before their first turning points.
We prove continuity of their first turning amplitudes, match them at a common nondegenerate maximum, and obtain a positive smooth pulse of the autonomous $W$-equation by gluing the two half-orbits. Transforming back yields a positive radial solution, unique for the selected parameter and origin data, of \eqref{f-elliptic-eqn} in $\mathbb{R}^n$ which satisfies \eqref{f-fr-origin-value} and \eqref{f-infty-limit}.

Note that when $n\ge 3$, $m=\frac{n-2}{n+2}$ and $\beta=0$, $\alpha=\frac{2\beta+1}{1-m}=\frac{n+2}{4}$, such function $f$ (cf. \cite{GP}, \cite{PS}) is given explicitly by
\begin{equation*}
 f(x)=\left(A+\frac{|x|^2}{4nA}\right)^{-\frac{n+2}{2}},\qquad A>0
\end{equation*}
which satisfies
\begin{equation*}
\Delta (f^m/m)+\frac{n+2}{4}f=0,\quad f>0, \quad\mbox{ in }\mathbb{R}^n
\end{equation*}
and
\begin{equation*}
\lim_{|x|\to\infty}|x|^{\frac{n-2}{m}}f(x)
 =(4nA)^{\frac{n+2}{2}}.
\end{equation*}
In this paper we will prove the following main results.

\begin{thm}[Maximal branch issuing from the origin]\label{thm-existence-self-similar1}
Let $n\ge 3$, $0<m<\frac{n-2}{n}$, $\eta_0>0$, $\rho_1>0$, $\alpha=\frac{2\beta+\rho_1}{1-m}$ and $\beta\in[\beta_-(\rho_1),\beta_+(\rho_1)]$. Then there exist a unique constant
$R_o=R_o(\beta,\eta _0)\in(0,\infty]$
and a unique maximal positive radial solution
\[
 f^o_\beta\in C^1([0,R_o))\cap C^2((0,R_o))
\]
of \eqref{f-elliptic-eqn}  in $B_{R_o}$ satisfying \eqref{f-fr-origin-value}. Its radial extension $x\mapsto f^o_\beta(|x|)$ belongs to $C^\infty(B_{R_o})$.
If $\beta=\beta_-(\rho _1)$, then $R_o=\infty$ and $f^o_\beta\equiv\eta _0$.  If
$\beta\in(\beta_-(\rho _1),\beta_+(\rho _1)]$, then
\begin{equation}\label{f-derivative-negative2}
f^o_{\beta,r}(r)<0\quad\forall 0<r<R_o.
\end{equation}
If $R_o<\infty$, then
\begin{equation}\label{origin-finite-zero}
 f^o_\beta(r)\searrow 0
 \qquad\text{as }r\nearrow R_o,
\end{equation}
and the endpoint flux is
\begin{equation}\label{origin-endpoint-flux}
 \begin{aligned}
 L_0
 &:=\lim_{r\nearrow R_o}
 r^{n-1}(f^o_{\beta})^{m-1}(r)f^o_{\beta,r}(r)\\
 & =-(\alpha-n\beta)
 \int_0^{R_o}\rho^{n-1}f^o_{\beta}(\rho)\,\dd\rho<0.
 \end{aligned}
\end{equation}
Consequently,
\begin{equation}\label{origin-finite-zero-asymptotic}
f^o_{\beta}(r)^m=-\frac{mL_0}{R_o^{n-1}}(R_o-r)+o(R_o-r),
 \qquad\mbox{ as }\quad r\nearrow R_o.
\end{equation}

\end{thm}

\begin{thm}[Maximal branch issuing from infinity]\label{thm-existence-self-similar2}
Let $n\ge 3$, $0<m<\frac{n-2}{n}$, $\eta>0$, $\rho_1>0$, $\alpha=\frac{2\beta+\rho_1}{1-m}$ and $\beta\in (-\infty,\beta_+(\rho_1)]$.
Let
\begin{equation}\label{q-and-tilde-def-intro}
\widetilde\beta=-\beta,\quad\mbox{and}\quad\widetilde\alpha=\alpha-\frac{n-2}{m}\,\beta=a_0(\beta_+(\rho_1)-\beta)\quad\mbox{ where }\quad a_0=\frac{n-2-nm}{(1-m)m}.
\end{equation}
Then there exists a unique constant $R_g=R_g(\beta,\eta)\in(0,\infty]$ and a unique maximal positive solution $g_\beta$ of the inverted equation
\begin{equation}\label{eq-inversion}
(g^m/m)_{rr}+\frac{n-1}{r}(g^m/m)_r+r^{\frac{n-2-nm}{m}-2}\left(\widetilde\alpha g+\widetilde\beta rg_r\right)=0,\quad g>0,
\end{equation}
in $[0,R_g)$ satisfying $g_\beta(0)=\eta$. If $0<m<\frac{n-2}{n+1}$, then $g_\beta\in C^1([0,R_g))\cap C^2((0,R_g))$ with $g_{\beta,r}(0)=0$. If $\frac{n-2}{n+1}\le m<\frac{n-2}{n}$, then 
$g_\beta\in C_{loc}^{0, \delta_0}([0,R_g))\cap C^2((0,R_g))$ with
\begin{equation}\label{eq-initial-m-large}
\lim_{r\to 0^+}r^{\delta_1} g_r(r)=-\zeta
\end{equation} 
where 
\begin{equation}\label{defn-delta0-1}
\delta_1=1-\frac{n-2-nm}{m}\in[0,1),\quad\delta_0=\frac{1-\delta_1}{2}=\frac{n-2-nm}{2m}\in(0,1/2]
\end{equation} 
and
\begin{equation}\label{zeta-defn}
\zeta:= \frac{\widetilde\alpha m\eta^{2-m}}{n-2-2m}.
\end{equation}
If $\beta=\beta_+(\rho _1)$, then $R_g=\infty$ and
$g_\beta\equiv\eta$.  If $\beta<\beta_+(\rho _1)$, then
\begin{equation}\label{g-monotonicity-intro}
 g_{\beta,r}(r)<0
 \qquad\forall 0<r<R_g.
\end{equation}
With the convention $R_g^{-1}=0$ when $R_g=\infty$, the function
\begin{equation}\label{outer-profile-definition-intro}
 f^{\infty}_\beta(r):=r^{-\frac{n-2}{m}}g_\beta(r^{-1})
 \qquad\forall r>R_g^{-1}
\end{equation}
is the unique maximal positive outer radial branch solution of \eqref{f-elliptic-eqn}  satisfying \eqref{f-infty-behaviour2}.  If $\beta<\beta_+(\rho _1)$, then 
\begin{equation}\label{outer-branch-inequalities-intro}
 f^{\infty}_\beta+\frac{m}{n-2}rf^{\infty}_{\beta,r}>0,
 \qquad
 0<f^{\infty}_\beta(r)<\eta r^{-\frac{n-2}{m}}\qquad\forall r>R_g^{-1},
\end{equation}
and
\begin{equation}\label{outer-derivative-limit-intro}
 \lim_{r\to\infty}r^{\frac{n-2}{m}+1}f^{\infty}_{\beta,r}(r)=-\frac{n-2}{m}\eta.
\end{equation}
If $R_g<\infty$, then $g_\beta(r)\searrow 0$ as $r\nearrow R_g$ and $f^{\infty}_\beta(r)\to 0$ as
$r\searrow R_g^{-1}$.

\end{thm}

\begin{thm}[Matched global pulse]\label{thm-existence-self-similar3}
Let $n\ge 3$, $0<m<\frac{n-2}{n}$, $\eta_0>0$ and $\rho_1>0$. Then there exists  $\beta\in(\beta_-(\rho_1),\beta_+(\rho_1))$ such that \eqref{f-elliptic-eqn} has a unique positive radially symmetric solution $f$ in $\mathbb{R}^n$ with $\alpha=\frac{2\beta+\rho_1}{1-m}$,
\begin{equation}\label{matched-origin-data}
 f(0)=\eta _0,
 \qquad f_r(0)=0,
 \qquad f_r(r)<0\quad\forall r>0
\end{equation} 
and satisfying \eqref{f-infty-limit} for some constant $C_*>0$ depending on  $n$, $m$, and is independent of $\rho_1$ and $\eta_0$. 
\end{thm}

Unless stated otherwise we will assume that $n\ge 3$, $0<m<\frac{n-2}{n}$, $\eta>0$, $\eta_0>0$, $\rho_1>0$, $\alpha=\frac{2\beta+\rho_1}{1-m}$ and $\beta\le\beta_+(\rho_1)$
for the rest of the paper. For any $R>0$, we let $B_R=\{x:|x|<R\}\subset\mathbb{R}^n$ with $B_{\infty}=\mathbb{R}^n$. Note that 
if $\beta\in (\beta_-(\rho_1),\beta_+(\rho_1)]$, then $\alpha>0$.  

\begin{rmk}[No uniqueness assertion for the matching parameter]\label{rmk:matching-parameter-nonunique}
The proof below establishes the existence of at least one matching value
$\beta$.  It does not prove that this value is unique.  
\end{rmk}

\begin{rmk}\label{rmk:alpha.n-beta}
Note that 
\begin{equation*}
\alpha>n\beta\quad\Leftrightarrow\quad\beta<\frac{\rho_1}{n-2-nm}.
\end{equation*}
\end{rmk}

The plan of the paper is as follows. In Section 2 we prove Theorem \ref{thm-existence-self-similar1} in its maximal-positive form. In Section 3 we prove the analogous exterior result, Theorem \ref{thm-existence-self-similar2}. In Section 4 we derive the slope--amplitude equations on the two first-turn half-branches and prove continuity of their turning levels without assuming global positivity of the unmatched profiles. In Section 5 we match and glue those half-branches to prove Theorem \ref{thm-existence-self-similar3}.

\section{The maximal positive branch from the origin}
\setcounter{equation}{0}
\setcounter{thm}{0}

In this section we will prove Theorem \ref{thm-existence-self-similar1}. We first observe that if $f$ is a radially symmetric solution of \eqref{f-elliptic-eqn} in $\mathbb{R}^n\setminus\{0\}$, then 
\begin{equation}\label{f-ode}
(f^m/m)_{rr}+\frac{n-1}{r}(f^m/m)_r+\alpha f+\beta rf_r=0,\quad f>0,\quad\forall r>0.
\end{equation}
Multiplying \eqref{f-ode} by $r^{n-1}$ and integrating over $(0,r)$, if $f$ satisfies \eqref{f-fr-origin-value}, we get
\begin{align*}
&r^{n-1}f(r)^{m-1}f_r(r)+\int_0^r\rho^{n-1}(\alpha f(\rho)+\beta \rho f_r(\rho))\,d\rho=0\quad\forall r>0\notag\\
\Rightarrow\quad&f_r(r)=-\frac{f(r)^{1-m}}{r^{n-1}}\int_0^r\rho^{n-1}(\alpha f(\rho)+\beta \rho f_r(\rho))\,d\rho\quad\forall r>0
\end{align*}
which suggests one to use fixed point argument to prove Theorem \ref{thm-existence-self-similar1}.

\begin{lem}[Local existence, uniqueness, and the origin expansion]\label{lem-local-existence}
Let $n\ge 2$, $0<m<1$, $\eta_0>0$, $\alpha,\beta\in\mathbb{R}$.
Then there exists a constant $\eps>0$ such that \eqref{f-ode} has a unique positive solution $f\in C^1([0,\eps])\cap C^2((0,\eps])$ which satisfies \eqref{f-fr-origin-value}. Moreover
\begin{equation}
 r^{n-1}f^{m-1}(r)f_r(r)
=-(\alpha-n\beta)\int_0^r\rho^{n-1}f(\rho)\,d\rho
   -\beta r^nf(r)\label{fr-representation-formula0}
\end{equation}
or equivalently
\begin{equation}\label{fr-representation-formula}
f_r(r)=-(\alpha-n\beta)\frac{f(r)^{1-m}}{r^{n-1}}\int_0^r\rho^{n-1}f(\rho)\,d\rho-\beta rf(r)^{2-m}
\end{equation}
holds for any $0<r\le\eps$ and
\begin{equation}\label{origin-second-order-expansion}
 \lim_{r\searrow0}\frac{f_r(r)}r
 =-\frac{\alpha}{n}\eta _0^{2-m}.
\end{equation}
The construction and the limit in \eqref{origin-second-order-expansion} are uniform when $(\alpha,\beta)$ ranges over a compact set.  The radial extension $x\mapsto f(|x|)$ is smooth in $B_{\eps}$ and $f$ satisfies \eqref{f-elliptic-eqn} in $B_{\eps}$.
\end{lem}
\begin{proof} 
We will use a modification of the proof of Lemma 2.2 of \cite{HuS} to prove this lemma. Let $0<\eps<1$. We define the Banach space 
\begin{equation*}
\cX_\eps:=\left\{(f,h): f, h\in C( [0,\eps])\,\mbox{ such that } ||(f,h)||_{\cX_\eps}<\infty\right\}
\end{equation*}
with a norm given by
\[||(f,h)||_{\cX_\eps}=\max\left\{\|f\|_{L^\infty([0, \eps])} ,\|h\|_{L^\infty\left([0, \eps]\right)} \right\}.\] 
For any $(f,h)\in \cX_\eps,$ we define  
\[\Phi(f,h):=\left(\Phi_1(f,h),\Phi_2(f,h)\right),\]
where 
\begin{equation}\label{Phi-defn}
\left\{\begin{aligned}
&\Phi_1(f,h)(r):=\eta_0+\int_0^r h(\rho)\,d\rho\quad\forall 0\le r\leq\eps,\\
&\Phi_2(f,h)(r):=-\frac{f(r)^{1-m}}{r^{n-1}}\int_0^r\rho^{n-1}(\alpha f(\rho)+\beta \rho h(\rho))\,d\rho\quad\forall 0< r\leq\eps,\quad \Phi_2(f,h)(0)=0.
\end{aligned}\right.
\end{equation}
Let 
\[
\cD_{\eps,\eta_0}:=\left\{ (f,h)\in \cX_\eps:  ||(f,h)-(\eta_0,0)||_{\cX_{\eps}}\leq \eta_0/2\right\}\quad\forall \eps\in (0,1).\]
Then $\cD_{\eps,\eta_0}$ is a closed subset of $\cX_\eps$.
For any $(f,h)\in \cD_{\eps,\eta_0}$,
\begin{equation}\label{h-l1-small}
\max_{0\leq r\leq\eps} \left| \int_0^r h(\rho)\,d\rho\right| \leq\frac{\eta_0\eps}{2}  \leq\eta_0/2
\end{equation} 
and 
\begin{align}\label{f-integral-small}
\frac{f^{1-m}(r)}{r^{n-1}}\int_0^r\rho^{n-1}\left(|\alpha|f(\rho)+|\beta\rho h(\rho)|\right)d\rho
\leq&\left(\frac{3\eta_0}{2}\right)^{1-m}\left(\frac{(3\eta_0/2)|\alpha|\,r}{n}+\frac{ |\beta|(\eta_0/2)r^2}{n+1}\right)\notag\\
\leq&(3\eta_0)^{2-m}(|\alpha |+|\beta|)r\notag\\
= &M_1 r\eta_0\leq M_1\eta_0\eps\le \eta_0/2\quad\forall 0< r\leq\eps
\end{align} 
if $\eps\in \left(0, \frac{1}{4M_1}\right)$ where 
\[
M_1=\max (1,3^{2-m}\eta_0^{1-m}(|\alpha |+|\beta|))
\]
since $\eta_0/2\le f(r)\le 3\eta_0/2$ for any $r\in[0,\eps]$. By \eqref{Phi-defn} and \eqref{f-integral-small} $\Phi_2(f,h)$ is continuous at $r=0$.
We will now fix $\eps\in (0, 1/(4M_1))$. Then by \eqref{Phi-defn}, \eqref{h-l1-small} and \eqref{f-integral-small}, $\Phi(\cD_{\eps,\eta_0})\subset \cD_{\eps,\eta_0}$. We will now prove that the map $\Phi:\cD_{\eps,\eta_0}\to \cD_{\eps,\eta_0}$ is Lipschitz continuous with a Lipschitz constant which is less than $1$. 
Let $(f_1,h_1),(f_2,h_2)\in \cD_{\eps,\eta_0}$ and $\delta:=||(f_1,h_1)-(f_2,h_2)||_{\cX_\eps}$. Then   
\begin{equation*}
\|\Phi_1(f_1,h_1)- \Phi_1(f_2,h_2)\|_{L^\infty([0,\eps])}\leq  \max_{0\leq r\leq\eps }\int_0^r |h_1(\rho)-h_2(\rho)|d\rho\leq \eps\delta\le\delta/4
\end{equation*}
and  for any $0< r\leq \eps$, 
\begin{align*}
&\left| \frac{f_1^{1-m}(r)}{r^{n-1}}\int_0^r\rho^{n-1}\left(\alpha f_1(\rho)+\beta\rho h_1(\rho)\right)d\rho
-\frac{f_2^{1-m}(r)}{r^{n-1}}\int_0^r\rho^{n-1}\left(\alpha f_2(\rho)+\beta\rho h_2(\rho)\right)d\rho 
  \right|\\
\leq&\frac{|f_1^{1-m}(r)-f_2^{1-m}(r)|}{r^{n-1}}\int_0^r\rho^{n-1}\left(|\alpha| f_1(\rho)+|\beta|\rho |h_1(\rho)|\right)d\rho
+\delta\cdot  \frac{f_2^{1-m}(r)}{r^{n-1}}\int_0^r\rho^{n-1}\left(|\alpha| +|\beta|\rho \right)d\rho\\
\leq& (1-m) \frac{|f_1(r)-f_2(r)|}{  \min (f_1(r),f_2(r))^m} \left(\frac{(3\eta_0/2)|\alpha|\,r}{n}+\frac{ |\beta|(\eta_0/2)r^2}{n+1}\right)\notag\\
&\qquad+\delta(3\eta_0/2)^{1-m}(|\alpha|+|\beta| )r   \\
\leq& (1-m) \frac{(3\eta_0/2)\delta}{ (\eta_0/2)^m} \left(|\alpha|r+|\beta|r^2\right)+M_1r\delta\\
\leq&((1-m)M_1+M_1)\eps\delta\quad\forall 0< r\leq\eps\notag\\
\leq&2M_1\eps\delta\notag\\
\le&\delta/2\quad\forall 0< r\leq\eps
\end{align*}
since $\eta_0/2\leq f_1(r), f_2(r)\leq 3\eta_0/2$ for any $r\in[0,\eps]$. Hence 
\begin{align*}
&\|\Phi_2(f_1,h_1)- \Phi_2(f_2,h_2)\|_{L^\infty([0,\eps])}\\
\leq&\max_{0 < r\leq\eps }\left| \frac{f_1^{1-m}(r)}{r^{n-1}}\int_0^r\rho^{n-1}\left(\alpha f_1(\rho)+\beta\rho h_1(\rho)\right)d\rho
-\frac{f_2^{1-m}(r)}{r^{n-1}}\int_0^r\rho^{n-1}\left(\alpha f_2(\rho)+\beta\rho h_2(\rho)\right)d\rho 
  \right|\\
\leq&\delta/2\quad\forall 0< r\leq\eps.
\end{align*}
Thus $\Phi$ is Lipschitz continuous on $\cD_{\eps,\eta_0}$ with Lipschitz constant $1/2$. 
Hence by the contraction map theorem there exists a unique fixed point $(f,h)=\Phi (f,h)$ in $\cD_{\eps,\eta_0}$. Then
\begin{align}\label{eq-f-h-representation}
&\left\{\begin{aligned}
f(r)=&\eta_0+\int_0^rh(\rho)\,d\rho\quad\forall 0\le r\le\eps\\
h(r)=&-\frac{ f^{1-m}(r)}{r^{n-1}}\int_0^r\rho^{n-1}(\alpha f(\rho)+\beta\rho h(\rho))\,d\rho \quad\forall 0<r\le\eps.
\end{aligned}\right.
\end{align}
Hence
\begin{align}
&f_r(r)=h(r)=-\frac{f^{1-m}(r)}{r^{n-1}}\int_0^r\rho^{n-1}(\alpha f(\rho)+\beta\rho f_r(\rho))\,d\rho  \quad\forall 0<r\le\eps\label{eqn-f-integral-representation}\\
\Rightarrow\quad&r^{n-1}(f^m/m)_r(r)=-\int_0^r\rho^{n-1}(\alpha f(\rho)+\beta\rho f_r(\rho))\,d\rho \quad\forall 0<r\le\eps.\label{eqn-f-integral-representation2}
\end{align}
By \eqref{eq-f-h-representation} and \eqref{eqn-f-integral-representation}, $f(0)=\eta_0$ and $f_r(r)$ is   continuously differentiable in $(0,\eps)$. By integrating the drift term in \eqref{eqn-f-integral-representation2} by parts we get \eqref{fr-representation-formula}.  Dividing \eqref{fr-representation-formula} by $r$ and letting  $r\to 0$, by \eqref{f-fr-origin-value} and l'Hospital's rule we get \eqref{origin-second-order-expansion}. 

Since 
$f_r=h\in C([0,\eps])$, by \eqref{eqn-f-integral-representation},
\begin{equation*}
|f_r(r)|\le\frac{C}{r^{n-1}}\int_0^r\rho^{n-1}\,d\rho\le C'r\to 0\quad\mbox{ as }r\to 0
\end{equation*}
and then $f$ belongs to $C^1([0,\eps])\cap C^2((0,\eps])$ and satisfies  \eqref{f-fr-origin-value}. 
Differentiating \eqref{eqn-f-integral-representation2} with respect to 
$r\in (0,\eps)$, we get that $f$ satisfies \eqref{f-ode}  in $(0,\eps)$. 

Suppose $\overline{f}\in C^1([0,\eps])\cap C^2((0,\eps])$ is another solution of \eqref{f-ode}  in $(0,\eps)$ that satisfies \eqref{f-fr-origin-value}. Then by \eqref{f-ode}  and \eqref{f-fr-origin-value},
\begin{equation}\label{overline-f-r-integral-eqn}
\overline{f}_r(r)=-\frac{\overline{f}(r)^{1-m}}{r^{n-1}}\int_0^r\rho^{n-1}(\alpha \overline{f}(\rho)+\beta \rho \overline{f}_r(\rho))\,d\rho\quad\forall 0<r<\eps.
\end{equation}
Let $\overline{h}=\overline{f}_r$. Then by \eqref{f-fr-origin-value} and \eqref{overline-f-r-integral-eqn},
\begin{equation}\label{bar-f-h-eqns}
\left\{\begin{aligned}
&\overline{f}(r):=\eta_0+\int_0^r \overline{h}(\rho)\,d\rho=\Phi_1(\overline{f},\overline{h})(r)\qquad\qquad\qquad\qquad\qquad\forall 0<r<\eps\\
&\overline{h}(r):=-\frac{\overline{f}(r)^{1-m}}{r^{n-1}}\int_0^r\rho^{n-1}(\alpha \overline{f}(\rho)+\beta \rho \overline{h}(\rho))\,d\rho=\Phi_2(\overline{f},\overline{h})(r)\quad\,\forall 0<r<\eps.
\end{aligned}\right.
\end{equation}
Let 
\begin{equation}\label{m2-3-defn}
M_2=\max(1,\|\overline{h}\|_{L^{\infty}(0,\eps)}),\quad M_3=\max\left(1,(3\eta_0)^{1-m}(3\eta_0|\alpha|+M_2|\beta|)\right)
\end{equation}
 and
\begin{equation}\label{epsilon1-defn}
\eps_1=\min \left(\frac{\eps}{2},\frac{\eta_0}{2M_2},\frac{\eta_0}{2M_3}\right).
\end{equation}
Then by \eqref{bar-f-h-eqns}, \eqref{m2-3-defn} and \eqref{epsilon1-defn},
\begin{equation}\label{f-bar-bds}
|\overline{f}(r)-\eta_0|\le M_2r\le\eta_0/2\quad\forall 0<r\le\eps_1
\end{equation}
and by \eqref{bar-f-h-eqns} and \eqref{f-bar-bds}
\begin{align}\label{h-bar-bd}
|\overline{h}(r)|\le&\frac{(3\eta_0/2)^{1-m}}{r^{n-1}}\int_0^r\rho^{n-1}((3\eta_0/2)|\alpha|+M_2|\beta| \rho) \,d\rho\notag\\
\le&M_3r\notag\\
\le&\eta_0/2\quad\forall 0<r\le\eps_1.
\end{align}
By \eqref{f-bar-bds} and \eqref{h-bar-bd}, $(\overline{f},\overline{h})\in\cD_{\eps_1,\eta_0}$. Since $\eps_1\leq\eps$, the invariance and contraction estimates proved above remain valid after restricting $\Phi$ to $\cD_{\eps_1,\eta_0}$; the restriction of $(f,h)$ to $[0,\eps_1]$ is its fixed point. Hence by \eqref{bar-f-h-eqns} and the uniqueness of the fixed point of the contraction map $\Phi:\cD_{\eps_1,\eta_0}\to \cD_{\eps_1,\eta_0}$, 
\begin{equation}\label{f-f-bar-near-0-equal}
(f(r),h(r))=(\overline{f}(r),\overline{h}(r))\quad\forall 0<r\le\eps_1.
\end{equation}
Now both $f$ and $\overline{f}$ are the solution of \eqref{f-ode}  in $[\eps_1,\eps)$ with $f(\eps_1)=\overline{f}(\eps_1)$ and $f_r(\eps_1)=\overline{f}_r(\eps_1)$. Then by standard ODE theory,
\begin{equation}\label{f-f-bar-away-0-equal}
f(r)=\overline{f}(r)\quad\forall \eps_1\le r<\eps.
\end{equation}
By \eqref{f-f-bar-near-0-equal} and \eqref{f-f-bar-away-0-equal}, $f(r)=\overline{f}(r)$ for any $0<r<\eps$.
Hence $f\in C^1([0,\eps])\cap C^2((0,\eps])$ is the unique solution of  \eqref{f-ode}  in $(0,\eps)$  which satisfies  \eqref{f-fr-origin-value}.
By \eqref{f-ode},

\begin{equation}\label{f-rr-eqn}
f_{rr}(r)+(m-1)\frac{f_r(r)^2}{f(r)}+(n-1)\frac{f_r(r)}{r}+\alpha f(r)^{2-m}+\beta rf(r)^{1-m}f_r(r)=0\quad\forall 0<r\le\eps.
\end{equation}
Letting $r\to 0$ in \eqref{f-rr-eqn}, by \eqref{f-fr-origin-value} and \eqref{origin-second-order-expansion} we have
\begin{equation*}
\lim_{r\to 0}f_{rr}(r)=-\alpha \eta_0^{2-m}-(n-1)\lim_{r\to 0}\frac{f_r(r)}{r}=-\alpha \eta_0^{2-m}+\alpha \frac{(n-1)}{n}\eta_0^{2-m}=- \frac{\alpha}{n}\eta_0^{2-m}.
\end{equation*}
Hence we can extend $f_{rr}$ to a continuous function on $[0,\eps]$ by letting 
\begin{equation}\label{f-rr-value-at-origin}
f_{rr}(0)=-\frac{\alpha}{n}\eta_0^{2-m}. 
\end{equation}
We now extend $f$ to a radially symmetric function on $B_{\eps}$ by letting $f(x)=f(|x|)$ for any $x\in B_{\eps}$. Then for any  $x\neq0$,
\begin{equation*}
\frac{\partial^2 f}{\partial x_i\partial x_j}(x)
= \frac{\partial^2 f}{\partial x_i\partial x_j}(|x|)
 =f_{rr}(|x|)\frac{x_i x_j}{|x|^2}
 +\frac{f_r(|x|)}{|x|}
 \left(\delta_{ij}-\frac{x_i x_j}{|x|^2}\right).
\end{equation*}
Hence
\begin{equation}\label{f-ij-eqn}
\left|\frac{\partial^2 f}{\partial x_i\partial x_j}(x)-\frac{f_r(|x|)}{|x|}
\delta_{ij}\right|
\le \left|f_{rr}(|x|)-\frac{f_r(|x|)}{|x|}\right|\left|\frac{x_i x_j}{|x|^2}\right|
\le\left|f_{rr}(|x|)-\frac{f_r(|x|)}{|x|}\right|
\end{equation}
Letting $0\ne |x|\to 0$ in \eqref{f-ij-eqn}, by \eqref{origin-second-order-expansion} and \eqref{f-rr-value-at-origin},
\begin{equation*}
\lim_{0\ne |x|\to 0}\left|\frac{\partial^2 f}{\partial x_i\partial x_j}(x)-\frac{f_r(|x|)}{|x|}
\delta_{ij}\right|=0.
\end{equation*}
Hence
\begin{equation}\label{f-ij-eqn2}
\lim_{0\ne |x|\to 0}\frac{\partial^2 f}{\partial x_i\partial x_j}(x)
=\lim_{0\ne |x|\to 0}\frac{f_r(|x|)}{|x|}\delta_{ij}=-\frac{\alpha}{n}\eta_0^{2-m}\delta_{ij}.
\end{equation}
By \eqref{f-ij-eqn2} we can extend $\frac{\partial^2 f}{\partial x_i\partial x_j}$ to a continuous function at the origin by letting
\begin{equation*}
\frac{\partial^2 f}{\partial x_i\partial x_j}(0)
=-\frac{\alpha}{n}\eta_0^{2-m}\delta_{ij}.
\end{equation*}
Then the function $x\mapsto f(x)$ is $C^2$ at the origin. Hence $f$ satisfies \eqref{f-elliptic-eqn} in $B_{\eps}$. Let $U=f^m/m$.  Since $f(0)=\eta_0>0$, then in a  small neighborhood of the origin $U>0$ and the equation \eqref{f-elliptic-eqn} becomes
\[
 \Delta U+\alpha(mU)^{1/m}
 +\beta(mU)^{(1-m)/m}x\cdot\nabla U=0.
\] 
whose lower-order coefficients are smooth because $U>0$.  Then by standard elliptic bootstrapping argument we get 
$U\in C^{\infty}(B_{\eps})$.  Hence $f\in C^{\infty}(B_{\eps})$.

Note that the construction of solution of \eqref{f-elliptic-eqn} in $B_{\eps}$ satisfying \eqref{f-fr-origin-value} and the limit in \eqref{origin-second-order-expansion} are uniform when $(\alpha,\beta)$ ranges over a compact set. Hence the lemma follows.

\end{proof}

\begin{lem}[Monotonicity while the profile is positive]\label{lem-local-f-derivative-negative}
Let $n\ge 3$, $0<m<\frac{n-2}{n}$, $\eta_0>0$, $\rho_1>0$, $\alpha\not =0$ and $\beta\in\mathbb{R}$. Let $R_o>0$ and $f\in C^1([0,R_o))\cap C^2((0,R_o))$ be a  positive solution of  \eqref{f-ode}, \eqref{f-fr-origin-value}, in $[0,R_o)$. Then $f$ satisfies 
\begin{equation*}
\left\{\begin{aligned}
&f_r(r)<0\quad\forall 0<r<R_o\quad\mbox{ if }\alpha>0\\ 
&f_r(r)>0\quad\forall 0<r<R_o\quad\mbox{ if }\alpha<0.
\end{aligned}
\right.
\end{equation*}
\end{lem}
\begin{proof}
Without loss of generality we may assume that $\alpha>0$. Then by \eqref{f-fr-origin-value}
there exists a constant $R_o'\in (0,R_o)$ such that
\begin{align}\label{f-expression1-ineqn}
&f(r)>\eta_0/2\mbox{ and }|\beta|rf_r(r)<\alpha\eta_0/2\quad\forall 0<r\le R_o'\notag\\
\Rightarrow\quad&\alpha f(r)+\beta rf_r(r)>0\quad\forall 0<r\le R_o'.
\end{align}
Since $f$ satisfies \eqref{eqn-f-integral-representation2}, by  \eqref{f-expression1-ineqn},
\begin{align}\label{f-derivative-negative-near-zero}
&r^{n-1}f(r)^{m-1}f_r(r)=r^{n-1}(f^m/m)_r(r)<0\quad\forall 0<r\le R_o'\notag\\
\Rightarrow\quad&f_r(r)<0\quad\forall 0<r\le R_o'.
\end{align}
Alternately by \eqref{origin-second-order-expansion} there exists a constant $R_o'>0$ such that \eqref{f-derivative-negative-near-zero} holds.
Let 
\begin{equation}
r_1=\sup\{b\in [R_o',R_o):f_r(r)<0\quad\forall 0<r<b\}.
\end{equation} 
 Then $r_1\in (R_o',R_o]$. Suppose $r_1<R_o$. Then
\begin{align}
&f_r(r)<0\quad\forall 0<r<r_1\quad\mbox{ and }\quad f_r(r_1)=0\label{f-fr-value1}\\
\Rightarrow\quad&f_{rr}(r_1)\ge 0.\label{frr-value}
\end{align}
By \eqref{f-ode},
\begin{equation}\label{f-ode2}
f^{m-1}f_{rr}+(m-1)f^{m-2}f_r^2+\frac{n-1}{r}f^{m-1}f_r=-\alpha f-\beta rf_r\quad\forall 0<r<R_o.
\end{equation}
Since $\alpha>0$, by putting $r=r_1$ in \eqref{f-ode2}, by \eqref{f-fr-value1} and \eqref{frr-value} we have
\begin{equation*}
0\le f(r_1)^{m-1}f_{rr}(r_1)=-\alpha f(r_1)<0
\end{equation*}
and a contradiction follows. Hence $r_1=R_o$ and the lemma follows.

\end{proof}

We are now ready to prove Theorem \ref{thm-existence-self-similar1}.

\noindent{\ni{\it Proof of Theorem \ref{thm-existence-self-similar1}:}} 

\noindent  
If $\beta=\beta_-(\rho_1)$, then $\alpha=0$ and the constant function $\eta_0$ is the unique positive solution of \eqref{f-elliptic-eqn} in $\mathbb{R}^n$ which satisfies \eqref{f-fr-origin-value}. Hence we may assume that $\beta_-(\rho_1)<\beta\le\beta_+(\rho_1)$. Then $\alpha>0$.

Let $\eps>0$ and $f^o_{\beta}\in C^1([0,\eps])\cap C^2((0,\eps])$ be a  positive solution of  \eqref{f-ode}, \eqref{f-fr-origin-value}, in $[0,\eps)$ given by Lemma \ref{lem-local-existence}. Let $(0,R_o)$ be the maximal interval of existence of  a positive solution $f^o_{\beta}\in C^1([0,R_o))\cap C^2((0,R_o))$ of  \eqref{f-ode}  satisfying  \eqref{f-fr-origin-value}.  Then $R_o\ge\eps$.
By Lemma \ref{lem-local-f-derivative-negative}, \eqref{f-derivative-negative2} holds.
Hence
\begin{equation}\label{f-bded}
0<f^o_{\beta}(r)<f^o_{\beta}(0)=\eta_0\quad\forall 0<r<R_o.
\end{equation}
Suppose $R_o<\infty$.  Then by \eqref{f-derivative-negative2},
\begin{equation}
\ell=\lim_{r\nearrow R_o}f^o_{\beta}(r)\mbox{ exists }\quad\mbox{ and }\quad \ell\geq0.
\end{equation}
Hence we can extend $f^o_{\beta}$ to a continuous function on $[0,R_o]$ by setting $f^o_{\beta}(R_o)=\ell$.
By the proof of Lemma \ref{lem-local-f-derivative-negative} \eqref{fr-representation-formula0}, \eqref{fr-representation-formula}, holds for any $0<r<R_o$. 
By \eqref{fr-representation-formula} and \eqref{f-bded},
\begin{equation*}
|f^o_{\beta,r}(r)|\le(\alpha+n|\beta|)\frac{\eta_0^{2-m}R_o}{n}+|\beta|\eta_0^{2-m}R_o\le\left(\frac{\alpha}{n}+2|\beta|\right)\eta_0^{2-m}R_o\quad\forall 0<r<R_o
\end{equation*}
and we can extend $f^o_{\beta,r}$ to a continuous function on $[0,R_o]$ by letting 
\begin{equation*}
f^o_{\beta,r}(R_o)=-(\alpha-n\beta)\frac{f^o_{\beta}(R_o)^{1-m}}{R_o^{n-1}}\int_0^{R_o}\rho^{n-1}f^o_{\beta}(\rho)\,d\rho-\beta R_of^o_{\beta}(R_o)^{2-m}.
\end{equation*} 
  If $\ell>0$, then by the regular
ODE continuation theorem we can extend $f^o_{\beta}$ to a solution of   \eqref{f-ode}, \eqref{f-fr-origin-value}, in $[0,R_1]$ for some constant $R_1>R_o$. This contradicts the maximality of $R_o$.  Thus $\ell=0$.  Letting $r\nearrow R_o$ in \eqref{fr-representation-formula0} by Remark \ref{rmk:alpha.n-beta} we get \eqref{origin-endpoint-flux}.  Hence by \eqref{origin-endpoint-flux},
\[
 ((f^o_{\beta})^m)_r=m (f^o_{\beta})^{m-1}f^o_{\beta,r}\longrightarrow \frac{mL_0}{R_o^{n-1}}\quad\mbox{ as }r\nearrow R_o.
\]
Thus \eqref{origin-finite-zero-asymptotic} holds and the theorem follows.

{\hfill$\square$\vspace{6pt}}

\begin{lem}[Local parameter dependence of the origin branch]\label{lem:f-beta-parameter-limit}
Let $n\ge 3$, $0<m<\frac{n-2}{n}$, $\eta_0>0$, $\rho_1>0$ and $\beta_0\in  [\beta_-(\rho_1),\beta_+(\rho_1)]$. Then  there exists a constant $\eps>0$ such that for any $\beta\in  [\beta_-(\rho_1),\beta_+(\rho_1)]$ there exists a  unique positive solution $f^o_{\beta}\in C^1([0,\eps])\cap C^2((0,\eps])$ of \eqref{f-ode} with $\alpha=\frac{2\beta +\rho_1}{1-m}$ in $(0,\eps)$ which satisfies \eqref{f-fr-origin-value} and
\begin{equation}\label{eq:f-beta-parameter-limit}
f^o_{\beta}(r)\to f^o_{\beta_0}(r)\quad\mbox{ and }\quad f^o_{\beta,r}(r)\to f^o_{\beta_0,r}(r)\quad\mbox{ uniformly in $[0,\eps]$}\quad\mbox{ as }\beta\to\beta_0.
\end{equation}
Moreover
\begin{equation}\label{uniform-origin-log-slope}
 \sup_{\beta_-(\rho_1)\le\beta\le\beta_+(\rho_1)}\abs{\frac{rf^o_{\beta,r}(r)}{f^o_\beta(r)}}
 \longrightarrow0
 \qquad\text{as }r\searrow 0.
\end{equation}
\end{lem}
\begin{proof}
By an argument similar to the proof of Lemma \ref{lem-local-existence} for any $\beta\in  [\beta_-(\rho_1),\beta_+(\rho_1)]$ there exists a constant $0<\eps<1$ such that there exists a  unique positive solution $f^o_{\beta}\in C^1([0,\eps])\cap C^2((0,\eps])$ of \eqref{f-ode} with $\alpha=\frac{2\beta +\rho_1}{1-m}$ in $(0,\eps)$ which satisfies \eqref{f-fr-origin-value} and
\begin{align}\label{eq-f-beta-representation}
&\left\{\begin{aligned}
f^o_{\beta}(r)=&\eta_0+\int_0^rf^o_{\beta,r}(\rho)\,d\rho\quad\forall 0\le r\le\eps\\
f^o_{\beta,r}(r)=&-\frac{ (f^o_{\beta})^{1-m}(r)}{r^{n-1}}\int_0^r\rho^{n-1}\left\{\left(\frac{2\beta+\rho_1}{1-m} \right)f^o_{\beta}(\rho)+\beta\rho f^o_{\beta,r}(\rho)\right\}\,d\rho \quad\forall 0<r\le\eps.
\end{aligned}\right.
\end{align}
with
\begin{align}\label{f-beta-local-bd1}
&|f^o_{\beta}(r)-\eta_0|\le\eta_0/2\quad\mbox{ and }\quad|f^o_{\beta,r}(r)|\le\eta_0/2\quad\forall 0\le r\le\eps\notag\\
\Rightarrow\quad&\eta_0/2\le f^o_{\beta}(r)\le3\eta_0/2\quad\mbox{ and }\quad|f^o_{\beta,r}(r)|\le\eta_0/2\quad\forall 0\le r\le\eps.
\end{align}
Note that 
\begin{equation}\label{beta-upper-bd1}
|\beta|\le\max(|\beta_-(\rho_1)|,\beta_+(\rho_1))=:k_1.
\end{equation}
Let
\[
 \mathcal N_\beta:=\max\{\|f^o_\beta-f^o_{\beta_0}\|_{L^\infty([0,\eps])},
 \|f^o_{\beta,r}-f^o_{\beta_0,r}\|_{L^\infty([0,\eps])}\}.
\]
Then by  \eqref{eq-f-beta-representation} for any $\beta,\beta_0\in  [\beta_-(\rho_1),\beta_+(\rho_1)]$ we have
\begin{align}\label{f-beta-beta0-norm}
&|f^o_{\beta}(r)-f^o_{\beta_0}(r)|\le\int_0^r|f^o_{\beta,r}(\rho)-f^o_{\beta_0,r}(\rho)|\,d\rho\le r\|f^o_{\beta,r}-f^o_{\beta_0,r}\|_{L^{\infty}([0,\eps])}\quad\forall 0\le r\le\eps\notag\\
\Rightarrow\quad&\|f^o_{\beta}-f^o_{\beta_0}\|_{L^{\infty}([0,\eps])}\le\eps\|f^o_{\beta,r}-f^o_{\beta_0,r}\|_{L^{\infty}([0,\eps])}.
\end{align}
By  the mean value theorem, \eqref{eq-f-beta-representation}, \eqref{beta-upper-bd1}, \eqref{f-beta-local-bd1} and \eqref{f-beta-beta0-norm},
\begin{align}\label{f-beta-beta0-difference2}
&|f^o_{\beta,r}(r)-f^o_{\beta_0,r}(r)|\notag\\
\le&\frac{ |(f^o_{\beta})^{1-m}(r)-(f^o_{\beta_0})^{1-m}(r)|}{r^{n-1}}\int_0^r\rho^{n-1}\left\{\left(\frac{2|\beta|+\rho_1}{1-m} \right)f^o_{\beta}(\rho)+|\beta|\rho |f^o_{\beta,r}(\rho)|\right\}\,d\rho\notag\\
&\quad + \frac{ (f^o_{\beta_0})^{1-m}(r)}{r^{n-1}}\int_0^r\rho^{n-1}\left\{\left(\frac{2|\beta-\beta_0|}{1-m} \right)f^o_{\beta}(\rho)+|\beta-\beta_0|\rho |f^o_{\beta,r}(\rho)|\right\}\,d\rho \notag\\
&\quad+ \frac{ (f^o_{\beta_0})^{1-m}(r)}{r^{n-1}}\int_0^r\rho^{n-1}\left\{\left(\frac{2|\beta_0|+\rho_1}{1-m} \right)|f^o_{\beta}(\rho)-f^o_{\beta_0}(\rho)|+|\beta_0|\rho |f^o_{\beta,r}(\rho)-f^o_{\beta_0,r}(\rho)|\right\}\,d\rho\notag\\
\le&3(k_1+\rho_1)\frac{ |f^o_{\beta}(r)-f^o_{\beta_0}(r)|}{(\eta_0/2)^mr^{n-1}}\int_0^r\rho^{n-1}((3\eta_0/2)+(\eta_0/2))\,d\rho 
+\frac{4(3\eta_0/2)^{2-m}}{1-m}|\beta-\beta_0|r\notag\\
&\quad +\left(\frac{3k_1+\rho_1}{1-m} \right)(3\eta_0/2)^{1-m}r\mathcal N_\beta\notag\\
\le&C_1r|\beta-\beta_0|+2C_1r\mathcal N_\beta\quad\forall 0<r\le\eps
\end{align}
where
\begin{equation*}
C_1=\max\left(1,\frac{12(k_1+\rho_1+1)(\eta_0^{1-m}+\eta_0^{2-m})}{1-m}\right).
\end{equation*}
We now choose $0<\eps<1/(4C_1)$. 
By \eqref{f-beta-beta0-norm} and \eqref{f-beta-beta0-difference2}, 
\begin{align*}
&\mathcal N_\beta\le C_1\eps|\beta-\beta_0|+2C_1\eps\mathcal N_\beta\le 2C_1\eps|\beta-\beta_0|+\frac{1}{2}\mathcal N_\beta\notag\\
\Rightarrow\quad&\mathcal N_\beta\le 4C_1\eps|\beta-\beta_0|
\end{align*}
and  \eqref{eq:f-beta-parameter-limit}  follows. 

We will now prove \eqref{uniform-origin-log-slope}. If $\beta=\beta_-(\rho_1)$, then $\alpha=0$ and $f^o_{\beta}(r)\equiv\eta_0$ for any $r\ge 0$. Hence 
\begin{equation*}
\abs{\frac{rf^o_{\beta_-(\rho_1),r}(r)}{f^o_{\beta_-(\rho_1)}(r)}}=0\quad\forall r\ge 0.
\end{equation*}
For $\beta\ne\beta_-(\rho_1)$, $\beta>\beta_-(\rho_1)$. Hence $\alpha>0$.
Then by Lemma \ref{lem-local-f-derivative-negative}, \eqref{f-derivative-negative2} holds and hence
\begin{equation}\label{f-bded2}
0<f^o_{\beta}(r)<f^o_{\beta}(0)=\eta_0\quad\forall 0<r\le \eps, \beta_-(\rho_1)<\beta\le\beta_+(\rho_1).
\end{equation}
Since $f^o_{\beta}$ satisfies  \eqref{fr-representation-formula}, by \eqref{f-bded2},
\begin{equation}\label{fr-bded1}
|(f^o_{\beta})_r(r)|\le Cr\quad\forall 0<r\le \eps,\beta_-(\rho_1)<\beta\le\beta_+(\rho_1)
\end{equation}
for some constant $C>0$. By \eqref{eq-f-beta-representation} and \eqref{fr-bded1},
\begin{equation}\label{f-lower-bd5}
f^o_{\beta}(r)\ge\eta_0-Cr^2/2\quad\forall 0<r\le \eps,\beta_-(\rho_1)<\beta\le\beta_+(\rho_1)
\end{equation}
By \eqref{fr-bded1} and \eqref{f-lower-bd5} we get \eqref{uniform-origin-log-slope} for any $\beta_-(\rho_1)<\beta\le\beta_+(\rho_1)$ and the lemma follows.
\end{proof}

\begin{thm}[Parameter dependence of the origin branch]\label{thm:parameter-limit-convergence}
Let $n\ge 3$, $0<m<\frac{n-2}{n}$, $\eta_0>0$, $\rho_1>0$, $\beta_0\in  [\beta_-(\rho_1),\beta_+(\rho_1)]$  and $\alpha_0=\frac{2\beta_0+\rho_1}{1-m}$. Let $R_o=R_o(\beta_0,\eta _0)$
and $f^o_{\beta_0}\in C^1([0,R_o))\cap C^2((0,R_o))$ be the unique maximal positive radial branch of \eqref{f-elliptic-eqn} in $(0,R_o)$ which satisfies \eqref{f-fr-origin-value} given by Theorem \ref{thm-existence-self-similar1}.
Let $0<R<R_o(\beta_0,\eta_0)$. Then there exists a constant $\delta_2>0$ such that for any $\beta\in I_{\delta_2}:=[\beta_0-\delta_2,\beta_0+\delta_2]\cap [\beta_-(\rho_1),\beta_+(\rho_1)]$ we have $R_o(\beta,\eta_0)>R$ and
\begin{equation}\label{eq:f-beta-parameter-limit2}
f^o_{\beta}(r)\to f^o_{\beta_0}(r)\quad\mbox{ and }\quad f^o_{\beta,r}(r)\to f^o_{\beta_0,r}(r)\quad\mbox{ uniformly in $[0,R]$}\quad\mbox{ as }\beta\to\beta_0.
\end{equation}
Moreover \eqref{uniform-origin-log-slope} holds.
\end{thm}
\begin{proof}
Let   $\eps>0$  be given by Lemma  \ref{lem:f-beta-parameter-limit}. Without loss of generality we may assume $R>\eps$. Then  by Lemma  \ref{lem:f-beta-parameter-limit} \eqref{eq:f-beta-parameter-limit} and \eqref{uniform-origin-log-slope} holds.
By \eqref{eq:f-beta-parameter-limit},
\begin{equation}\label{eq:f-beta-parameter-limit3}
f^o_{\beta}(\eps)\to f^o_{\beta_0}(\eps)\quad\mbox{ and }\quad f^o_{\beta,r}(\eps)\to f^o_{\beta_0,r}(\eps)\quad\mbox{ as }\beta\to\beta_0.
\end{equation}
Since $f^o_{\beta_0}$ has a positive minimum on $[\eps,R]$, by \eqref{eq:f-beta-parameter-limit3} and the 
continuous dependence of solutions of  regular first-order system on $[\eps,R]$ with respect to its parameters that there exists a constant $\delta_2>0$ such that for any $\beta\in I_{\delta_2}:=[\beta_0-\delta_2,\beta_0+\delta_2]\cap [\beta_-(\rho_1),\beta_+(\rho_1)]$ we have $R_o(\beta,\eta_0)>R$ 
 and
\begin{equation}\label{eq:f-beta-parameter-limit5}
f^o_{\beta}(r)\to f^o_{\beta_0}(r)\quad\mbox{ and }\quad f^o_{\beta,r}(r)\to f^o_{\beta_0,r}(r)\quad\mbox{ uniformly in $[\eps,R]$}\quad\mbox{ as }\beta\to\beta_0.
\end{equation}
Hence by \eqref{eq:f-beta-parameter-limit} and \eqref{eq:f-beta-parameter-limit5}, we get \eqref{eq:f-beta-parameter-limit2} and the theorem follows.

\end{proof}

\section{Existence of solutions originating from infinity}
\setcounter{equation}{0}
\setcounter{thm}{0}

In this section we will use a modification of the technique of \cite{HuS} to prove Theorem \ref{thm-existence-self-similar2}.  
We first note that as observed by K.M.~Hui and Soojung Kim in \cite{HuS} the function  $f$ is a radially symmetric solution of \eqref{f-elliptic-eqn} in $\mathbb{R}^n\setminus\{0\}$ or equivalently 
\eqref{f-ode} in $(0,\infty)$ if and only if 
\begin{equation}\label{eq-def-g-from-f}
g(r):=r^{-\frac{n-2}{m}}f(r^{-1}),\quad r=|x|>0,
\end{equation}   
satisfies \eqref{eq-inversion} or equivalently
\begin{equation}\label{g-eqn7}
(r^{n-1}(g^m/m)_r)_r=-r^{n+\frac{n-2-nm}{m}-3}\left(\widetilde\alpha g(r)+\widetilde\beta r g_r(r)\right)\quad\forall r=|x|>0 
\end{equation}
with \eqref{q-and-tilde-def-intro} hold.
Moreover \eqref{f-infty-behaviour2} holds if and only if $g$ is continuous at $r=0$ with
\begin{equation}\label{g-value-at-origin}
g(0)=\eta.
\end{equation}
Hence the study of the existence of radially symmetric solution of \eqref{f-elliptic-eqn} satisfying
\eqref{f-infty-behaviour2} is equivalent to the study of the existence of the solution of \eqref{eq-inversion}, \eqref{g-value-at-origin}. Note that for $\alpha=\frac{2\beta+\rho_1}{1-m}$ and $\rho_1>0$,
\begin{equation}\label{eq-tilde-alpha-beta-1}
\widetilde\alpha=\alpha-\frac{n-2}{m}\,\beta>0\quad\Leftrightarrow\quad\beta<\frac{m\rho_1}{n-2-nm}=\beta_+(\rho_1).
\end{equation}

We now recall some results of \cite{HuS}.

\begin{lem}[Local inverted branch, Lemma 2.2 and Lemma 2.3 of \cite{HuS}]\label{lem-local-existence-inversion-m-small}
Let $n\geq3$,  $\widetilde \alpha,\widetilde\beta\in\mathbb{R}$ and $\eta>0$. Then the following holds.
\begin{enumerate}[(i)]

\item 
When $0<m<\frac{n-2}{n+1}$, there exists a constant $\eps>0$ such that \eqref{eq-inversion} has a unique positive solution $g\in C^1([0,\eps])\cap C^2((0,\eps])$ in $(0,\eps)$ which satisfies \eqref{g-value-at-origin} and
\begin{equation}\label{eq-initial-m-small}
g_r(0)=0.
\end{equation} 

\item 
When $\frac{n-2}{n+1}\leq m< \frac{n-2}{n}$, there exists a constant $\eps>0$ such that \eqref{eq-inversion} has a unique positive solution $g\in C^{0, \delta_0}([0,\eps])\cap C^2((0,\eps])$ in $(0,\eps)$ which satisfies \eqref{g-value-at-origin} and \eqref{eq-initial-m-large}
where $\delta_0$, $\delta_1$ and $\zeta$ are given by 
\eqref{defn-delta0-1} and
\eqref{zeta-defn}.

\end{enumerate}

\end{lem}

\begin{lem}[Automatic endpoint condition at infinity]\label{lem:automatic-inverted-endpoint}
Let $n\geq3,$  $0<m< \frac{n-2}{n}$, $\eta>0$, $R>0$.
Let $g\in C^2((0,R))$ be a positive solution of
\eqref{eq-inversion} which satisfies 
\begin{equation}\label{g-0-near-zero}
g(r)\to\eta\quad\mbox{ as }r\searrow 0.
\end{equation}  
Then
\begin{equation}\label{automatic-g-flux-identity}
 \begin{aligned}
 r^{n-1}g^{m-1}(r)g_r(r)
 ={}&-\widetilde\beta r^{n+\frac{n-2-nm}{m}-2}g(r)\\
 &+\left(\widetilde\beta\left(\frac{n-2-2m}{m}\right)-\widetilde\alpha\right)
   \int_0^r\rho^{n+\frac{n-2-nm}{m}-3}g(\rho)\,\dd\rho\quad\forall 0<r<R,
 \end{aligned}
\end{equation}
\begin{equation}\label{gr-bd-near-origin}
|g_r(r)|\le CM^{2-m}r^{\frac{n-2-nm}{m}-1}\quad\forall 0<r<R/2
\end{equation}
where $M:=\max_{0\le r\le R/2}g(r)<\infty$, $C>0$ is some constant   and 
 \begin{equation}\label{automatic-g-derivative-bound}
\lim_{r\to 0^+}rg_r(r)=0
\end{equation} 
holds.

\end{lem}

\begin{proof}
We first claim that  there exists a sequence $\{\xi_i\}_{i=1}^{\infty}\subset (0,R)$, $\xi_j\to 0$ as $j\to\infty$,  such that
\begin{equation}\label{rgr-sequence-zero}
\xi_i g_r(\xi_i)\to 0\quad\mbox { as }i\to\infty.
\end{equation}
To prove the claim we let $\{r_j\}_{j=1}^{\infty}\subset (0,R/2)$ be a sequence such that $r_j\to 0$ as $j\to\infty$.
Then by \eqref{g-0-near-zero} and the mean value theorem for any $j\in\mathbb{Z}^+$ there exists a constant $\xi_j\in (r_j, 2r_j)$ such that
\begin{align*}
&g(2r_j)-g(r_j)=r_jg_r(\xi_j)\notag\\
\Rightarrow\quad&|\xi_jg_r(\xi_j)|\le 2|r_jg_r(\xi_j)|\le2|g(2r_j)-g(r_j)|\to 0\quad\mbox{ as }j\to\infty
\end{align*}
and the claim \eqref{rgr-sequence-zero} follows. Integrating \eqref{g-eqn7} over $(\xi_i,r)$, 
\begin{align}
r^{n-1}g(r)^{m-1}g_r(r)=&\xi_i^{n-1}(g^m/m)_r(\xi_i)-\int_{\xi_i}^r\rho^{n+\frac{n-2-nm}{m}-3}\left(\widetilde\alpha g(\rho)+\widetilde\beta \rho g_r(\rho)\right)\,d\rho\quad\forall \xi_i<r<R, i\in\mathbb{Z}^+\label{g-integral-eqn0}\\
=&\xi_i^{n-1}g(\xi_i)^{m-1}g_r(\xi_i)-\widetilde\beta r^{n+\frac{n-2-nm}{m}-2}g(r)+\widetilde \beta \xi_i^{n+\frac{n-2-nm}{m}-2}g(\xi_i)\notag\\
&\qquad +C_2\int_{\xi_i}^r\rho^{n+\frac{n-2-nm}{m}-3}g(\rho)\,d\rho\quad\forall \xi_i<r<R, i\in\mathbb{Z}^+\label{g-integral-eqn2}
\end{align}
where
\begin{equation}\label{c2-defn}
C_2=\widetilde\beta\left(\frac{n-2-2m}{m}\right)-\widetilde\alpha.
\end{equation}
Letting $i\to\infty$ in \eqref{g-integral-eqn2}, by \eqref{rgr-sequence-zero} we get
\eqref{automatic-g-flux-identity}.
By \eqref{g-0-near-zero} and \eqref{automatic-g-flux-identity},
\begin{align}
|rg_r(r)|=&g^{1-m}(r)\left(|\widetilde\beta| r^{\frac{n-2-nm}{m}}g(r)+\frac{|C_2|}{r^{n-2}}
   \int_0^r\rho^{n+\frac{n-2-nm}{m}-3}g(\rho)\,\dd\rho\right)\quad\forall 0<r<R\label{rgr-bd0}\\
   \le& CM^{2-m}r^{\frac{n-2-nm}{m}}\quad\forall 0<r<R/2\label{rgr-bd}
 \end{align}
where \begin{equation*}
M:=\max_{0\le r\le R/2}g(r)<\infty
\end{equation*} 
and  $C>0$ is some constant. Hence \eqref{gr-bd-near-origin} follows. Letting $r\to 0$ in \eqref{rgr-bd} we get \eqref{automatic-g-derivative-bound}
and the lemma follows.
\end{proof}

\begin{lem}[Monotonicity of the inverted branch]\label{g-derivative-sign-lem}
Let $n\geq3,$  $0<m< \frac{n-2}{n}$, $\eta>0$, $R>0$ and $\widetilde\alpha\ne 0$, $\widetilde\beta\in\mathbb{R}.$ Let $g\in C^2((0,R))$ be a positive solution  of \eqref{eq-inversion} in $(0,R)$ which satisfies \eqref{g-0-near-zero}. Then 
\begin{equation}\label{g-integral-eqn}
r^{n-1}(g^m/m)_r(r)=-\int_0^r\rho^{n+\frac{n-2-nm}{m}-3}\left(\widetilde\alpha g(\rho)+\widetilde\beta \rho g_r(\rho)\right)\,d\rho\quad\forall 0<r<R
\end{equation}
holds. Moreover the following holds.
\begin{enumerate}[(i)]

\item If $\widetilde\alpha>0$, then 
\begin{equation}\label{g-derivative-negative}
g_r(r)<0\quad\forall 0<r<R
\end{equation}
and

\item If
$\widetilde\alpha<0$, then 
\begin{equation}\label{g-derivative-positive}
g_r(r)>0\quad\forall 0<r<R.
\end{equation}
\end{enumerate}

\end{lem}

\begin{proof}
When $\widetilde{\alpha}>0$, $\widetilde\beta\not=0$, $\widetilde\alpha/\widetilde\beta\leq \frac{n-2}{m}$ and \eqref{automatic-g-derivative-bound} holds, \eqref{g-derivative-negative} is proved in Lemma 2.1 of \cite{HuS}. We will give a different simple proof for the general case $\widetilde\alpha\ne 0$, $\widetilde\beta\in\mathbb{R}$, here. We first observe that by Lemma \ref{lem:automatic-inverted-endpoint} \eqref{automatic-g-derivative-bound} holds.
Since $\widetilde\alpha\ne 0$, we may assume without loss of generality that $\widetilde\alpha>0$.
By  \eqref{g-0-near-zero} and \eqref{automatic-g-derivative-bound} there exists a constant $R_o'\in (0,R)$ such that
\begin{align}\label{g-gr-ineqn1}
&|\widetilde{\beta}rg_r(r)|<\widetilde{\alpha}\eta/2\quad\mbox{ and }g(r)>\eta/2\quad\forall 0<r\le R_o'\notag\\
\Rightarrow\quad&\widetilde{\alpha}g(r)+\widetilde{\beta}rg_r(r)>0\qquad\qquad\qquad\quad\forall 0<r\le R_o'.
\end{align} 
Integrating \eqref{g-eqn7} over $(0,r)$, $0<r<R$, by   \eqref{automatic-g-derivative-bound} and \eqref{g-gr-ineqn1}, we have
\begin{align}
&r^{n-1}(g^m/m)_r(r)=-\int_0^r\rho^{n+\frac{n-2-nm}{m}-3}\left(\widetilde\alpha g(\rho)+\widetilde\beta \rho g_r(\rho)\right)\,d\rho\quad\forall 0<r<R\notag\\
\Rightarrow\quad &r^{n-1}g(r)^{m-1}g_r(r)<0\quad\forall 0<r\le R_o'\notag\\
\Rightarrow\quad &g_r(r)<0\quad\forall 0<r\le R_o'.\notag
\end{align}
and \eqref{g-integral-eqn} follows.
Let $r_1=\sup\{0<b<R:g_r(r)<0\quad\forall 0<r<b\}$. Then $r_1\in (R_o',R]$. Suppose $r_1<R$. Then
\begin{align}
&g_r(r)<0\quad\forall 0<r<r_1\quad\mbox{ and }\quad g_r(r_1)=0\label{g-ineqn5}\\
\Rightarrow\quad &g_{rr}(r_1)\ge 0.\label{g-ineqn6}
\end{align}
Now by \eqref{eq-inversion},
\begin{equation}\label{g-eqn2}
g^{m-1}g_{rr}+(m-1)g^{m-2}g_r^2+\frac{n-1}{r}g^{m-1}g_r=-r^{\frac{n-2-nm}{m}-2}\left(\widetilde\alpha g+\widetilde\beta rg_r\right)\quad\forall 0<r<R.
\end{equation}
Putting $r=r_1$ in \eqref{g-eqn2}, by \eqref{g-ineqn5} and \eqref{g-ineqn6} we get
\begin{equation*}
0\le g(r_1)^{m-1}g_{rr}(r_1)=-\widetilde\alpha r_1^{\frac{n-2-nm}{m}-2} g(r_1)<0
\end{equation*}
and a contradiction follows. Hence $r_1=R$ and \eqref{g-derivative-negative} follows.
By a similar argument when $\widetilde\alpha<0$, we get \eqref{g-derivative-positive} and the lemma follows.
\end{proof}

\begin{lem}[Uniform endpoint condition]\label{lem:uniform-inverted-endpoint2}
Let $n\geq3$, $0<m<\frac{n-2}{n}$, $\eta>0$, and $R>0$.  Let $K\subset(0,\infty)\times\mathbb{R}$ be compact.  For every $(\widetilde\alpha,\widetilde\beta)\in K$, let $g_{\widetilde\alpha,\widetilde\beta}\in C^2((0,R))$ be a positive solution of \eqref{eq-inversion} satisfying
\[
 g_{\widetilde\alpha,\widetilde\beta}(r)\longrightarrow\eta
 \qquad\text{as }r\searrow0.
\]
Then there is a constant $C_K>0$, independent of $(\widetilde\alpha,\widetilde\beta)\in K$, such that
\begin{equation}\label{rgr-bd6}
 \left|r(g_{\widetilde\alpha,\widetilde\beta})_r(r)\right|
 \leq C_K\eta^{2-m}r^{\frac{n-2-nm}{m}}
 \qquad\forall 0<r<R/2.
\end{equation}
In particular, $r(g_{\widetilde\alpha,\widetilde\beta})_r(r)\to0$ as $r\searrow0$, uniformly for $(\widetilde\alpha,\widetilde\beta)\in K$.
\end{lem}
\begin{proof}
Since $\widetilde\alpha>0$ on $K$, by Lemma \ref{g-derivative-sign-lem},
\begin{equation}\label{g-upper-bd6}
0<g_{\widetilde\alpha,\widetilde\beta}(r)<\eta\quad\forall 0<r<R.
\end{equation}
Then  by the proof of \eqref{gr-bd-near-origin} and \eqref{g-upper-bd6}, we get \eqref{rgr-bd6}
for some constant $C_K>0$. Letting $r\searrow 0$ in \eqref{rgr-bd6} we get $r(g_{\widetilde\alpha,\widetilde\beta})_r(r)\to0$ as $r\searrow0$, uniformly for $(\widetilde\alpha,\widetilde\beta)\in K$ and the lemma follows.

\end{proof}

The following theorem is an extension of Theorem 2.4 of \cite{HuS}.

\begin{thm}[Maximal inverted branch] \label{thm-existence-inversion}
Let $n\geq3,$  $0< m< \frac{n-2}{n}$, $\eta>0$, $\widetilde\alpha\ge 0$ and $\widetilde\beta\in\mathbb{R}.$
\begin{enumerate}[(a)] 
\item If $\widetilde\alpha>0$ and $0<m<\frac{n-2}{n+1}$, then there are a unique constant $R_g\in(0,\infty]$ and a unique maximal positive solution $g\in C^1([0,R_g))\cap C^2((0,R_g))$ of \eqref{eq-inversion} satisfying \eqref{g-value-at-origin} and \eqref{eq-initial-m-small}.
\item If $\widetilde\alpha>0$ and $\frac{n-2}{n+1}\leq m< \frac{n-2}{n}$, then there are a unique constant $R_g\in(0,\infty]$ and a unique maximal positive solution $g\in C_{loc}^{0, \delta_0}([0,R_g))\cap C^2((0,R_g))$ satisfying \eqref{g-value-at-origin} and \eqref{eq-initial-m-large}, where $\delta_0$ and $\delta_1$ are given by \eqref{defn-delta0-1}.
\item If $\widetilde\alpha=0$, then $R_g=\infty$ and $g\equiv\eta$ is the unique maximal positive solution in the endpoint class of Lemma \ref{lem-local-existence-inversion-m-small}: it satisfies \eqref{g-value-at-origin} and \eqref{eq-initial-m-small} when $0<m<\frac{n-2}{n+1}$, and it satisfies \eqref{g-value-at-origin} and \eqref{eq-initial-m-large} with $\zeta=0$ when $\frac{n-2}{n+1}\le m<\frac{n-2}{n}$. 
\end{enumerate}
Moreover \eqref{automatic-g-derivative-bound} holds and, if $\widetilde\alpha>0$, then
\begin{equation}\label{g-derivative-negative2}
g_r(r)<0\quad\forall 0<r<R_g.
\end{equation}
If $R_g<\infty$, then $g(r)\searrow0$ as $r\nearrow R_g$.

\end{thm}
\begin{proof}
If $\widetilde\alpha=0$, the constant function $g\equiv\eta$ satisfies the appropriate endpoint condition.  Its local uniqueness follows from Lemma \ref{lem-local-existence-inversion-m-small}, and regular ODE uniqueness away from $r=0$ extends this identity to all $r>0$.  Thus part (c) holds, and we may now assume $\widetilde\alpha>0$.
By Lemma \ref{lem-local-existence-inversion-m-small} and a continuation argument either (a) or (b) of 
Theorem \ref{thm-existence-inversion} holds. By Lemma \ref{lem:automatic-inverted-endpoint} \eqref{automatic-g-derivative-bound} holds and if $\widetilde\alpha>0$, then by Lemma \ref{g-derivative-sign-lem} \eqref{g-derivative-negative2} holds.

By \eqref{g-derivative-negative2},
\begin{equation}\label{g-bded}
0<g(r)<g(0)=\eta\quad\forall 0<r<R_g.
\end{equation}
Suppose $R_g<\infty$.   Then by \eqref{g-derivative-negative2},
\begin{equation}
\ell=\lim_{r\nearrow R_g}g(r)\mbox{ exists }\quad\mbox{ and }\quad \ell\geq0.
\end{equation}
Hence  we can extend $g$ to a continuous function on $[0,R_g]$ by letting $g(R_g)=\ell$.
By Lemma \ref{lem:automatic-inverted-endpoint}, \eqref{automatic-g-flux-identity}
holds for any $0<r<R_g$. 
By \eqref{automatic-g-flux-identity} and \eqref{g-bded},
 \begin{equation*}
|g_r(r)|
\le C\eta^{2-m}\max_{R_g/2\le r\le R_g}r^{\frac{n-2-nm}{m}-1}
\end{equation*}
for some constant $C>0$.
Hence we can extend $g_r$ to a continuous function on $[0,R_g]$ by letting 
\begin{equation*}
g_r(R_g)
=\left(-\widetilde\beta R_g^{\frac{n-2-nm}{m}-1}\ell+C_2R_g^{1-n}
\int_0^{R_g}\rho^{n+\frac{n-2-nm}{m}-3}g(\rho)\,\dd\rho\right)\ell^{1-m}
\end{equation*}
  If $\ell>0$, then by the regular
ODE continuation theorem we can extend $g$ to a solution of   \eqref{eq-inversion} in $(0,R_1)$ for some $R_1>R_g$ which contradicts the maximality of $R_g$. Hence $\ell=0$. Thus $g(r)\searrow 0$ as $r\nearrow R_g$. 

\end{proof}

\begin{lem}\label{lem-existence-self-similar2}
Let $n\ge 3$, $0<m<\frac{n-2}{n}$, $\eta>0$, $\rho_1>0$, $\alpha=\frac{2\beta+\rho_1}{1-m}$, and $\beta\leq\beta_+(\rho_1)$.  Suppose $f\in C^2((R_Q,\infty))$ is a positive radial solution of \eqref{f-ode} in $(R_Q,\infty)$ for some $R_Q\in[0,\infty)$ and satisfies \eqref{f-infty-behaviour2}.  With the convention $0^{-1}=\infty$, let $R_g=R_Q^{-1}$ and define $g$ by \eqref{eq-def-g-from-f} for $0<r<R_g$, with $g(0)=\eta$.  Let $\widetilde\alpha$, $\widetilde\beta$, $\delta_0$, and $\delta_1$ be given by \eqref{q-and-tilde-def-intro} and \eqref{defn-delta0-1}.  Then $g$ satisfies \eqref{eq-inversion} in $(0,R_g)$ and \eqref{g-integral-eqn}.  Moreover,
\begin{equation}\label{eq-g-class}
\left\{\begin{aligned}
&g\in C^1([0,R_g))\cap C^2((0,R_g))&&\text{if }0<m<\frac{n-2}{n+1},\\
&g\in C_{\mathrm{loc}}^{0,\delta_0}([0,R_g))\cap C^2((0,R_g))&&\text{if }\frac{n-2}{n+1}\leq m<\frac{n-2}{n},
\end{aligned}\right.
\end{equation}
and
\begin{equation}\label{eq-g(0)-g'(0)}
\left\{\begin{aligned}
&g(0)=\eta,\qquad g_r(0)=0
&&\text{if }0<m<\frac{n-2}{n+1},\\
&g(0)=\eta,\qquad
 \lim_{r\to0^+}r^{\delta_1}g_r(r)
 =-\frac{\widetilde\alpha m\eta^{2-m}}{n-2-2m}
&&\text{if }\frac{n-2}{n+1}\leq m<\frac{n-2}{n}.
\end{aligned}\right.
\end{equation}
\end{lem}
\begin{proof}
We will use a modification of the proof of Lemma 3.3 of \cite{HuS} to prove the lemma.
By  \eqref{f-infty-behaviour2}, \eqref{f-ode} and \eqref{eq-def-g-from-f}, $g$ satisfies \eqref{eq-inversion} in $(0,R_g)$ and \eqref{g-0-near-zero} holds. Hence we can extend $g$ to a continuous function on $[0,R_g)$ by letting $g(0)=\eta$. Hence $g\in C([0,R_g))\cap C^2((0,R_g))$.
By Lemma \ref{lem:automatic-inverted-endpoint} and
 Lemma \ref{g-derivative-sign-lem}  \eqref{automatic-g-flux-identity}, \eqref{gr-bd-near-origin} and  \eqref{automatic-g-derivative-bound} hold.
Let $C_2$ be given by \eqref{c2-defn}.

We now divide the proof of \eqref{eq-g-class} and \eqref{eq-g(0)-g'(0)}  into two cases.

\noindent $\underline{\text{\bf Case 1}}$: $0<m<\frac{n-2}{n+1}$

\noindent By \eqref{automatic-g-flux-identity} and l'Hospital's rule, 
\begin{align}\label{g-derivative=0-1}
\lim_{r\to 0}g_r(r)=&\lim_{r\to 0}g(r)^{1-m}\cdot\left(-\widetilde\beta \lim_{r\to 0} r^{\frac{n-2-(n+1)m}{m}}g(r) +\lim_{r\to 0}\frac{C_2}{r^{n-1}}\int_0^r\rho^{n+\frac{n-2-nm}{m}-3}g(\rho)\,d\rho\right)\notag\\
=&\eta^{1-m}\left(-\widetilde\beta\eta\lim_{r\to 0} r^{\frac{n-2-(n+1)m}{m}} +C_2\lim_{r\to 0}\frac{r^{n+\frac{n-2-nm}{m}-3}g(r)}{(n-1)r^{n-2}}\right)\notag\\
=&\left(-\widetilde\beta+\frac{C_2}{n-1}\right)\eta^{2-m}\lim_{r\to 0}r^{\frac{n-2-(n+1)m}{m}}\notag\\
=&0.
\end{align}
Then by \eqref{g-derivative=0-1} $g_r(r)$ can be extended to a continuous function on $[0,R_g)$ by setting $g_r(0)=0$.

\noindent $\underline{\text{\bf Case 2}}$: $\frac{n-2}{n+1}\le m<\frac{n-2}{n}$

\noindent By \eqref{c2-defn}, \eqref{automatic-g-flux-identity} and l'Hospital's rule, 
\begin{align}\label{g-derivative=0-2}
\lim_{r\to 0}r^{\delta_1}g_r(r)=&\lim_{r\to 0}g(r)^{1-m}\cdot\left(-\widetilde\beta \lim_{r\to 0}g(r) +\lim_{r\to 0}\frac{C_2}{r^{n-2+\frac{n-2-nm}{m}}}\int_0^r\rho^{n+\frac{n-2-nm}{m}-3}g(\rho)\,d\rho\right)\notag\\
=&\eta^{1-m}\left(-\widetilde\beta\eta+C_2\lim_{r\to 0}\frac{r^{n+\frac{n-2-nm}{m}-3}g(r)}{(n-2+\frac{n-2-nm}{m})r^{n-3+\frac{n-2-nm}{m}}}\right)\notag\\
=&\eta^{2-m}\left(-\widetilde\beta+\frac{mC_2}{n-2-2m}\right)\notag\\
=&-\frac{\widetilde\alpha m\eta^{2-m}}{n-2-2m}.
\end{align}
Let $\widetilde h(r):=r^{\delta_1}g_r(r)$ for any $0<r<R_g$. Then by \eqref{g-derivative=0-2} $\widetilde h$ can be extended to a continuous function on $[0,R_g)$ by setting 
\begin{equation*}
\widetilde h(0)=-\frac{\widetilde\alpha m\eta^{2-m}}{n-2-2m}.
\end{equation*}
Moreover by \eqref{gr-bd-near-origin} with $R=R_g/2$ and the mean value theorem,
\[
 |g(r)-g(s)|\leq C|r^\frac{n-2-nm}{m}-s^\frac{n-2-nm}{m}|\leq C'|r-s|^\frac{n-2-nm}{m}\leq C'|r-s|^{\delta_0}
 \qquad\forall 0\leq s,r\le R_g/2
\]
for some constants $C>0$, $C'>0$, $C''>0$ since $\delta_0=q/2$. Hence $g\in C_{\mathrm{loc}}^{0,\delta_0}([0,R_g))\cap C^2((0,R_g))$.
Hence  $g$ satisfies  \eqref{g-integral-eqn}, \eqref{eq-g-class}  and \eqref{eq-g(0)-g'(0)}  and the lemma follows.

\end{proof}

The following lemma is an extension of Lemma 3.1, Lemma 3.2 and Corollary 3.4 of \cite{HuS}.

\begin{lem}\label{lem-existence-self-similar}
Let $n\ge 3$, $0<m<\frac{n-2}{n}$, $\eta>0$, $\rho_1>0$,  $\alpha=\frac{2\beta+\rho_1}{1-m}$ and $\beta\le\beta_+(\rho_1)$.  
Then there are a unique $R_Q\in[0,\infty)$ and a unique maximal positive radial solution $f$ of \eqref{f-ode} in $(R_Q,\infty)$ which satisfies \eqref{f-infty-behaviour2}. Moreover if $\beta<\beta_+(\rho_1)$, then
\begin{equation}\label{f-ineqn1}
f(r)+\frac{m}{n-2}rf_r(r)>0\quad\forall r>R_Q,
\end{equation}
\begin{equation}\label{f-upper-bd1}
f(r)<\eta r^{-\frac{n-2}{m}}\quad\forall r>R_Q,
\end{equation}
and
\begin{equation}\label{fr-infty-limit}
\lim_{r\to\infty}r^{\frac{n-2}{m}+1}f_r(r)=-\frac{n-2}{m}\eta.
\end{equation}
If $0<\beta<\beta_+(\rho_1)$, then
\begin{equation}\label{f-ineqn2}
\alpha f(r)+\beta rf_r(r)>0\quad\forall r>R_Q.
\end{equation}
At the endpoint $\beta=\beta_+(\rho_1)$, one has
\begin{equation}\label{f-endpoint-betaplus-equality}
 f(r)=\eta r^{-\frac{n-2}{m}},\qquad
 \alpha f(r)+\beta r f_r(r)=0\quad\forall r>0.
\end{equation}
\end{lem}
\begin{proof}
Let $\widetilde\alpha$ and $\widetilde\beta$ be given by \eqref{q-and-tilde-def-intro}.
Since $\beta\le\beta_+(\rho_1)$, $\widetilde{\alpha}\ge 0$ and if $\beta<\beta_+(\rho_1)$, then $\widetilde\alpha>0$.
By Theorem \ref{thm-existence-inversion}, there are a unique constant $R_g\in(0,\infty]$ and a unique maximal positive solution $g$ of \eqref{eq-inversion} in $[0,R_g)$ which satisfies \eqref{eq-g-class}, \eqref{eq-g(0)-g'(0)},
and   \eqref{automatic-g-derivative-bound},
where $\delta_0$ and $\delta_1$ are given by \eqref{defn-delta0-1}.  Let $R_Q=0$ if $R_g=\infty$ and $R_Q=R_g^{-1}$ otherwise, and let
\begin{equation}\label{f-defn}
f(r):=r^{-\frac{n-2}{m}}g(r^{-1}),\quad r>R_Q.
\end{equation}
By \eqref{eq-inversion}, \eqref{g-value-at-origin} and \eqref{f-defn}, $f$ satisfies \eqref{f-ode} in $(R_Q,\infty)$ and \eqref{f-infty-behaviour2} holds.
To prove uniqueness, let $\bar f$ be any other positive exterior branch
with the same limit \eqref{f-infty-behaviour2}, and let $\bar g$ be its
inversion.  If $\beta<\beta_+(\rho_1)$, by Lemma
\ref{lem-existence-self-similar2} $\bar g$ is in the same endpoint
class as $g$. When $\beta=\beta_+(\rho_1)$, $\bar g\equiv g\equiv\eta$.  Then local uniqueness in Lemma
\ref{lem-local-existence-inversion-m-small} gives $\bar g=g$ near zero,
and regular ODE uniqueness continues the equality throughout their
common positive interval.  Maximality then gives the same endpoint and
hence $\bar f=f$.  Thus $f$ is the unique maximal positive exterior
branch satisfying \eqref{f-infty-behaviour2}.
We now let $\beta<\beta_+(\rho_1)$. 
By  \eqref{g-derivative-negative2} and \eqref{f-defn} we get \eqref{f-ineqn1}. 
By \eqref{g-derivative-negative2}, \eqref{eq-g(0)-g'(0)}  and \eqref{f-defn}, 
\begin{align*}
&g(r)<g(0)=\eta\quad\forall 0<r<R_g\\
\Rightarrow\quad&r^{\frac{n-2}{m}}f(r)<\eta\quad\forall r>R_Q
\end{align*}
and \eqref{f-upper-bd1} follows.
By \eqref{f-infty-behaviour2}, \eqref{automatic-g-derivative-bound} and
\eqref{f-defn} we get \eqref{fr-infty-limit}.

Finally suppose $0<\beta<\beta_+(\rho_1)$.  Then
$\widetilde\alpha>0$, $\widetilde\beta=-\beta<0$ and
$g_r<0$ by \eqref{g-derivative-negative2}.  Hence
\[
 \widetilde\alpha g(r)+\widetilde\beta r g_r(r)>0
 \quad\forall 0<r<R_g.
\]
By direct differentiation of \eqref{f-defn} we get
\[ \alpha f(r)+\beta r f_r(r)
 =r^{-\frac{n-2}{m}}
  \bigl(\widetilde\alpha g(s)+\widetilde\beta s g_r(s)\bigr)>0\quad\forall 0<s=r^{-1}<R_g
\]
and \eqref{f-ineqn2} follows.  When
$\beta=\beta_+(\rho_1)$, Theorem \ref{thm-existence-inversion} gives
$g\equiv\eta$, and \eqref{f-endpoint-betaplus-equality} follows.  This proves the lemma.

\end{proof}

By Lemma \ref{lem-existence-self-similar2}, Lemma \ref{lem-existence-self-similar},  we get Theorem \ref{thm-existence-self-similar2}.

\begin{lem}\label{lem:g-beta-parameter-limit}
Let $n\ge 3$, $0<m<\frac{n-2}{n+1}$, $\eta>0$, $\rho_1>0$ and $\beta_0\in  [\beta_-(\rho_1),\beta_+(\rho_1)]$. Then there exists a constant $\eps>0$ such that for any $\beta\in  [\beta_-(\rho_1),\beta_+(\rho_1)]$  \eqref{eq-inversion} has a unique positive solution $g_{\beta}\in C^1([0,\eps])\cap C^2((0,\eps])$ in $(0,\eps)$ with $\alpha=\frac{2\beta +\rho_1}{1-m}$ and $\widetilde\alpha$, $\widetilde\beta$,  given by \eqref{q-and-tilde-def-intro}, which satisfies \eqref{g-value-at-origin} and \eqref{eq-initial-m-small} and 
\begin{equation}\label{eq:g-beta-parameter-limit1}
g_{\beta}(r)\to g_{\beta_0}(r)\quad\mbox{ and }\quad g_{\beta,r}(r)\to g_{\beta_0,r}(r)\quad\mbox{ uniformly in $[0,\eps]$}\quad\mbox{ as }\beta\to\beta_0.
\end{equation}
Moreover, 
\begin{equation}\label{uniform-infinity-log-slope}
 \sup_{\beta_-(\rho_1)\le\beta\le\beta_+(\rho_1)}\abs{\frac{r g_{\beta,r}(r)}{g_\beta(r)}}
 \longrightarrow0
 \qquad\text{as }r\searrow0.
\end{equation}
\end{lem}

\begin{proof}
By an argument similar to the proof of Lemma 2.2 of \cite{HuS} there exists a constant $0<\eps<1$ such that for any $\beta\in  [\beta_-(\rho_1),\beta_+(\rho_1)]$ \eqref{eq-inversion} has a unique solution $g_{\beta}\in C^1([0,\eps])\cap C^2((0,\eps])$ in $(0,\eps)$ which satisfies \eqref{g-value-at-origin},
\eqref{eq-initial-m-small}, 
\begin{equation}\label{eq-g-beta-representation}
\left\{
\begin{aligned}
&g_{\beta}(r):=\eta+\int_0^r g_{\beta,r}(\rho)\,d\rho\quad\forall 0\le r\le\eps\\
&g_{\beta,r}(r):=-\frac{g_{\beta}^{1-m}(r)}{r^{n-1}}\int_0^r \rho^{ n-3+ \frac{n-2-nm}{m}}\left\{\widetilde\alpha g_{\beta}(\rho)+\widetilde\beta \rho g_{\beta,r}(\rho)\right\}\,d\rho\quad\forall 0<r\le\eps
\end{aligned}\right.
\end{equation}
and 
\begin{align}\label{g-beta-local-bd1}
&|g_{\beta}(r)-\eta|\le\eta/2\quad\mbox{ and }\quad|g_{\beta,r}(r)|\le\eta/2\quad\forall 0\le r\le\eps\notag\\
\Rightarrow\quad&\eta/2\le g_{\beta}(r)\le3\eta/2\quad\mbox{ and }\quad|g_{\beta,r}(r)|\le\eta/2\quad\forall 0\le r\le\eps.
\end{align}
Then by  \eqref{eq-g-beta-representation} for any $\beta, \beta_0\in  [\beta_-(\rho_1),\beta_+(\rho_1)]$ we have
\begin{align}\label{g-beta-beta0-norm}
&|g_{\beta}(r)-g_{\beta_0}(r)|\le\int_0^r|g_{\beta,r}(\rho)-g_{\beta_0,r}(\rho)|\,d\rho\le r\|g_{\beta,r}-g_{\beta_0,r}\|_{L^{\infty}([0,\eps])}\quad\forall 0\le r\le\eps\notag\\
\Rightarrow\quad&\|g_{\beta}-g_{\beta_0}\|_{L^{\infty}([0,\eps])}\le\eps\|g_{\beta,r}-g_{\beta_0,r}\|_{L^{\infty}([0,\eps])}.
\end{align}
Let
\begin{equation}\label{eq:G-beta-norm}
 \mathcal G_\beta:=\max\left\{
 \|g_\beta-g_{\beta_0}\|_{L^\infty([0,\eps])},
 \|g_{\beta,r}-g_{\beta_0,r}\|_{L^\infty([0,\eps])}
 \right\}.
\end{equation}
Since $\beta\in  [\beta_-(\rho_1),\beta_+(\rho_1)]$ , \eqref{beta-upper-bd1} holds.
By the mean value theorem, \eqref{q-and-tilde-def-intro}, \eqref{beta-upper-bd1}, \eqref{eq-g-beta-representation}, \eqref{g-beta-local-bd1} and \eqref{g-beta-beta0-norm}.

\begin{align}\label{g-beta-beta0-difference2}
&|g_{\beta,r}(r)-g_{\beta_0,r}(r)|\notag\\
\le&\frac{ |g_{\beta}^{1-m}(r)-g_{\beta_0}^{1-m}(r)|}{r^{n-1}}\int_0^r\rho^{ n-3+ \frac{n-2-nm}{m}}\left\{a_0(\beta_+(\rho_1)+|\beta|)g_{\beta}(\rho)+|\beta|\rho |g_{\beta,r}(\rho)|\right\}\,d\rho\notag\\
&\quad + \frac{ g_{\beta_0}^{1-m}(r)}{r^{n-1}}\int_0^r\rho^{ n-3+ \frac{n-2-nm}{m}}\left\{a_0|\beta-\beta_0|g_{\beta}(\rho)+|\beta-\beta_0|\rho |g_{\beta,r}(\rho)|\right\}\,d\rho \notag\\
&\quad+ \frac{ g_{\beta_0}^{1-m}(r)}{r^{n-1}}\int_0^r\rho^{ n-3+ \frac{n-2-nm}{m}}\left\{a_0|\beta_+(\rho_1)-\beta_0||g_{\beta}(\rho)-g_{\beta_0}(\rho)|+|\beta_0|\rho |g_{\beta,r}(\rho)-g_{\beta_0,r}(\rho)|\right\}\,d\rho\notag\\
\le&2(a_0+1)(\beta_+(\rho_1)+k_1)\cdot\frac{ |g_{\beta}(r)-g_{\beta_0}(r)|}{(\eta/2)^mr^{n-1}}\int_0^r\rho^{ n-3+ \frac{n-2-nm}{m}}((3\eta/2)+(\eta/2))\,d\rho\notag\\
&\quad +\frac{a_0+1}{n-2-2m}(3\eta/2)^{2-m}r^{\frac{n-2-(n+1)m}{m}}|\beta-\beta_0|\notag\\
 &\quad +\frac{2(a_0+1)}{n-2-2m}(3\eta/2)^{1-m}(\beta_+(\rho_1)+|\beta_0|)r^{\frac{n-2-(n+1)m}{m}}\mathcal G_\beta\notag\\
\le&C_3r^{\frac{n-2-(n+1)m}{m}}\|g_{\beta,r}-g_{\beta_0,r}\|_{L^{\infty}([0,\eps])}+C_3r^{\frac{n-2-(n+1)m}{m}}|\beta-\beta_0|  +C_3r^{\frac{n-2-(n+1)m}{m}}\mathcal G_\beta\notag\\
\le&C_3r^{\frac{n-2-(n+1)m}{m}}|\beta-\beta_0|  +2C_3r^{\frac{n-2-(n+1)m}{m}}\mathcal G_\beta\quad\forall 0<r<\eps
\end{align}
where
\begin{equation*}
C_3=\max\left(1,\frac{8(a_0+1)(\beta_+(\rho_1)+k_1+|\beta_0|+1)(\eta^{1-m}+\eta^{2-m})}{n-2-2m}\right).
\end{equation*}
We now choose $0<\eps<\min\left((8C_3)^{-\frac{m}{n-2-(n+1)m}},1/4\right)$. Then by \eqref{g-beta-beta0-norm} and \eqref{g-beta-beta0-difference2}, 
\begin{align*}
\mathcal G_\beta
\le &C_3\eps^{\frac{n-2-(n+1)m}{m}}|\beta-\beta_0|+2C_3\eps^{\frac{n-2-(n+1)m}{m}}\mathcal G_\beta+\eps\mathcal G_\beta\notag\\
\le&C_3\eps^{\frac{n-2-(n+1)m}{m}}|\beta-\beta_0|+\frac{1}{2}\mathcal G_\beta\notag\\
\Rightarrow\quad\mathcal G_\beta\le&2C_3\eps^{\frac{n-2-(n+1)m}{m}}|\beta-\beta_0|.
\end{align*}
and  \eqref{eq:g-beta-parameter-limit1}  follows.
Finally, Lemma \ref{g-derivative-sign-lem} gives $g_{\beta,r}<0$ on $(0,\eps)$ when $\beta<\beta_+(\rho_1)$, whereas Theorem \ref{thm-existence-inversion} gives $g_{\beta_+}\equiv\eta$.  Consequently
\begin{equation}\label{g-bded2}
0<g_{\beta}(r)\le g_{\beta}(0)=\eta
\quad\forall 0<r\le \eps,\quad
\beta\in[\beta_-(\rho_1),\beta_+(\rho_1)].
\end{equation}
Since $g_{\beta}$ satisfies  \eqref{automatic-g-flux-identity}, by \eqref{automatic-g-flux-identity}  and \eqref{g-bded2},
\begin{equation}\label{gr-bded1}
|rg_{\beta,r}(r)|\le C\eta^{2-m}r^{\frac{n-2-nm}{m}}\quad\forall 0<r\le \eps,\beta\in  [\beta_-(\rho_1),\beta_+(\rho_1)]
\end{equation}
for some constant $C>0$. By \eqref{eq-g-beta-representation} and \eqref{gr-bded1},
\begin{equation}\label{g-lower-bd5}
g_{\beta}(r)\ge\eta-C'\eta^{2-m}r^{\frac{n-2-nm}{m}}\quad\forall 0<r\le \eps,\beta\in  [\beta_-(\rho_1),\beta_+(\rho_1)]
\end{equation}
for some constant $C'>0$.
By \eqref{gr-bded1} and \eqref{g-lower-bd5} we get \eqref{uniform-infinity-log-slope}
and the lemma follows.
\end{proof}

\begin{lem}\label{lem:g-beta-parameter-limit2}
Let $n\ge 3$, $\frac{n-2}{n+1}\leq m< \frac{n-2}{n}$, $\eta>0$, $\rho_1>0$ and $\beta_0\in  [\beta_-(\rho_1),\beta_+(\rho_1)]$. Then there exists a constant $\eps>0$ such that for any $\beta\in  [\beta_-(\rho_1),\beta_+(\rho_1)]$   \eqref{eq-inversion} has a unique positive solution $g_{\beta}\in C^{0, \delta_0}([0,\eps])\cap C^2((0,\eps])$  with $\alpha=\frac{2\beta +\rho_1}{1-m}$ and $\widetilde\alpha$, $\widetilde\beta$,  given by \eqref{q-and-tilde-def-intro}, which satisfies \eqref{g-value-at-origin} and \eqref{eq-initial-m-large} where $\delta_0$, $\delta_1$ and $\zeta$ are given by \eqref{defn-delta0-1} and \eqref{zeta-defn}.
Then
\begin{equation}\label{eq:g-beta-parameter-limit4}
g_{\beta}(r)\to g_{\beta_0}(r)\quad\mbox{ and }\quad r^{\delta_1}g_{\beta,r}(r)\to r^{\delta_1}g_{\beta_0,r}(r)\quad\mbox{ uniformly in $[0,\eps]$}\quad\mbox{ as }\beta\to\beta_0.
\end{equation}
Moreover \eqref{uniform-infinity-log-slope} holds.
\end{lem}
\begin{proof}
For any $\beta\in[\beta_-(\rho_1),\beta_+(\rho_1)]$, let
\begin{equation}\label{eq:zeta-beta-uniform}
 \zeta_{\beta}:=
 \frac{a_0(\beta_+(\rho_1)-\beta)m\eta^{2-m}}{n-2-2m},
 \qquad
 \zeta_*:=\max_{\beta\in[\beta_-(\rho_1),\beta_+(\rho_1)]}\zeta_{\beta}
 =\frac{a_0(\beta_+(\rho_1)-\beta_-(\rho_1))m\eta^{2-m}}{n-2-2m}.
\end{equation}
By an argument analogous to the proof of Lemma 2.3 of \cite{HuS}, there exists a constant $\eps>0$ such that for every $\beta\in[\beta_-(\rho_1),\beta_+(\rho_1)]$, equation \eqref{eq-inversion} has a unique solution $g_{\beta}\in C^{0,\delta_0}([0,\eps])\cap C^2((0,\eps])$ satisfying \eqref{g-value-at-origin}, \eqref{eq-initial-m-large}, \eqref{eq-g-beta-representation}, and
\begin{align}\label{g-beta-local-bd2}
&|g_{\beta}(r)-\eta|\le\eta/2,
 \qquad
 |r^{\delta_1}g_{\beta,r}(r)+\zeta_{\beta}|\le\eta/2
 \quad\forall 0\le r\le\eps\notag\\
\Rightarrow\quad&\eta/2\le g_{\beta}(r)\le3\eta/2,
 \qquad
 r^{\delta_1}|g_{\beta,r}(r)|\le\zeta_*+\eta/2
 \quad\forall 0\le r\le\eps.
\end{align}
Let $\beta\in  [\beta_-(\rho_1),\beta_+(\rho_1)]$. Then by  \eqref{eq-g-beta-representation} ,
\begin{equation*}
|g_{\beta}(r)-g_{\beta_0}(r)|\le\|r^{\delta_1}(g_{\beta,r}-g_{\beta_0,r})\|_{L^{\infty}([0,\eps])}\int_0^r\rho^{-\delta_1}\,d\rho
=\frac{r^{1-\delta_1}}{1-\delta_1}\|r^{\delta_1}(g_{\beta,r}-g_{\beta_0,r})\|_{L^{\infty}([0,\eps])}
\end{equation*}
for any $0\le r\le\eps$. Hence
\begin{equation}\label{g-beta-beta0-norm2}
\|g_{\beta}-g_{\beta_0}\|_{L^{\infty}([0,\eps])}
\le\frac{\eps^{1-\delta_1}}{1-\delta_1}\|r^{\delta_1}(g_{\beta,r}-g_{\beta_0,r})\|_{L^{\infty}([0,\eps])}.
\end{equation}
Since $\beta\in  [\beta_-(\rho_1),\beta_+(\rho_1)]$ , \eqref{beta-upper-bd1} holds.
By the mean value theorem, \eqref{q-and-tilde-def-intro}, \eqref{beta-upper-bd1}, \eqref{eq-g-beta-representation}, \eqref{g-beta-local-bd2} and \eqref{g-beta-beta0-norm2},

\begin{align}\label{g-beta-beta0-difference3}
&r^{\delta_1}|g_{\beta,r}(r)-g_{\beta_0,r}(r)|\notag\\
\le&\frac{ |g_{\beta}^{1-m}(r)-g_{\beta_0}^{1-m}(r)|}{r^{n-2+\frac{n-2-nm}{m}}}\int_0^r\rho^{ n-3+ \frac{n-2-nm}{m}}\left\{a_0(\beta_+(\rho_1)+|\beta|)g_{\beta}(\rho)+|\beta|\rho |g_{\beta,r}(\rho)|\right\}\,d\rho\notag\\
&\quad + \frac{ g_{\beta_0}^{1-m}(r)}{r^{n-2+\frac{n-2-nm}{m}}}\int_0^r\rho^{ n-3+ \frac{n-2-nm}{m}}\left\{a_0|\beta-\beta_0|g_{\beta}(\rho)+|\beta-\beta_0|\rho |g_{\beta,r}(\rho)|\right\}\,d\rho \notag\\
&\quad+ \frac{ g_{\beta_0}^{1-m}(r)}{r^{n-2+\frac{n-2-nm}{m}}}\int_0^r\rho^{ n-3+ \frac{n-2-nm}{m}}\left\{a_0|\beta_+(\rho_1)-\beta_0||g_{\beta}(\rho)-g_{\beta_0}(\rho)|+|\beta_0|\rho |g_{\beta,r}(\rho)-g_{\beta_0,r}(\rho)|\right\}\,d\rho\notag\\
\le&(a_0+1)(\beta_+(\rho_1)+k_1)\cdot\frac{ |g_{\beta}(r)-g_{\beta_0}(r)|}{(\eta/2)^m{r^{n-2+\frac{n-2-nm}{m}}}}\int_0^r\rho^{ n-3+ \frac{n-2-nm}{m}}((3\eta/2)+\zeta_*+\eta/2)\,d\rho\notag\\
&\quad +\frac{a_0+1}{n-2-nm}(3\eta/2)^{1-m}\left(2\eta+\zeta_*\right)|\beta-\beta_0|
 +\frac{a_0|\beta_0|(3\eta/2)^{1-m}}{n-2-nm}\|g_{\beta}-g_{\beta_0}\|_{L^{\infty}([0,\eps])}\notag\\
& +\frac{2(a_0+1)(\beta_+(\rho_1)+|\beta_0|)(3\eta/2)^{1-m}}{r^{n-2+\frac{n-2-nm}{m}}}\|r^{\delta_1}(g_{\beta,r}-g_{\beta_0,r})\|_{L^{\infty}([0,\eps])}\int_0^r\rho^{ n-2+ \frac{n-2-nm}{m}-\delta_1}\,d\rho\notag\\
\le&2C_3\|g_{\beta}-g_{\beta_0}\|_{L^{\infty}([0,\eps])}+C_3|\beta-\beta_0|+C_3r^{1-\delta_1}\|r^{\delta_1}(g_{\beta,r}-g_{\beta_0,r})\|_{L^{\infty}([0,\eps])}\notag\\
\le&\frac{2C_3\eps^{1-\delta_1}}{1-\delta_1}\|r^{\delta_1}(g_{\beta,r}-g_{\beta_0,r})\|_{L^{\infty}([0,\eps])}+C_3|\beta-\beta_0|  +C_3\eps^{1-\delta_1}\|r^{\delta_1}(g_{\beta,r}-g_{\beta_0,r})\|_{L^{\infty}([0,\eps])}\notag\\
\le&C_4\eps^{1-\delta_1}\|r^{\delta_1}(g_{\beta,r}-g_{\beta_0,r})\|_{L^{\infty}([0,\eps])}+C_3|\beta-\beta_0|\quad\forall 0<r<\eps
\end{align}
where
\begin{equation*}
C_3=\frac{6(a_0+1)(\beta_+(\rho_1)+k_1+|\beta_0|+1)(\eta^{1-m}+\eta^{2-m}+\zeta_*\eta^{1-m}+\zeta_*\eta^{-m})}{n-2-nm}
\end{equation*}
and
\begin{equation*}
C_4=C_3\left(\frac{3-\delta_1}{1-\delta_1}\right).
\end{equation*}
We now choose $0<\eps<(4C_4)^{-\frac{1}{1-\delta_1}}$. Then by \eqref{g-beta-beta0-difference3},
\begin{align}\label{g-beta-r-beta-ineqn2}
&\|r^{\delta_1}(g_{\beta,r}-g_{\beta_0,r})\|_{L^{\infty}([0,\eps])}\le (1/4)\|r^{\delta_1}(g_{\beta,r}-g_{\beta_0,r})\|_{L^{\infty}([0,\eps])}+C_3|\beta-\beta_0|\notag\\
\Rightarrow\quad&\|r^{\delta_1}(g_{\beta,r}-g_{\beta_0,r})\|_{L^{\infty}([0,\eps])}\le 2C_3|\beta-\beta_0|.
\end{align}
By \eqref{g-beta-beta0-norm2} and \eqref{g-beta-r-beta-ineqn2}, we get \eqref{eq:g-beta-parameter-limit4}.
Finally by the same argument as Lemma \ref{lem:g-beta-parameter-limit}, \eqref{gr-bded1} and \eqref{g-lower-bd5} hold. Then by \eqref{gr-bded1} and \eqref{g-lower-bd5} we get \eqref{uniform-infinity-log-slope}
  hold, and the lemma follows.

\end{proof}

The preceding local estimates, combined with continuous dependence for the regular first-order system associated with \eqref{eq-inversion} away from $r=0$, give the following two global compact-dependence statements.

\begin{thm}\label{thm:parameter-limit-convergence-g-eqn}
Let $n\ge 3$, $0<m<\frac{n-2}{n+1}$, $\eta>0$,  $\rho_1>0$, $\beta_-(\rho_1)\le\beta_0\le\beta_+(\rho_1)$  and $\alpha_0=\frac{2\beta_0+\rho_1}{1-m}$. Let $R_g=R_g(\beta_0,\eta)\in(0,\infty]$
and $g_{\beta_0}\in C^1([0,R_g))\cap C^2((0,R_g))$   be the unique maximal positive branch of \eqref{eq-inversion} in $[0,R_g)$ given by Theorem \ref{thm-existence-self-similar2} with $\alpha=\alpha_0$, $\beta=\beta_0$, and $\widetilde\alpha_0$, $\widetilde\beta_0$,  given by \eqref{q-and-tilde-def-intro}, which satisfies \eqref{g-value-at-origin} and
\eqref{eq-initial-m-small}.  If $0<R<R_g(\beta_0,\eta)$, then there exists a constant $\delta_2>0$ such that $R_g(\beta,\eta)>R$ for all $\beta\in I_{\delta_2}:=[\beta_0-\delta_2,\beta_0+\delta_2]\cap [\beta_-(\rho_1),\beta_+(\rho_1)]$, and
\begin{equation}\label{eq:g-beta-parameter-limit2}
g_{\beta}(r)\to g_{\beta_0}(r)\quad\mbox{ and }\quad g_{\beta,r}(r)\to g_{\beta_0,r}(r)\quad\mbox{ uniformly in $[0,R]$ }\quad\mbox{ as }\beta\to\beta_0.
\end{equation}
Moreover \eqref{uniform-infinity-log-slope} holds.

\end{thm}

\begin{proof}
Let   $\eps>0$ and $g_{\beta}$, $\beta\in  [\beta_-(\rho_1),\beta_+(\rho_1)]$, be as given by Lemma  \ref{lem:g-beta-parameter-limit}. Without loss of generality we may assume $R>\eps$. Then  by Lemma  \ref{lem:g-beta-parameter-limit},  \eqref{eq:g-beta-parameter-limit1} and \eqref{uniform-infinity-log-slope}  holds.
By \eqref{eq:g-beta-parameter-limit1},
\begin{equation}\label{eq:g-beta-parameter-limit7}
g_\beta(\eps)\to g_{\beta_0}(\eps),
 \quad\mbox{ and }\quad
 g_{\beta,r}(\eps)\to g_{\beta_0,r}(\eps)\quad\mbox{ as }\beta\to\beta_0.
\end{equation}
Since $g_{\beta_0}$ has a positive minimum on $[\eps,R]$, by \eqref{eq:g-beta-parameter-limit7} and the 
continuous dependence of solutions $(g,g_r)$ of  regular first-order system on $[\eps,R]$ with respect to its parameters that there exists a constant $\delta_2>0$ such that for any $\beta\in I_{\delta_2}:=[\beta_0-\delta_2,\beta_0+\delta_2]\cap [\beta_-(\rho_1),\beta_+(\rho_1)]$ we have $R_g(\beta,\eta)>R$ and
\begin{equation}\label{eq:g-beta-parameter-limit8}
(g_\beta,g_{\beta,r})\longrightarrow
 (g_{\beta_0},g_{\beta_0,r})
 \quad\text{uniformly on }[\eps,R]\quad\mbox{ as }\beta\to\beta_0.
\end{equation}
By \eqref{eq:g-beta-parameter-limit4} and \eqref{eq:g-beta-parameter-limit8} we get \eqref{eq:g-beta-parameter-limit2}.  
\end{proof}

\begin{thm}\label{thm:parameter-limit-convergence-g-eqn2}
Let $n\ge 3$, $\frac{n-2}{n+1}\le m<\frac{n-2}{n}$, $\eta>0$,  $\rho_1>0$, $\beta_-(\rho_1)\le\beta_0\le\beta_+(\rho_1)$  and $\alpha_0=\frac{2\beta_0+\rho_1}{1-m}$. Let $R_g=R_g(\beta_0,\eta)\in(0,\infty]$
and $g_{\beta_0}\in C^{0, \delta_0}([0,R_g))\cap C^2((0,R_g))$ be the maximal positive branch given by Theorem \ref{thm-existence-self-similar2} with $\alpha=\alpha_0$, $\beta=\beta_0$, and $\widetilde\alpha_0$, $\widetilde\beta_0$ given by \eqref{q-and-tilde-def-intro}.  Suppose that $g_{\beta_0}$ satisfies \eqref{g-value-at-origin} and \eqref{eq-initial-m-large}, where $\delta_0$, $\delta_1$, and $\zeta$ are given by \eqref{defn-delta0-1} and \eqref{zeta-defn}.
If $0<R<R_g(\beta_0,\eta)$, then there exists a constant $\delta_2>0$ such that $R_g(\beta,\eta)>R$ for all $\beta\in I_{\delta_2}:=[\beta_0-\delta_2,\beta_0+\delta_2]\cap [\beta_-(\rho_1),\beta_+(\rho_1)]$, and
\begin{equation}\label{eq:g-beta-parameter-limit-high-m}
 g_{\beta}\longrightarrow g_{\beta_0}
 \quad\text{and}\quad
 r^{\delta_1}g_{\beta,r}\longrightarrow r^{\delta_1}g_{\beta_0,r}
 \quad\text{uniformly on }[0,R]
 \quad\text{as }\beta\to\beta_0.
\end{equation}
Moreover, for every $\varepsilon\in(0,R)$,
\begin{equation}\label{eq:g-beta-parameter-limit-high-m-away}
 g_{\beta,r}\longrightarrow g_{\beta_0,r}
 \quad\text{uniformly on }[\varepsilon,R]
 \quad\text{as }\beta\to\beta_0,
\end{equation}
and \eqref{uniform-infinity-log-slope} holds.
\end{thm}
\begin{proof}
Let   $\eps>0$ and $g_{\beta}$, $\beta\in  [\beta_-(\rho_1),\beta_+(\rho_1)]$, be as given by Lemma  \ref{lem:g-beta-parameter-limit2}. Without loss of generality we may assume $R>\eps$. Then  by Lemma  \ref{lem:g-beta-parameter-limit2},  \eqref{uniform-infinity-log-slope}  and \eqref{eq:g-beta-parameter-limit4} holds.
By \eqref{eq:g-beta-parameter-limit4}, \eqref{eq:g-beta-parameter-limit7} holds.
Since $g_{\beta_0}$ has a positive minimum on $[\eps,R]$, by \eqref{eq:g-beta-parameter-limit7} and the 
continuous dependence of solutions $(g,g_r)$ of  regular first-order system on $[\eps,R]$ with respect to its parameters that there exists a constant $\delta_2>0$ such that for any $\beta\in I_{\delta_2}:=[\beta_0-\delta_2,\beta_0+\delta_2]\cap [\beta_-(\rho_1),\beta_+(\rho_1)]$ we have $R_g(\beta,\eta)>R$ and \eqref{eq:g-beta-parameter-limit8} holds.
By \eqref{eq:g-beta-parameter-limit4} and \eqref{eq:g-beta-parameter-limit8} we get \eqref{eq:g-beta-parameter-limit-high-m} and \eqref{eq:g-beta-parameter-limit-high-m-away}.  
\end{proof}

\section{Slope-amplitude equation and its properties}
\setcounter{equation}{0}
\setcounter{thm}{0}

In this section we will  derive the slope-amplitude equation that corresponds to the radial solution of \eqref{f-elliptic-eqn} and obtain its various properties.  We will first transform a solution $f$ of \eqref{f-ode} into an equation for its slope as a function of the amplitude.
 Let
\begin{equation}\label{eq:ab-sigma}
 a:=\frac{2m}{1-m},\quad b:=n-2-a,\quad\sigma:=\frac{1-m}{m}\quad\mbox{ and }\quad E:=\frac{m}{1-m}.
\end{equation}
The assumption \(0<m<(n-2)/n\) gives
\begin{equation}\label{eq:ab-identities}
 a>0,\quad b=\frac{n-2-nm}{1-m}>0,\quad a+b=n-2,\quad a\sigma=2.
\end{equation}
Unless stated otherwise  we will now assume $\rho_1=1$, $\alpha=\frac{2\beta+1}{1-m}$ and $\beta\in[\beta_-,\beta_+]$ where
\begin{equation}\label{eq:beta-endpoints}
\beta_-:=\beta_-(1)=-\frac{E}{a}=-\frac12\quad\mbox{ and }\quad\beta_+:=\beta_+(1)=\frac{E}{b}=\frac{m}{n-2-nm}.
\end{equation}
For a positive radial function \(f=f(r)\), let
\begin{equation}
 t:=\log r,
 \qquad
 W(t):=r^a f(r)^m=e^{at}f(e^t)^m.
 \label{eq:W-def}
\end{equation}
Equivalently,
\begin{equation}\label{eq:f-from-W}
 f(r)=r^{-\frac{2}{1-m}}W(\log r)^{1/m}.
\end{equation}

\begin{lem}\label{lem:lienard-reduction}
Let \(\rho _1=1\) and \(\alpha=\frac{2\beta+1}{1-m}\).  Then \(f\) solves
\eqref{f-ode} if and only if \(W\) defined by \eqref{eq:W-def} solves
\begin{equation}\label{eq:lienard}
 W''+\bigl(b-a+\beta W^\sigma\bigr)W'
      -abW+EW^{1+\sigma}=0.
\end{equation}
\end{lem}
\begin{proof}
By  \eqref{eq:f-from-W}, 
\begin{equation}\label{fm-laplacian}
 (f^m)_{rr}
 +\frac{n-1}{r}(f^m)_r
 =r^{-\frac{2}{1-m}}
 \left(W''+(n-2-2a)W'+a(a-n+2)W\right).
\end{equation}
and
\begin{equation}\label{eq:r-f-r-eqn}
rf_r=r^{-\frac{2}{1-m}}\left(-\frac{2}{1-m}W^{1/m}+\frac1mW^{1/m-1}W'\right).
\end{equation}
Hence by  \eqref{f-ode}, \eqref{eq:f-from-W}, \eqref{fm-laplacian} and  \eqref{eq:r-f-r-eqn}
we get \eqref{eq:lienard}.  Since every step is reversible, the lemma follows.
\end{proof}

We next  define
\begin{equation}\label{eq:h-H}
h(s):=-abs+Es^{1+\sigma},\quad H(s):=-\frac{ab}{2}s^2+\frac{E}{2+\sigma}s^{2+\sigma}\quad\forall s\ge 0.
\end{equation}
Then \(H'(s)=h(s)\) on $[0,\infty)$ and the unique positive zero of \(h\) is
\begin{equation}\label{eq:s0}
 s_0:=\left(\frac{ab}{E}\right)^{1/\sigma}=\left(\frac{2(n-2-nm)}{1-m}\right)^{m/(1-m)}>0
\end{equation}
with
\begin{equation}\label{eq:H-s0}
H(s_0)=-\frac{ab\sigma}{2(2+\sigma)}s_0^2<0.
\end{equation}
Moreover
\begin{equation}\label{eq:h-sign-s0}
 h(s)<0\quad\forall 0<s<s_0,
 \quad h(s_0)=0,
 \quad h(s)>0\quad\forall s>s_0.
\end{equation}
For any $\beta\in[\beta_-,\beta_+]$ and $\alpha=\frac{2\beta+1}{1-m}$, let $f^o_\beta$ be the maximal origin branch with $\eta_0=1$ defined on $[0,R_o)$, $R_o=R_o(\beta,1)$, given by Theorem \ref{thm-existence-self-similar1}. For any $\beta\in(-\infty,\beta_+]$ and the same normalized relation for $\alpha$, let $f^{\infty}_\beta$ be the maximal outer branch with $\eta=1$ defined on $(R_g^{-1},\infty)$, $R_g:=R_g(\beta,1)$, given by Theorem \ref{thm-existence-self-similar2}.  Let
\[
  T^{\mathrm o}_\beta:=\log R_{\mathrm o}(\beta,1)\quad\mbox{ if }R_{\mathrm o}(\beta,1)<\infty\quad\mbox{ and }\quad T^{\mathrm o}_\beta:=\infty\quad\mbox{ if }R_{\mathrm o}(\beta,1)=\infty
\]
and
\[
  T^{\infty}_\beta:=-\log R_g(\beta,1)\quad\mbox{ if }R_g(\beta,1)<\infty\quad\mbox{ and }\quad T^{\infty}_\beta:=-\infty\quad\mbox{ if }R_g(\beta,1)=\infty.
\]
Let $g_{\beta}$ be given by  Theorem \ref{thm-existence-self-similar2} which is related to $f^{\infty}_\beta$ by \eqref{outer-profile-definition-intro}. Let
\begin{equation}\label{eq:W-origin-infinity}
 W^o_\beta(t)=e^{at}f^o_\beta(e^t)^m\quad\forall t<T^{\mathrm o}_\beta,
 \qquad
 W^{\infty}_\beta(t)=e^{at}f^{\infty}_\beta(e^t)^m\quad\forall t>T^{\infty}_\beta
\end{equation}
where $a$ is given by \eqref{eq:ab-sigma}.
By Lemma \ref{lem:lienard-reduction} and \eqref{eq:W-origin-infinity}, $W^o_\beta$ and $W^{\infty}_\beta$ satisfy \eqref{eq:lienard} in $(-\infty,T^{\mathrm o}_\beta)$ and $(T^{\infty}_\beta,\infty)$ respectively.
By \eqref{outer-profile-definition-intro} and \eqref{eq:W-origin-infinity},
\begin{equation}\label{eq:W-infinity1}
W^{\infty}_\beta(t)=e^{at}(e^{-\frac{(n-2)}{m}t}g_\beta(e^{-t}))^m=e^{-bt}g_\beta(e^{-t})^m\quad\forall t>T^{\infty}_\beta
\end{equation}
where $b$ is given by \eqref{eq:ab-sigma}.
Differentiating \eqref{eq:W-origin-infinity} and \eqref{eq:W-infinity1} with respect to $t$, we get
\begin{equation}\label{eq:w0-beta-eqn}
(W^o_\beta)'(t)=ae^{at}f^o_\beta(e^t)^m+me^{(a+1)t}f^o_\beta(e^t)^{m-1}f^o_{\beta,r}(e^t)=W^o_\beta(t)\left(a+m\frac{e^tf^o_{\beta,r}(e^t)}{f^o_\beta(e^t)}\right)\quad\forall t<T^{\mathrm o}_\beta
\end{equation}
and
\begin{equation}\label{eq:w-infty-beta-eqn}
(W^{\infty}_\beta)'(t)=-be^{-bt}g_\beta(e^{-t})^m-me^{-(b+1)t}g_\beta(e^{-t})^{m-1}g_{\beta,r}(e^{-t})=W^{\infty}_\beta(t)\left(-b-m\frac{e^{-t}g_{\beta,r}(e^{-t})}{g_\beta(e^{-t})}\right)\quad\forall t>T^{\infty}_\beta.
\end{equation}
By Theorem \ref{thm:parameter-limit-convergence}, Theorem \ref{thm:parameter-limit-convergence-g-eqn}, Theorem \ref{thm:parameter-limit-convergence-g-eqn2}, \eqref{eq:w0-beta-eqn} and \eqref{eq:w-infty-beta-eqn}, the endpoint asymptotics are
\begin{equation}\label{W-endpoint-asymptotics}
 \frac{(W^o_\beta)'(t)}{W^o_\beta(t)}\to a
 \quad\mbox{ as } t\to-\infty,
 \quad\mbox{ and }\quad
 \frac{(W^{\infty}_\beta)'(t)}{W^{\infty}_\beta(t)}\to-b
 \quad\mbox{ as }t\to\infty,
\end{equation}
uniformly for $\beta\in[\beta_-,\beta_+]$. We now define the first-turn times
\begin{align}
 \tau_P(\beta)&:=\sup\left\{\tau<T^{\mathrm o}_\beta:
 (W^o_\beta)'(t)>0\text{ for all }t<\tau\right\},
 \label{tau-P-definition}\\
 \tau_Q(\beta)&:=\inf\left\{\tau>T^{\infty}_\beta:
 (W^{\infty}_\beta)'(t)<0\text{ for all }t>\tau\right\}.
 \label{tau-Q-definition}
\end{align}
The sets are nonempty by \eqref{W-endpoint-asymptotics}.  Let
\begin{equation}\label{turning-amplitudes-definition}
 M_P(\beta):=\lim_{t\uparrow\tau_P(\beta)}W^o_\beta(t),
 \qquad
 M_Q(\beta):=\lim_{t\downarrow\tau_Q(\beta)}W^{\infty}_\beta(t).
\end{equation}
Only these half-branches are used below.

\begin{lem}[First-turn alternatives]\label{lem:first-turn-alternatives}
For $\beta\in(\beta_-,\beta_+]$, the origin half-branch $W^o_{\beta}$ is strictly increasing on
$(-\infty,\tau_P(\beta))$.  Either $\tau_P(\beta)<\infty$ and
\begin{equation}\label{P-finite-turn-data}
 (W^o_\beta)'(\tau_P)=0,
 \qquad M_P(\beta)\geq s_0,
\end{equation}
or $\tau_P(\beta)=\infty$.  When $\beta=\beta_-$,
\begin{equation}\label{P-endpoint-exact-orbit}
 W^o_{\beta_-}(t)=e^{at},
 \qquad M_P(\beta_-)=\infty.
\end{equation}

For $\beta\in[\beta_-,\beta_+)$, the outer half-branch $W^{\infty}_{\beta}$ is strictly decreasing on
$(\tau_Q(\beta),\infty)$.  Either $\tau_Q(\beta)>-\infty$ and
\begin{equation}\label{Q-finite-turn-data}
 (W^{\infty}_\beta)'(\tau_Q)=0,
 \qquad M_Q(\beta)\geq s_0,
\end{equation}
or $\tau_Q(\beta)=-\infty$.  When $\beta=\beta_+$,
\begin{equation}\label{Q-endpoint-exact-orbit}
 W^{\infty}_{\beta_+}(t)=e^{-bt},
 \qquad M_Q(\beta_+)=\infty.
\end{equation}
\end{lem}
\begin{proof}
Let $\beta\in(\beta_-,\beta_+]$ and suppose $\tau_P(\beta)<\infty$.  We first claim that $\tau_P(\beta)<T^{\mathrm o}_\beta$.  Otherwise $\tau_P(\beta)=T^{\mathrm o}_\beta<\infty$, so the finite-zero alternative in Theorem \ref{thm-existence-self-similar1} gives $W^o_\beta(t)\to0$ as $t\uparrow T^{\mathrm o}_\beta$.  This is impossible because $W^o_\beta$ is positive and strictly increasing on $(-\infty,\tau_P(\beta))$.  Thus the branch is defined on a neighbourhood of $\tau_P(\beta)$, and the definition of the first turn gives
\begin{equation}\label{W1''-limit-pt-sign}
 (W^o_\beta)'(\tau_P(\beta))=0\quad\mbox{ and }\quad (W^o_\beta)''(\tau_P(\beta))\leq0.
\end{equation}
By Lemma \ref{lem:lienard-reduction}, $W^o_\beta$ satisfies
\eqref{eq:lienard} in a neighbourhood of $\tau_P(\beta)$. Putting $W=W^o_\beta$ and  $t=\tau_P(\beta)$ in \eqref{eq:lienard}, by \eqref{eq:h-sign-s0} and \eqref{W1''-limit-pt-sign},
\begin{align*}
&0\ge (W^o_\beta)''(\tau_P(\beta))=-h(W^o_\beta(\tau_P(\beta)))\notag\\
\Rightarrow\quad&W^o_\beta(\tau_P(\beta))\ge s_0
\end{align*}
and \eqref{P-finite-turn-data} follows. When $\beta=\beta_-$, $f^o_{\beta}\equiv 1$ and hence \eqref{P-endpoint-exact-orbit} holds.

Suppose $\beta\in[\beta_-,\beta_+)$ and $\tau_Q(\beta)>-\infty$.  We first claim that $T^{\infty}_\beta<\tau_Q(\beta)$.  If equality held, then $R_g<\infty$ and Theorem \ref{thm-existence-self-similar2} would give $W^{\infty}_\beta(t)\to0$ as $t\downarrow T^{\infty}_\beta$.  This is incompatible with positivity and strict decrease on $(\tau_Q(\beta),\infty)$.  Hence the branch is defined on a neighbourhood of $\tau_Q(\beta)$, and
\begin{equation}\label{P-finite-turn-data0}
 (W^{\infty}_\beta)'(\tau_Q)=0\quad\mbox{ and }\quad (W^{\infty}_\beta)''(\tau_Q)\le0.
\end{equation}
By Lemma \ref{lem:lienard-reduction}, $W^{\infty}_\beta$ satisfies
\eqref{eq:lienard} in a neighbourhood of $\tau_Q(\beta)$. Putting $W=W^{\infty}_\beta$ and $t=\tau_Q(\beta)$ in \eqref{eq:lienard}, by \eqref{eq:h-sign-s0} and \eqref{P-finite-turn-data0},
\begin{align*}
&0\ge (W^{\infty}_\beta)''(\tau_Q(\beta))=-h(W^{\infty}_\beta(\tau_Q(\beta)))\notag\\
\Rightarrow\quad&W^{\infty}_\beta(\tau_Q(\beta))\ge s_0
\end{align*}
and \eqref{Q-finite-turn-data} follows. When $\beta=\beta_+$, $g_{\beta}\equiv 1$ and hence \eqref{Q-endpoint-exact-orbit} holds.

\end{proof}

For \(\beta\in [\beta_-,\beta_+]\), the function $W^o_\beta: (-\infty,\tau_P)\to (0,M_P(\beta))$  is strictly increasing and has an inverse  
\begin{equation*}
t=(W^o_\beta)^{-1}:(0,M_P(\beta))\to (-\infty,\tau_P(\beta)).
\end{equation*}
Hence on the increasing origin half-branch $W^o_\beta$ we define
\begin{equation}\label{P-beta-defn}
 P_\beta(s):=(W^o_\beta)'
 \bigl((W^o_\beta)^{-1}(s)\bigr),
 \quad\forall 0<s<M_P(\beta),
\end{equation}
Then 
\begin{equation}\label{P-defn}
(W^o_\beta)'(t)=P_{\beta}(W^o_\beta(t))>0\quad\forall t<\tau_P(\beta).
\end{equation}
Since $W^o_\beta$ satisfies \eqref{eq:lienard}, by \eqref{P-defn} we get
\begin{equation}\label{eq:P-equation}
P_\beta P_\beta'+c_\beta(s)P_\beta+h(s)=0\quad\forall 0<s<M_P(\beta)
\end{equation}
where
\begin{equation}\label{eq:c-beta}
 c_\beta(s):=b-a+\beta s^\sigma.
\end{equation}
Similarly the function $W^{\infty}_\beta: (\tau_Q,\infty)\to (0,M_Q(\beta))$  is strictly decreasing and has an inverse  
\begin{equation*}
t=(W^{\infty}_\beta)^{-1}:(0,M_Q(\beta))\to (\tau_Q(\beta),\infty).
\end{equation*}
Hence on the decreasing outer half-branch $W^{\infty}_\beta$ we define
\begin{equation}\label{Q-beta-defn}
 Q_\beta(s):=-(W^{\infty}_\beta)'
 \bigl((W^{\infty}_\beta)^{-1}(s)\bigr),
 \quad\forall 0<s<M_Q(\beta).
\end{equation}
Then
\begin{equation}\label{Q-defn}
-(W^{\infty}_\beta)'(t)=Q_{\beta}(W^{\infty}_\beta(t))>0\quad\forall t>\tau_Q(\beta).
\end{equation}
Since $W^{\infty}_\beta$ satisfies \eqref{eq:lienard}, by \eqref{Q-defn} we get
\begin{equation}\label{eq:Q-equation}
 Q_\beta Q_\beta'-c_\beta(s)Q_\beta+h(s)=0\quad\forall 0<s<M_Q(\beta).
\end{equation}
Note that by \eqref{f-fr-origin-value} and \eqref{eq:W-origin-infinity},
\begin{equation}\label{eq:w1-limit-zero}
W^o_\beta(t)\to 0\quad\mbox{ as }t\to-\infty.
\end{equation}
and by \eqref{g-value-at-origin} and \eqref{eq:W-infinity1},
\begin{equation}\label{eq:W-infinity2}
W^{\infty}_\beta(t)\to0\quad\mbox{ as }t\to\infty.
\end{equation}
 By Theorem \ref{thm:parameter-limit-convergence}, Theorem \ref{thm:parameter-limit-convergence-g-eqn}, Theorem \ref{thm:parameter-limit-convergence-g-eqn2}, \eqref{eq:w0-beta-eqn}, \eqref{eq:w-infty-beta-eqn}, \eqref{eq:w1-limit-zero} and  \eqref{eq:W-infinity1},
\begin{equation}\label{eq:w1'-limit-zero}
(W^o_\beta)'(t)\to 0\quad\mbox{ as }t\to-\infty
\end{equation}
and 
\begin{equation}\label{eq:w2'-limit-zero}
(W^{\infty}_\beta)'(t)\to 0\quad\mbox{ as }t\to\infty.
\end{equation}
By \eqref{P-defn}, \eqref{Q-defn}, \eqref{eq:w1-limit-zero}, \eqref{eq:W-infinity2}, \eqref{eq:w1'-limit-zero} and \eqref{eq:w2'-limit-zero}, we can extend $P_{\beta}$ to a continuous function on $[0,M_P(\beta))$ by letting  $P_{\beta}(0)=0$ and  we can extend $Q_{\beta}$ to a continuous function on $[0,M_Q(\beta))$ by letting  $Q_{\beta}(0)=0$.

A global matched profile is obtained by joining the increasing half-orbit described by $P_\beta$ to the decreasing half-orbit described by $Q_\beta$ at a common first-turn amplitude.  The joined orbit is a pulse solution of the autonomous equation \eqref{eq:lienard}; the functions $P_\beta$ and $Q_\beta$ themselves satisfy the two distinct slope--amplitude equations \eqref{eq:P-equation} and \eqref{eq:Q-equation}.  We now establish various properties of $P_{\beta}$ and $Q_{\beta}$ for any $\beta\in(\beta_-,\beta_+)$.

\begin{lem}[Linear slope bounds]\label{lem:linear-comparison}
For $\beta\in(\beta_-,\beta_+]$,
\begin{equation}\label{eq:P-below-line}
 0<P_\beta(s)<as\qquad\forall 0<s<M_P(\beta),
\end{equation}
and for $\beta\in[\beta_-,\beta_+)$,
\begin{equation}\label{eq:Q-below-line}
 0<Q_\beta(s)<bs\qquad\forall 0<s<M_Q(\beta).
\end{equation}
Furthermore
\begin{equation}\label{P-Q-zero-slopes}
 \lim_{s\searrow0}\frac{P_\beta(s)}s=a,
 \qquad
 \lim_{s\searrow0}\frac{Q_\beta(s)}s=b,\quad\mbox{ uniformly for $\beta\in[\beta_-,\beta_+]$}
\end{equation}
 and
\begin{equation}\label{endpoint-linear-slopes}
 P_{\beta_-}(s)=as,
 \qquad
 Q_{\beta_+}(s)=bs
 \qquad\forall s>0
\end{equation}
where $a$ and $b$ are given by \eqref{eq:ab-sigma}.
\end{lem}

\begin{proof}
For every \(\beta\in(\beta_-,\beta_+]\), by \eqref{f-derivative-negative2},  \eqref{eq:w0-beta-eqn} and \eqref{P-defn} we get \eqref{eq:P-below-line}. 

We now let \(\beta\in[\beta_-,\beta_+)\). 
By \eqref{g-monotonicity-intro}, \eqref{Q-defn}, and \eqref{eq:w-infty-beta-eqn},
\[
 \frac{Q_{\beta}(W^{\infty}_\beta(t))}{W^{\infty}_\beta(t)}=-\frac{(W^{\infty}_\beta)'(t)}{W^{\infty}_\beta(t)}
 =b+m\frac{e^{-t} g_{\beta,r}(e^{-t})}{g_\beta(e^{-t})}<b\quad\forall t>\tau_Q(\beta)
\]
and \eqref{eq:Q-below-line} follows. By \eqref{W-endpoint-asymptotics} we get \eqref{P-Q-zero-slopes}. By \eqref{P-endpoint-exact-orbit} and \eqref{Q-endpoint-exact-orbit} we get \eqref{endpoint-linear-slopes}.

\end{proof}

\begin{lem}[Finite turning amplitudes and endpoint expansions]\label{lem:turning-finite}
Let $s_0$ be given by \eqref{eq:s0}. Then the following hold.
\begin{enumerate}[(i)]

\item For every $\beta\in(\beta_-,\beta_+]$,
\begin{equation}\label{eq:MP-finite}
 s_0\leq M_P(\beta)<\infty,
\end{equation}
and $P_\beta$ can be extended to a continuous function on $[0,M_P(\beta)]$ by setting $P_{\beta}(M_P(\beta))=0$.
\item For every $\beta\in[\beta_-,\beta_+)$,
\begin{equation}\label{eq:MQ-finite}
 s_0\leq M_Q(\beta)<\infty.
\end{equation}
and $Q_\beta$ can be extended to a continuous function on $[0,M_Q(\beta)]$ by setting $Q_{\beta}(M_Q(\beta))=0$.

\item If $M_P(\beta)>s_0$, then
\begin{equation}\label{eq:P-turning-expansion}
 P_\beta(s)=
 \sqrt{2h(M_P(\beta))(M_P(\beta)-s)}
 +o\bigl(\sqrt{M_P(\beta)-s}\bigr).
\end{equation}

\item If $M_Q(\beta)>s_0$, then
\begin{equation}\label{eq:Q-turning-expansion}
 Q_\beta(s)=
 \sqrt{2h(M_Q(\beta))(M_Q(\beta)-s)}
 +o\bigl(\sqrt{M_Q(\beta)-s}\bigr).
\end{equation}

\item If $M_P(\beta)=s_0$, then $\tau_P(\beta)=\infty$ and
$W^o_\beta(t)\nearrow s_0$ as $t\to\infty$.

\item If $M_Q(\beta)=s_0$, then $\tau_Q(\beta)=-\infty$ and
$W^{\infty}_\beta(t)\to s_0$ as $t\to-\infty$.

\end{enumerate}

\end{lem}

\begin{proof}
\noindent $\text{\bf Proof of (i)}$: We first prove the right-hand side inequality of (i).
Suppose there exists $\beta\in (\beta_-,\beta_+]$ such that $M_P(\beta)=\infty$. By \eqref{eq:P-equation},
\begin{equation}\label{eq:P-square}
 (P_\beta^2)'(s)=-2(b-a+\beta s^\sigma)P_\beta(s)+2abs-2Es^{1+\sigma}\quad\forall s>0.
\end{equation}
We now divide the proof into two cases.

\noindent $\underline{\text{\bf Case 1}}$: \(0\le \beta\le\beta_+\).

\noindent By \eqref{eq:P-below-line} and \eqref{eq:P-square},
\begin{equation}\label{P2-ineqn2}
 (P_\beta^2)'(s)\leq C_1 s-2Es^{1+\sigma}=s(C_1 -2Es^{\sigma})
 \quad\forall s> 0
\end{equation}
where $C_1=2a^2+2ab>0$. 

\noindent $\underline{\text{\bf Case 2}}$: \(\beta_-<\beta<0\).

\noindent
Then by \eqref{eq:beta-endpoints}, \eqref{eq:P-below-line} and \eqref{eq:P-square},
\begin{equation}\label{P2-ineqn3}
(P_\beta^2)'(s)\leq C_1s-2(E+a\beta)s^{1+\sigma}=s(C_1-2a(\beta-\beta_-)s^{\sigma})
\quad\forall s> 0.
\end{equation}
Let 
\begin{equation}\label{s1-defn}
s_1=\max\left(1,(C_1/E)^{1/\sigma},\left(C_1/[a(\beta-\beta_-)]\right)^{1/\sigma}\right).
\end{equation}
 Then by \eqref{P2-ineqn2}, \eqref{P2-ineqn3} and \eqref{s1-defn},
\begin{equation*}
 (P_\beta^2)'(s)\le -C_1\quad\forall s\ge s_1
\end{equation*}
which implies that $P_\beta(s)^2<0$ for sufficiently large $s$ and a contradiction follows. 
Hence the right-hand inequality of (i) holds.

We will now prove the left-hand side inequality of (i). Observe that  by \eqref{eq:P-equation} and \eqref{eq:P-below-line}, 
\begin{equation*}
P_\beta(s)\le C\quad\forall 0\le s<M_P(\beta)
\end{equation*}
and
\begin{align*}
 &|(P_\beta^2)'(s)|\le C\quad\forall 0<s<M_P(\beta)\notag\\
 \Rightarrow\quad&|P_\beta(s)^2-P_{\beta}(s')^2|\le C|s-s'|\quad\forall 0<s,s'<M_P(\beta)
\end{align*}
for some constant $C>0$. Hence
\begin{align*}
&c_1:=\lim_{s\nearrow M_P(\beta)}P_{\beta}(s)^2\quad\mbox{ exists}\notag\\
\Rightarrow\quad&\sqrt{c_1}:=\lim_{s\nearrow M_P(\beta)}P_{\beta}(s)\quad\mbox{ exists}.
\end{align*}
We claim that $c_1=0$. Suppose $c_1>0$. Then there exists $s_1\in (0,M_P(\beta))$ such that
\begin{equation*}
P_{\beta}(s)\ge\sqrt{c_1}/2\quad\forall s_1\le s<\lim_{t\nearrow \tau_P(\beta)}W^o_{\beta}(t).
\end{equation*}
where $\tau_P(\beta)$ is given by \eqref{tau-P-definition}.
Let $t_1$ be such that $W^o_{\beta}(t_1)=s_1$. Then
\begin{align}
&(W^o_{\beta})'(t)=P_{\beta}(W^o_{\beta}(t))\ge\sqrt{c_1}/2\quad\forall t_1\le t<\tau_P(\beta)\label{Wo-beta'-uniform-lower-bd}\\
\Rightarrow\quad &M_P(\beta)>W^o_{\beta}(t)\ge W^o_{\beta}(t_1)+\frac{\sqrt{c_1}}{2}(t-t_1)\quad\forall t_1\le t<\tau_P(\beta).\label{W1'-uniform-lower-bd}
\end{align}
If $\tau_P(\beta)=\infty$, we can let $t\to\infty$ in \eqref{W1'-uniform-lower-bd}. Then the right-hand side of \eqref{W1'-uniform-lower-bd} goes to infinity and a contradiction follows. Hence $\tau_P(\beta)<\infty$. Letting $t\nearrow \tau_P(\beta)$ in \eqref{Wo-beta'-uniform-lower-bd},
we get
\begin{equation*}
(W^o_{\beta})'(\tau_P(\beta))\ge\sqrt{c_1}/2.
\end{equation*}
Since $M_P(\beta)>0$ and $(W^o_{\beta})'(\tau_P(\beta))>0$, the regular ODE continuation theorem extends the orbit beyond $\tau_P(\beta)$ if necessary. By continuity of $(W^o_{\beta})'$, there exists $t_2>\tau_P(\beta)$ such that
\[
 (W^o_{\beta})'(t)\ge\sqrt{c_1}/4>0
 \qquad\forall \tau_P(\beta)\le t\le t_2.
\]
This contradicts the definition of $\tau_P(\beta)$ as the supremum of the times up to which $(W^o_{\beta})'$ remains positive. Hence $c_1=0$. 
Then $P_\beta$ can be extended to a continuous function on $[0,M_P(\beta)]$ by setting 
\begin{equation}\label{P-beta-right-hand-limit}
P_{\beta}(M_P(\beta))=0.
\end{equation}
By \eqref{P-beta-right-hand-limit},
\begin{equation}\label{W1'-infty-limit-zero}
(W^o_{\beta})'(t)=P_{\beta}(W^o_{\beta}(t))\to 0\quad\mbox{ as }t\nearrow \tau_P(\beta).
\end{equation}
Suppose $M_P(\beta)<s_0$. Suppose $\tau_P(\beta)<\infty$. Then by \eqref{W1'-infty-limit-zero},
\begin{align}
&(W^o_{\beta})'(\tau_P(\beta))=P_{\beta}(W^o_{\beta}(\tau_P(\beta)))=0\label{W0-beta'-infty-zero}\\
\Rightarrow\quad&(W^o_{\beta})''(\tau_P(\beta))\le 0.\label{W0-beta''-negative}
\end{align}
Putting $W=W^o_\beta$ and  $t=\tau_P(\beta)$ in \eqref{eq:lienard},  by \eqref{W0-beta'-infty-zero},
\begin{equation*}
(W^o_\beta)''(\tau_P(\beta))=-h(M_P(\beta))>0
\end{equation*}
which contradicts \eqref{W0-beta''-negative}. Hence $\tau_P(\beta)=\infty$.
Putting $W=W^o_\beta$ and letting $t\to\infty$ in \eqref{eq:lienard}, we get
\begin{equation*}
\lim_{t\to\infty}(W^o_\beta)''(t)=-h(M_P(\beta))=:2c_2>0.
\end{equation*}
Thus there exists a constant $t_2$ such that
\begin{equation*}
(W^o_\beta)''(t)>c_2\quad\forall t\ge t_2\quad\Rightarrow\quad P_{\beta}(W^o_\beta(t))=(W^o_\beta)'(t)>(W^o_\beta)'(t_2)+c_2(t-t_2)\to\infty
\mbox{ as }t\to\infty 
\end{equation*}
which contradicts the fact that $P_{\beta}(s)\le as\le aM_P(\beta)$ for any $s=W^o_\beta(t)\in (0,M_P(\beta))$ by \eqref{eq:P-below-line}.
Hence $M_P(\beta)\ge s_0$ and the left-hand side inequality of (i) holds. 

\noindent $\text{\bf Proof of (ii)}$: We first prove the right-hand side inequality of (ii). Suppose there exists $\beta\in [\beta_-,\beta_+)$ such that $M_Q(\beta)=\infty$. Then by \eqref{eq:Q-equation},
\begin{equation}\label{eq:Q-square}
 (Q_\beta^2)'(s)=2(b-a+\beta s^\sigma)Q_\beta(s)+2abs-2Es^{1+\sigma}\quad\forall s> 0.
\end{equation}
We now divide the proof into two cases.

\noindent $\underline{\text{\bf Case 1}}$: \(\beta_-\le\beta\le 0\).

By \eqref{eq:Q-below-line} and \eqref{eq:Q-square},
\begin{equation}\label{Q2-ineqn2}
 (Q_\beta^2)'(s)\leq C_2 s-2Es^{1+\sigma}=s(C_2 -2Es^{\sigma})
 \quad\forall s> 0
\end{equation}
where $C_2=2b^2+2ab$.

\noindent $\underline{\text{\bf Case 2}}$:  \(0<\beta<\beta_+\).

\noindent By \eqref{eq:beta-endpoints}, \eqref{eq:Q-below-line} and \eqref{eq:Q-square},
\begin{equation}\label{Q2-ineqn3}
(Q_\beta^2)'(s)
 \leq C_2s-2(E-b\beta)s^{1+\sigma}=s(C_2-2b(\beta_+-\beta)s^{\sigma})\quad\forall s> 0.
\end{equation}
Let 
\begin{equation}\label{s2-defn}
s_2=\max\left(1,(C_2/E)^{1/\sigma},\left(C_2/[b(\beta_+-\beta)]\right)^{1/\sigma}\right).
\end{equation}
 Then by \eqref{Q2-ineqn2}, \eqref{Q2-ineqn3} and \eqref{s2-defn},
\begin{equation*}
 (Q_\beta^2)'(s)\le -C_2\quad\forall s\ge s_2
\end{equation*}
which implies that $Q_\beta(s)^2<0$ for sufficiently large $s$ and a contradiction follows. 
Hence the right-hand inequality of (ii) holds.

By an argument similar to the  proof of (i) the left-hand side inequality of (ii) follows. 
Moreover $Q_\beta$ can be extended to a continuous function on $[0,M_Q(\beta)]$ by setting 
\begin{equation}\label{Q-beta-right-hand-limit}
Q_{\beta}(M_Q(\beta))=0.
\end{equation}

\noindent $\text{\bf Proof of (iii)}$: Suppose  \(M_P(\beta)>s_0\) for some $\beta\in (\beta_-,\beta_+]$. Then by  \eqref{eq:h-sign-s0}, \eqref{eq:P-square} and \eqref{P-beta-right-hand-limit},
\begin{equation*}
(P_\beta^2)'(M_P(\beta))=-2h(M_P(\beta))<0.
\end{equation*}
This together with the Taylor formula then gives (iii). 

Then by  \eqref{eq:h-sign-s0}, \eqref{eq:Q-square}, \eqref{Q-beta-right-hand-limit} and  a similar argument we get (iv).

\noindent $\text{\bf Proof of (v)}$: Suppose $M_P(\beta)=s_0$. If $\tau_P(\beta)<\infty$, then
\begin{equation*}
W^o_\beta(\tau_P(\beta))=s_0.
\end{equation*}
Let $t_1<\tau_P(\beta)$ and 
\begin{equation*}
F(t)=\mbox{exp}\,\left(\int_{t_1}^tc_{\beta}(W^o_\beta(z))\,dz\right)\quad\forall t\le \tau_P(\beta).
\end{equation*}
Multiplying \eqref{eq:lienard} with $W=W^o_\beta$ by $F(t)$ and integrating over $(t_1,\tau_P(\beta))$,
\begin{equation}\label{W1'-eqn1}
F(\tau_P(\beta))(W^o_\beta)'(\tau_P(\beta))=(W^o_\beta)'(t_1)-\int_{t_1}^{\tau_P(\beta)}h(W^o_\beta(t))F(t)\,dt
\end{equation}
Since $W^o_\beta(t)<s_0$ for any $t<\tau_P(\beta)$, $h(W^o_\beta(t))<0$ for any $t<\tau_P(\beta)$. Hence
\begin{equation}\label{h-F-integral}
\int_{t_1}^{\tau_P(\beta)}h(W^o_\beta(t))F(t)\,dt<0.
\end{equation}
Since $(W^o_\beta)'(t_1)>0$, by \eqref{W1'-eqn1} and \eqref{h-F-integral},
\begin{equation*}
F(\tau_P)(W^o_\beta)'(\tau_P(\beta))>0\quad\Rightarrow\quad P_{\beta}(s_0)=P_{\beta}(W^o_\beta(\tau_P(\beta)))=(W^o_\beta)'(\tau_P(\beta))>0
\end{equation*}
which contradicts \eqref{P-beta-right-hand-limit}. Hence $\tau_P(\beta)=\infty$ and (v) follows.
By a similar argument we get (vi) and the lemma follows.

\end{proof}

\begin{lem}[Sign of $c_\beta(s)$ at a degenerate endpoint]\label{lem:degenerate-sign}
For any $\beta\in\mathbb{R}$ let $c_{\beta}$ be given by \eqref{eq:c-beta}. Then the following holds.

\begin{enumerate}[(i)]

\item Let $\beta\in (\beta_-,\beta_+]$ be such that \(M_P(\beta)=s_0\). Then either
\begin{equation}
 c_\beta(s_0)>0,
 \qquad\hbox{or}\qquad
 c_\beta(s_0)=0\ \hbox{ and }\ \beta\leq0.
 \label{eq:P-degenerate-sign}
\end{equation}

\item Let $\beta\in [\beta_-,\beta_+)$ be such that \(M_Q(\beta)=s_0\). Then either
\begin{equation}
 c_\beta(s_0)<0,
 \qquad\hbox{or}\qquad
 c_\beta(s_0)=0\ \hbox{ and }\ \beta\geq0.
 \label{eq:Q-degenerate-sign}
\end{equation}

\item There is no $\beta\in (\beta_-,\beta_+)$  satisfying \(M_P(\beta)=M_Q(\beta)=s_0\).

\end{enumerate}
\end{lem}
\begin{proof}
Let
\begin{equation*}
 \calE(t):=\frac12(W^o_{\beta})'(t)^2+H(W^o_{\beta}(t))\quad\forall t<\tau_P(\beta).
\end{equation*}
where $H$ is given by \eqref{eq:h-H}. Since $W^o_{\beta}$ satisfies  \eqref{eq:lienard} in $(-\infty,\tau_P(\beta))$,  
\begin{equation}
 \calE'(t)=-c_\beta(W^o_{\beta}(t))(W^o_{\beta})'(t)^2\quad\forall t<\tau_P(\beta).
 \label{eq:energy-derivative}
\end{equation}
 We next observe that 
\begin{equation}\label{eq:W1<s0}
0<W^o_{\beta}(t)<s_0,\quad (W^o_{\beta})'(t)>0\quad\forall t<\tau_P(\beta).
\end{equation}
Since 
\begin{equation*}
H'(s)=h(s)\quad\forall s\ge 0,
\end{equation*}
by \eqref{eq:h-sign-s0} \(H\) is strictly decreasing on \((0,s_0)\). Hence by \eqref{eq:W1<s0},
\begin{equation}\label{eq:E1-lower-bd}
\calE(t)=\frac12(W^o_{\beta})'(t)^2+H(W^o_{\beta}(t))>H(s_0)\quad\forall t<\tau_P(\beta).
\end{equation}
Since 
\begin{equation*}
W^o_{\beta}(t)\to s_0\quad\mbox{ and }\quad (W^o_{\beta})'(t)=P_{\beta}(W^o_{\beta}(t))\to P_{\beta}(s_0)=0\quad\mbox{ as }
t\nearrow\tau_P(\beta), 
\end{equation*}
we have
\begin{equation}\label{E1-limit}
 \calE(t)\longrightarrow H(s_0)\quad\mbox{ as }t\nearrow\tau_P(\beta).
\end{equation}
\noindent $\text{\bf Proof of (i)}$:
Suppose  
\begin{equation}\label{eq:c-beta-s0-negative}
c_\beta(s_0)<0
\end{equation}
holds.  By continuity of
\(c_\beta\), there exists a constant \(\delta>0\)  such that
\begin{equation}\label{eq:c-beta-negative}
c_\beta(s)<0\quad\forall s_0-\delta<s\leq s_0.
\end{equation}
By \eqref{eq:W1<s0} there exists a constant \(T_1\) such that
\(W^o_{\beta}(t)\in(s_0-\delta,s_0)\) for any \(t\geq T_1\).  Hence by 
\eqref{eq:energy-derivative}, \eqref{eq:E1-lower-bd}, \eqref{E1-limit} and \eqref{eq:c-beta-negative},
\begin{align*}
&\calE'(t)=-c_\beta(W^o_{\beta}(t))(W^o_{\beta})'(t)^2>0\quad\forall t\geq T_1\notag\\
\Rightarrow\quad&\calE(t)\ge\calE(T_1)>H(s_0)\quad\forall t\geq T_1\notag\\
\Rightarrow\quad&H(s_0)\ge\calE(T_1)>H(s_0)\quad\mbox{ as }t\nearrow\tau_P(\beta)
\end{align*}
and a contradiction follows. Thus \eqref{eq:c-beta-s0-negative} does not hold.

Suppose now
\begin{equation}\label{eq:c-beta-s0=0-case}
c_\beta(s_0)=0\quad\mbox{ and }\quad\beta>0.
\end{equation}
Then
\[
 c_\beta(s)=c_\beta(s)-c_\beta(s_0)=\beta\bigl(s^\sigma-s_0^\sigma\bigr)<0\quad\forall 0<s<s_0.
\]
Hence we may repeat the previous argument to conclude that the case \eqref{eq:c-beta-s0=0-case}
is also not possible. Hence \eqref{eq:P-degenerate-sign} holds.

By an argument similar to the proof of (i) we get (ii).

\noindent $\text{\bf Proof of (iii)}$:
Suppose that there exists  $\beta\in (\beta_-,\beta_+)$  such that  \(M_P(\beta)=M_Q(\beta)=s_0\). Then by (i) and (ii) we have
\begin{align*}
&c_\beta(s_0)=0\quad\mbox{ and }\quad\beta=0\notag\\
\Rightarrow\quad&b-a=0\notag\\
\Rightarrow\quad&\calE'(t)=0\quad\forall t<\tau_P(\beta)\notag\\
\Rightarrow\quad&\calE(t_1)=\calE(t_2)\quad\forall t_1,t_2\le\tau_P(\beta)\notag\\
\Rightarrow\quad&0=H(s_0)\quad\mbox{ as }t_1\to-\infty, t_2\to\tau_P(\beta).
\end{align*}
Since $H(s_0)<0$, a contradiction follows. Thus no such $\beta$ exists and (iii) follows.

\end{proof}

\begin{lem}[Compact dependence of the slope half-branches]\label{lem:slope-compact-dependence}
Let $\beta _0\in[\beta_-,\beta_+]$.  
\begin{enumerate}[(i)]
\item For any $0<R<M_P(\beta _0)$, there exists a constant $\delta_2>0$ such that for any $\beta\in I_{\delta_2}:=[\beta_0-\delta_2,\beta_0+\delta_2]\cap [\beta_-,\beta_+]$ we have $R<M_P(\beta)$ and
\begin{equation}\label{P-compact-convergence}
 P_\beta\to P_{\beta _0}
 \quad\text{uniformly on }[0,R]\quad\mbox{ as }\beta\to\beta_0.
\end{equation}
The function $M_P(\beta)$ is lower semicontinuous on $[\beta_-,\beta_+]$. That is 
\begin{equation}\label{eq:MP-lower-s-cont}
\liminf_{\beta\to\beta_0}M_P(\beta)\ge M_P(\beta_0)\quad\forall\beta_0\in [\beta_-,\beta_+].
\end{equation}

\item For any $0<R<M_Q(\beta _0)$, there exists a constant $\delta_2>0$ such that for any $\beta\in I_{\delta_2}$ we have $R<M_Q(\beta)$ and
\begin{equation}\label{Q-compact-convergence}
 Q_\beta\to Q_{\beta _0}
 \quad\text{uniformly on }[0,R]\quad\mbox{ as }\beta\to\beta_0.
\end{equation}
 The function $M_Q(\beta)$ is lower semicontinuous on $[\beta_-,\beta_+]$. That is 
\begin{equation}\label{eq:MQ-lower-s-cont}
\liminf_{\beta\to\beta_0}M_Q(\beta)\ge M_Q(\beta_0)\quad\forall\beta_0\in [\beta_-,\beta_+].
\end{equation}

\end{enumerate}

\end{lem}

\begin{proof}
Let \(0<R<M_P(\beta_0)\).  Then \(P_{\beta_0}(s)>0\) for \(0<s\leq R\). By Theorem \ref{thm-existence-self-similar1} and Lemma \ref{lem:f-beta-parameter-limit} 
there exists a constant $\eps>0$ such that for any $\beta\in[\beta_-,\beta_+]$ there exists a unique maximal positive solution $f^o_{\beta}\in C^1([0,R_o(\beta,1)))\cap C^2((0,R_o(\beta,1)))$ of \eqref{f-ode} with $\alpha=\frac{2\beta+1}{1-m}$, $3\eps<R_o(\beta,1)$, in $(0,R_o(\beta,1))$ which satisfies \eqref{f-fr-origin-value}  and 
\eqref{eq:f-beta-parameter-limit} with $\eps$ being replaced by $2\eps$. Let $W^o_{\beta}(t)$ be given by \eqref{eq:W-origin-infinity}. Then  \eqref{eq:w0-beta-eqn} holds.
By replacing $\eps>0$ with a smaller value if necessary we may assume without loss of generality that 
\begin{equation*}
W^o_{\beta_0}(\log(2\eps))<R.
\end{equation*}
By Lemma \ref{lem:linear-comparison} \eqref{P-Q-zero-slopes} holds. Hence there exists a constant $t_0<\log 2\eps$ such that
\begin{equation}\label{W'-positive10}
\frac{(W^o_{\beta})'(t)}{W^o_{\beta}(t)}>\frac{a}{2}>0\quad\forall t\le t_0,\beta\in [\beta_-,\beta_+].
\end{equation}
On the other hand by \eqref{eq:f-beta-parameter-limit}, \eqref{eq:W-origin-infinity} and \eqref{eq:w0-beta-eqn},
\begin{equation}\label{eq:W-beta-beta0-limit}
W^o_{\beta}(t)\to W^o_{\beta_0}(t)\quad\mbox{and}\quad (W^o_{\beta})'(t)\to (W^o_{\beta_0})'(t)\quad\mbox{uniformly on $(-\infty,\log(2\eps)]$ as $\beta\to\beta_0$}. 
\end{equation}
Hence  there exists a constant $\delta_2>0$ such that
\begin{equation}\label{eq:W'-beta-uniformly-positive}
(W^o_{\beta})'(t)>0\quad\forall t_0\le t\le\log(2\eps)\quad\mbox{ and }\quad W^o_{\beta_0}(\log\eps)<W^o_{\beta}(\log(2\eps))<(W^o_{\beta_0}(\log(2\eps))+R)/2
\end{equation}
for any $\beta\in I_{\delta_2}$. By \eqref{W'-positive10} and \eqref{eq:W'-beta-uniformly-positive},
\begin{equation}\label{eq:W'-beta-uniformly-positive11}
(W^o_{\beta})'(t)>0\quad\forall  t\le\log(2\eps),\beta\in I_{\delta_2}.
\end{equation}
Let $s_{\eps}:=W^o_{\beta_0}(\log\eps)$.
Then the map $W^o_{\beta_0}:(-\infty,\log\eps]\to (0,s_{\eps}]$ has an inverse $(W^o_{\beta_0})^{-1}$. We will now assume that 
$\beta\in I_{\delta_2}$. Then the map $W^o_{\beta}:(-\infty,\log(2\eps)]\to(0,W^o_{\beta}(\log(2\eps))]$ has an inverse $(W^o_{\beta})^{-1}$ on $(0,W^o_{\beta}(\log(2\eps))]$. Moreover by \eqref{eq:W'-beta-uniformly-positive},
\begin{equation}\label{s-variepsilon-bd}
0<s_{\eps}:=W^o_{\beta_0}(\log\eps)<W^o_{\beta}(\log(2\eps))<M_P(\beta)\quad\forall\beta\in I_{\delta_2}.
\end{equation} 
Since $P_{\beta}$ satisfies \eqref{eq:P-below-line} and \eqref{eq:P-square} in $(0,M_P(\beta))$, by \eqref{eq:W'-beta-uniformly-positive} we have
\begin{align*}
|(P_{\beta}^2)'(s)|&\le C\quad\mbox{ in }(0,W^o_{\beta}(\log(2\eps))]\quad\forall\beta\in I_{\delta_2}
\end{align*}
for some constant $C>0$. Hence
\begin{equation}\label{eq:P-beta-lip}
|P_{\beta}(s)^2-P_{\beta}(s')^2|\le C|s-s'|\quad\forall s,s'\in [0,W^o_{\beta}(\log(2\eps))]\quad\forall\beta\in I_{\delta_2}.
\end{equation}
For any $0<s\le s_{\eps}$, let $t=(W^o_{\beta_0})^{-1}(s)$. Then $t\le\log\eps$ and by \eqref{eq:W-beta-beta0-limit},
\eqref{eq:W'-beta-uniformly-positive}, \eqref{s-variepsilon-bd} and \eqref{eq:P-beta-lip},
\begin{align}\label{P-beta-beta0-limit2}
|P_{\beta}(s)^2-P_{\beta_0}(s)^2|=&|P_{\beta}(W^o_{\beta_0}(t))^2-P_{\beta_0}(W^o_{\beta_0}(t))^2|\notag\\
\le&|P_{\beta}(W^o_{\beta_0}(t))^2-P_{\beta}(W^o_{\beta}(t))^2|+|P_{\beta}(W^o_{\beta}(t))^2-P_{\beta_0}(W^o_{\beta_0}(t))^2|\notag\\
\le&C|W^o_{\beta_0}(t)-W^o_{\beta}(t)|+|(W^o_{\beta})'(t)^2-(W^o_{\beta_0})'(t)^2|\quad\forall s=W^o_{\beta_0}(t)\in  [0,s_{\eps}],\beta\in I_{\delta_2}\notag\\
\to&0\quad\mbox{ uniformly on $[0,s_{\eps}]$ as $\beta\to\beta_0$}.
\end{align}
By \eqref{P-beta-beta0-limit2},
\begin{equation}\label{P-beta-beta0-pt-limit}
P_{\beta}(s_{\eps})\to P_{\beta_0}(s_{\eps})\quad\mbox{ as }\beta\to\beta_0.
\end{equation}  
We now claim that $R<M_P(\beta)$ for $\beta$ sufficiently close to $\beta_0$ and
\begin{equation}\label{P-beta-beta0-uniform-limit}
P_{\beta}(s)\to P_{\beta_0}(s)\quad\mbox{ uniformly on $[s_{\eps},R]$ as }\beta\to\beta_0.
\end{equation}
In order to prove this claim we let
\[
 m_R:=\min_{s_{\eps}\leq s\leq R}P_{\beta_0}(s).
\]
Then $m_R>0$. Let
\[
 \mathcal K:=
 \left\{(s,p,\beta):
 s_{\eps}\leq s\leq R,
 \abs{p-P_{\beta_0}(s)}\leq m_R/2,
 \abs{\beta-\beta_0}\leq1
 \right\}.
\]
Then \(p\geq m_R/2\) for any $(s,p,\beta)\in\mathcal K$.  Since $P_{\beta}$ satisfies \eqref{eq:P-equation},
\begin{equation}\label{eq:P-equation2}
P_{\beta}'(s)=-\bigl(b-a+\beta s^\sigma\bigr)-\frac{h(s)}{P_{\beta}(s)}=-c_{\beta}(s)-\frac{h(s)}{P_{\beta}(s)}=:F(s,P_{\beta}(s),\beta)\quad\mbox{ in }(0,M_P(\beta)).
\end{equation}  Hence \(F\) is continuous on \(\mathcal K\) and 
by the mean-value theorem,
\[
 \abs{F(s,p,\beta)-F(s,q,\beta)}
 \leq L_R\abs{p-q}\quad\forall (s,p,\beta),(s,q,\beta)\in\mathcal K.
\]
where
\[
 L_R:=\frac{4\max_{s_{\eps}\leq s\leq R}\abs{h(s)}}{m_R^2}.
\]
Moreover
\[
 \abs{F(s,p,\beta)-F(s,p,\beta_0)}
 =\abs{\beta-\beta_0}s^\sigma
 \leq R^\sigma\abs{\beta-\beta_0}\quad\forall (s,p,\beta),(s,p,\beta_0)\in\mathcal K.
\]
Let
\[
 e_{\beta}(s):=P_{\beta}(s)-P_{\beta_0}(s).
\]
Since \(e_{\beta}(s_{\eps})\to0\) as $\beta\to\beta_0$, by choosing $\delta_2>0$ to be sufficiently small if necessary we may assume that
\begin{equation*}
\abs{e_{\beta}(s_{\eps})}<m_R/4\quad\forall \beta\in I_{\delta_2}.
\end{equation*}  
Let
\begin{equation*}
\bar s_\beta=\sup\{s_{\eps}<s_1'\le R:\abs{e_{\beta}(s)}<m_R/2\quad\mbox{ and }\quad P_{\beta}(s)>0\quad\forall s_{\eps}\le s\le s_1'\}.
\end{equation*}
Suppose $\bar s_\beta<R$. Then
\begin{equation*}
\abs{e_{\beta}(s)}<m_R/2\quad\mbox{ and }\quad P_{\beta}(s)>0\quad\forall s_{\eps}\le s<\bar s_\beta 
\end{equation*} 
and either 
\begin{equation}\label{tube-end-pts-value}
\abs{e_{\beta}(\bar s_\beta)}=m_R/2\quad\mbox{ or }\quad P_{\beta}(\bar s_\beta)=0.
\end{equation}
By \eqref{eq:P-equation2},
\[
 \begin{aligned}
 |e_\beta(s)|
 &\leq |e_\beta(s_\eps)|
 +L_R\int_{s_\eps}^{s}|e_\beta(\tau)|\dd\tau
 +R^\sigma|\beta-\beta_0|(s-s_\eps)
 \end{aligned}
 \quad\forall s_\eps\leq s\leq\bar s_\beta.
\]
Gronwall's inequality gives
\begin{equation}\label{P-beta-beta0-ineqn}
 \sup_{s_\eps\leq s\leq\bar s_\beta}|e_\beta(s)|
 \leq
 \left(|e_\beta(s_\eps)|
 +R^\sigma(R-s_\eps)|\beta-\beta_0|\right)
 e^{L_R(R-s_\eps)}.
\end{equation}
As \(\beta\to\beta_0\), the right-hand side of \eqref{P-beta-beta0-ineqn} tends to zero because
\(e_{\beta}(s_{\eps})\to 0\) as $\beta\to\beta_0$.  Hence by choosing $\delta_2>0$ to be sufficiently small we may assume that the right hand side of \eqref{P-beta-beta0-ineqn} is strictly smaller than \(m_R/4\) for any $\beta\in I_{\delta_2}$. Hence 
\begin{equation}\label{end-pt-estimate1}
|e_\beta(s)|\le \frac{m_R}{4}\quad\mbox{ and }\quad P_{\beta}(s)\ge P_{\beta_0}(s)-\frac{m_R}{4}\ge\frac{3m_R}{4}>0\quad\forall s_\eps\leq s\leq\bar s_\beta,\beta\in I_{\delta_2}
\end{equation}
which contradicts \eqref{tube-end-pts-value}. Hence $\bar s_\beta=R$ and \eqref{P-beta-beta0-ineqn}, \eqref{end-pt-estimate1}, holds with $\bar s_\beta=R$.  
Hence $M_P(\beta)>R$ for any $\beta\in I_{\delta_2}$. Letting $\beta\to\beta_0$ in \eqref{P-beta-beta0-ineqn} we get
\eqref{P-beta-beta0-uniform-limit} holds.  By \eqref{P-beta-beta0-limit2} and \eqref{P-beta-beta0-uniform-limit} we get \eqref{P-compact-convergence}.   Finally by \eqref{end-pt-estimate1},
\begin{align*} 
&\liminf_{\beta\to\beta_0}M_P(\beta)
 \geq R\quad\forall 0<R<M_P(\beta_0)\notag\\
\Rightarrow\quad& \liminf_{\beta\to\beta_0}M_P(\beta)
 \geq M_P(\beta_0)\quad\mbox{ as }R\to M_P(\beta_0)
\end{align*}
and (i) follows. By Lemma \ref{lem:g-beta-parameter-limit} or Lemma \ref{lem:g-beta-parameter-limit2}, according to the range of $m$, together with \eqref{eq:W-infinity1}, \eqref{eq:Q-equation}, \eqref{P-Q-zero-slopes} and a similar argument we get (ii) and the lemma follows.

\end{proof}

\begin{lem}[Endpoint divergence]\label{lem:endpoint-divergence}
One has
\begin{equation}\label{eq:endpoint-divergence1}
 \lim_{\beta\searrow\beta_-}M_P(\beta)=\infty,
 \qquad
 \lim_{\beta\nearrow\beta_+}M_Q(\beta)=\infty.
\end{equation}
\end{lem}

\begin{proof}
By Lemma \ref{lem:linear-comparison}, $M_P(\beta_-)=\infty$ and $M_Q(\beta_+)=\infty$.  Fix $R>0$.  Applying Lemma \ref{lem:slope-compact-dependence}(i) at $\beta_0=\beta_-$ gives $M_P(\beta)>R$ for all $\beta$ sufficiently close to $\beta_-$ from the right.  Hence
\begin{equation}\label{MP-beta-limit6}
\liminf_{\beta\searrow\beta_-}M_P(\beta)\ge R.
\end{equation}
Since $R$ is arbitrary, $M_P(\beta)\to\infty$ as $\beta\searrow\beta_-$.  Applying part (ii) at $\beta_0=\beta_+$ gives the analogous conclusion $M_Q(\beta)\to\infty$ as $\beta\nearrow\beta_+$, and the lemma follows.

\end{proof}

\begin{lem}[Continuity of the turning amplitudes]\label{lem:turning-continuity}
The following holds.
\begin{enumerate}[(i)]

\item  The function  $M_P(\beta)$ is continuous on $(\beta_-,\beta_+]$.

\item The function $M_Q(\beta)$ is continuous on $[\beta_-,\beta_+)$.

\end{enumerate}

\end{lem}

\begin{proof}
Let $\beta_0\in (\beta_-,\beta_+]$. By \eqref{eq:MP-finite}, we can divide the proof of (i)  into two cases.

\noindent $\underline{\text{\bf Case 1}}$: $M_P(\beta_0)>s_0$.

\noindent Let
\[
 M_0:=M_P(\beta_0),
 \qquad
 Y_\beta(s):=P_\beta(s)^2.
\]
Let \(\varepsilon>0\) and 
\begin{equation}
 0<r_1<\min\left(\varepsilon,\frac{M_0-s_0}{2}\right).
 \label{eq:choose-r-nondegenerate}
\end{equation}
Then
\[
 [M_0-r_1,M_0+r_1]\subset(s_0,\infty).
\]
Then by \eqref{eq:P-square}, 
\begin{equation}\label{eq:YP-square-continuity}
 Y_\beta'(s)=-2c_\beta(s)\sqrt{Y_\beta(s)}-2h(s)\quad\mbox{ in }(0,M_P(\beta)).
\end{equation}
where $c_{\beta}$ and $h$ are given by \eqref{eq:c-beta} and \eqref{eq:h-H} respectively.
By Lemma \ref{lem:slope-compact-dependence}, $M_P(\beta)$ is lower semicontinuous at $\beta_0$. Hence in order to prove continuity of $M_P(\beta)$ at $\beta_0$
 it suffices to prove that $M_P(\beta)$ is upper semicontinuous at \(\beta=\beta_0\).
Since \(h\) is continuous and strictly positive on \((s_0,\infty)\),
there is a constant \(\gamma>0\) such that
\begin{equation}\label{eq:g-positive-near-M0}
 h(s)\geq\gamma
 \quad\forall M_0-r_1\leq s\leq M_0+r_1.
\end{equation}
Note that there exists a constant $C_1>0$ such that
\begin{equation}
 \abs{c_\beta(s)}\leq C_1
 \quad\forall \beta\in (\beta_-,\beta_+],\ M_0-r_1\leq s\leq M_0+r_1.
 \label{eq:c-uniform-bound-near-M0}
\end{equation}
Let
\begin{equation}
 p_*:=\frac{\gamma}{2(C_1+1)}.
 \label{eq:p-star-continuity}
\end{equation}
Since \(P_{\beta_0}(s)\to0\) as \(s\nearrow M_0\), there exists a constant 
$s_1\in(M_0-r_1/2,M_0)$ such that
\begin{equation}
 P_{\beta_0}(s_1)<\frac{p_*}{2},
 \qquad
 P_{\beta_0}(s_1)^2<\frac{\gamma r_1}{8}.
 \label{eq:P0-small-at-s1}
\end{equation}
By \eqref{P-compact-convergence},
\begin{equation}\label{eq:P-beta-beta0-limit5}
P_\beta(s_1)\longrightarrow P_{\beta_0}(s_1)
 \quad\mbox{ as }\beta\to\beta_0.
\end{equation}
By \eqref{P-compact-convergence}, \eqref{eq:P0-small-at-s1} and \eqref{eq:P-beta-beta0-limit5} there exists a constant $\delta_2>0$  such  that  for any $\beta\in I_{\delta_2}:=[\beta_0-\delta_2,\beta_0+\delta_2]\cap (\beta_-,\beta_+]$, 
\begin{equation*}
P_{\beta}(s)>0\quad\forall 0<s\le\max (s_1,M_0-(r_1/4))\quad\Rightarrow\quad M_P(\beta)>\max (s_1,M_0-(r_1/4))
\end{equation*}
and
\begin{equation}
 P_\beta(s_1)<p_*,
 \qquad
 P_\beta(s_1)^2<\frac{\gamma r_1}{4}.
 \label{eq:Pbeta-small-at-s1}
\end{equation}
We will now assume $\beta\in I_{\delta_2}$ for the rest of the proof.

We claim that $M_P(\beta)\leq M_0+r_1$.  Suppose the claim does not hold. Then $M_P(\beta)>M_0+r_1$.  Hence $P_\beta$ is positive on $[s_1,M_0+r_1]$.  Let
\[
 s_2:=\sup\{s'\in(s_1,M_0+r_1]:
 P_\beta(s)<p_*\text{ for every }s\in[s_1,s')\}.
\]
If $s_2<M_0+r_1$, then $P_\beta(s_2)=p_*$ and therefore $Y_\beta'(s_2)\geq0$.  On the other hand by \eqref{eq:YP-square-continuity}, \eqref{eq:g-positive-near-M0}, \eqref{eq:c-uniform-bound-near-M0} and \eqref{eq:p-star-continuity},
\begin{equation}\label{eq:Y-strictly-decreasing}
 Y_\beta'(s)\leq2C_1p_*-2\gamma\leq-\gamma
 \qquad\forall s_1\leq s\leq s_2,
\end{equation}
which is a contradiction at $s=s_2$.  Hence $s_2=M_0+r_1$, and \eqref{eq:Y-strictly-decreasing} holds throughout $[s_1,M_0+r_1]$.  Set
\[
 \bar s:=s_1+\frac{Y_\beta(s_1)}{\gamma}.
\]
By \eqref{eq:Pbeta-small-at-s1}, $\bar s<s_1+r_1/4<M_0+r_1$.  Integrating \eqref{eq:Y-strictly-decreasing} from $s_1$ to $\bar s$ yields
\[
 Y_\beta(\bar s)\leq Y_\beta(s_1)-\gamma(\bar s-s_1)=0,
\]
contradicting the strict positivity of $P_\beta$ on this interval.  Thus $M_P(\beta)\leq M_0+r_1<M_0+\varepsilon$, and
 \begin{align}\label{eq:MP-upper-semicontinuity}
&M_P(\beta)<M_0+\varepsilon\quad\forall \beta\in I_{\delta_2}\notag\\
 \Rightarrow\quad&\limsup_{\beta\to\beta_0}M_P(\beta)\leq M_0+\varepsilon\notag\\
 \Rightarrow\quad&\limsup_{\beta\to\beta_0}M_P(\beta)\leq M_0\quad\mbox{ as }\eps\to 0.
\end{align}
This together with \eqref{eq:MP-lower-s-cont} of Lemma \ref{lem:slope-compact-dependence} implies that 
\[
 M_P(\beta)\longrightarrow M_P(\beta_0)
 \quad\mbox{ as }\beta\to\beta_0.
\]

\noindent $\underline{\text{\bf Case 2}}$: $M_P(\beta_0)=s_0$.

\noindent Suppose \(M_P(\beta)\) is not continuous at $\beta_0$.
Since $M_P(\beta)\ge s_0$ for any $\beta\in (\beta_-,\beta_+]$, there exist a constant $\delta>0$ and a sequence $\{\beta_j\}_{j=1}^{\infty}\subset (\beta_-,\beta_+]$, $\beta_j\to\beta_0$ as $j\to\infty$, such that
\begin{equation}\label{M-P-beta-sequence-lower-bd}
 M_P(\beta_j)\geq s_0+2\delta\quad\forall j\in\mathbb{Z}^+.
\end{equation}
Let
\[
 S:=s_0+\delta.
\]
Then \(S<M_P(\beta_j)\) for any $j\in\mathbb{Z}^+$.  Hence every branch
\(P_{\beta_j}\) is defined and strictly positive on the same interval
\((0,S]\).  Let
\[
 Y_j(s):=P_{\beta_j}(s)^2\quad\quad\forall 0\le s\leq S.
\]
By Lemma \ref{lem:linear-comparison},
\begin{align}\label{eq:Yj-upper-bd}
&0<P_{\beta_j}(s)<as\quad\forall 0<s\leq S,j\in\mathbb{Z}^+\notag\\
\Rightarrow\quad&0<Y_j(s)\leq a^2s^2\leq a^2S^2
 \quad\forall 0<s\leq S,j\in\mathbb{Z}^+.
\end{align}
By \eqref{eq:YP-square-continuity} and \eqref{eq:Yj-upper-bd},
\begin{align*}
&|Y_j'(s)|=\left|-2c_{\beta_j}(s)\sqrt{Y_j(s)}-2h(s)\right|\le C_2
\quad\forall 0<s\leq S,j\in\mathbb{Z}^+\notag\\
\Rightarrow\quad& \abs{Y_j(s)-Y_j(r)}\leq C_2\abs{s-r}\quad\forall s,r\in [0,S],j\in\mathbb{Z}^+
\end{align*}
where $C_2>0$ is some constant.
By the Arzel\`a--Ascoli theorem, the sequence $\{Y_j\}_{j=1}^{\infty}$ has a subsequence which we may assume without loss of generality to be the sequence $\{Y_j\}_{j=1}^{\infty}$ itself, converges uniformly on $[0,S]$ to some  function \(Y\in C([0,S])\) as $j\to\infty$. Then
\begin{equation}\label{eq:Y-nonnegative}
Y(s)\geq0\qquad\forall 0\leq s\leq S
\end{equation}
and
\[
 \abs{Y(s)-Y(r)}\leq C_2\abs{s-r}
 \qquad\forall 0\leq r,s\leq S.
\]
By Lemma \ref{lem:slope-compact-dependence} \eqref{P-compact-convergence} holds for any $0<R<M_P(\beta_0)=s_0$.
Hence
\begin{align}\label{eq:Yj-limit}
&Y_j=P_{\beta_j}^2\longrightarrow P_{\beta_0}^2
 \qquad\text{uniformly on }[0,R]\quad\mbox{ as }j\to\infty\notag\\
 \Rightarrow\quad&Y(s)=P_{\beta_0}(s)^2\quad\forall 0\leq s\leq R
\end{align}
for any $0<R<M_P(\beta_0)=s_0$. Letting $R\nearrow M_P(\beta_0)$ in \eqref{eq:Yj-limit},
\begin{equation*}
Y(s)=P_{\beta_0}(s)^2\quad\forall 0\leq s<s_0.
\end{equation*}
Since $Y$ is continuous and $P_{\beta_0}(s)\to0$ as $s\nearrow s_0$, we have
\begin{equation}\label{eq:Y-s0-zero-continuity}
Y(s_0)=0.
\end{equation}
By \eqref{eq:YP-square-continuity},
\begin{align}\label{eq:Yj-eqn2}
&Y_j'(s)=-2c_{\beta_j}(s)\sqrt{Y_j(s)}-2h(s)\quad\mbox{ in }(0,S]\quad\forall j\in\mathbb{Z}^+\notag\\
\Rightarrow\quad& Y_j(s)-Y_j(s_0)=\int_{s_0}^{s}
 \left[-2c_{\beta_j}(\tau)\sqrt{Y_j(\tau)}-2h(\tau)\right]\dd\tau\quad\forall s_0<s\le S.
\end{align}
Since
\[
 \left|\sqrt{Y_j(s)}-\sqrt{Y(s)}\right|\le\sqrt{\left|Y_j(s)-Y(s)\right|}\to 0
 \quad\text{ uniformly on }[s_0,S]\quad\mbox{ as }j\to\infty,
\]
letting $j\to\infty$ in \eqref{eq:Yj-eqn2}, by \eqref{eq:Y-s0-zero-continuity} we get
\begin{equation}\label{eq:Y-limit-P}
 Y(s)=\int_{s_0}^s\left[-2c_{\beta_0}(\tau)\sqrt{Y(\tau)}-2h(\tau)\right]\dd\tau
 \quad\forall s_0\leq s\leq S.
\end{equation}
By Lemma \ref{lem:degenerate-sign}, \eqref{eq:P-degenerate-sign} holds. First suppose that
\begin{equation}\label{eq:c-beta-s0-positive}
c_{\beta_0}(s_0)>0.
\end{equation}
By continuity there exists a constant $s_1\in (s_0,S)$ such that
\begin{equation}
c_{\beta_0}(\tau)>0\quad\forall s_0\le\tau\le s_1.
\end{equation} 
  Since \(h(\tau)>0\) for any \(\tau>s_0\), the integrand in
\eqref{eq:Y-limit-P} is strictly negative for $s_0<s\le s_1$.  Consequently \(Y(s)<0\) for any $s_0<s\le s_1$ which contradicts  \eqref{eq:Y-nonnegative}. Hence the case \eqref{eq:c-beta-s0-positive} is not possible. 

It remains to consider the case
\begin{equation}\label{eq:c-beta0=0-beta-negative}
 c_{\beta_0}(s_0)=0,\qquad \beta_0\leq0.
\end{equation}
Since \(c_{\beta_0}\) is continuously differentiable near  \(s_0\), by the mean-value theorem there exists a constant \(C_0\geq0\) such that
\[
 \abs{c_{\beta_0}(s)}\leq C_0(s-s_0)
\]
for \(s> s_0\) close to \(s_0\).  Furthermore,
\[
 h(s_0)=0,
 \qquad
 h'(s_0)=-ab+(1+\sigma)Es_0^\sigma=\sigma ab>0,
\]
because \(Es_0^\sigma=ab\).  Therefore
\[
 h(s)\geq\frac{\sigma ab}{2}(s-s_0).
\]
for \(s> s_0\) close to \(s_0\).
Since \(Y\) is \(C_2\)-Lipschitz and \(Y(s_0)=0\),
\[
 0\leq Y(s)\leq C_2(s-s_0)\quad\mbox{ and }\quad
 \sqrt{Y(s)}\leq\sqrt{C_2}\,(s-s_0)^{1/2}\quad\forall s_0<s\le S.
\]
Thus
\begin{align*}
-2c_{\beta_0}(s)\sqrt{Y(s)}-2h(s)\leq\left(2C_0\sqrt{C_2}\,(s-s_0)^{1/2}
 -\sigma ab\right)(s-s_0)<0
\end{align*}
for \(s>s_0\) sufficiently close to \(s_0\).  Equation
\eqref{eq:Y-limit-P} again gives \(Y(s)<0\) for \(s>s_0\) sufficiently close to \(s_0\), contradicting
\eqref{eq:Y-nonnegative}. Hence the case \eqref{eq:c-beta0=0-beta-negative} is also not possible.
 Therefore  there does not exist any  sequence $\{\beta_j\}_{j=1}^{\infty}\subset (\beta_-,\beta_+]$, $\beta_j\to\beta_0$ as $j\to\infty$, that satisfies \eqref{M-P-beta-sequence-lower-bd}.
We conclude that
\[
 M_P(\beta)\longrightarrow s_0
        \quad\text{as }\beta\to\beta_0.
\]
By case 1 and case 2, we get (i). By \eqref{eq:MQ-finite} and a similar argument we get (ii) and the lemma follows.
\end{proof}

\section{Existence of solutions with fixed decay rate at infinity}
\setcounter{equation}{0}
\setcounter{thm}{0}

In this section we will use a matching argument to prove Theorem \ref{thm-existence-self-similar3}. We will prove Theorem \ref{thm-existence-self-similar3} for the case $\rho_1=1$ first. The general case of $\rho_1>0$ will then follow by a rescaling argument.

\begin{prop}[Matching of the first turning levels]\label{prop:eigenvalue}
There exists a constant  $\widehat\beta\in(\beta_-,\beta_+)$
such that
\begin{equation}
 M_P(\widehat\beta)=M_Q(\widehat\beta)=:s_*>s_0.
\label{eq:matched-M}
\end{equation}
\end{prop}
\begin{proof}
By (ii) of Lemma \ref{lem:turning-finite}, \eqref{eq:endpoint-divergence1} of Lemma \ref{lem:endpoint-divergence} and (ii) of Lemma \ref{lem:turning-continuity}, there exists a constant $\beta_1\in (\beta_-,\beta_+)$ sufficiently close to $\beta_-$ such that
\begin{equation}\label{eq:left-order}
 M_P(\beta_1)>M_Q(\beta_1).
\end{equation}
Similarly by (i) of Lemma \ref{lem:turning-finite}, \eqref{eq:endpoint-divergence1} of Lemma \ref{lem:endpoint-divergence} and (i) of Lemma \ref{lem:turning-continuity}, there exists a constant $\beta_2\in (\beta_-,\beta_+)$ sufficiently close to $\beta_+$ such that
\begin{equation}\label{eq:right-order}
 M_P(\beta_2)<M_Q(\beta_2).
\end{equation}
Since by Lemma \ref{lem:turning-continuity} the functions $M_P(\beta)$ and $M_Q(\beta)$ are continuous on $(\beta_-,\beta_+)$, by \eqref{eq:left-order}, \eqref{eq:right-order}, the intermediate value theorem, \eqref{eq:MP-finite}, \eqref{eq:MQ-finite} and (iii) of Lemma \ref{lem:degenerate-sign}, there exists a constant $\widehat\beta\in(\beta_1,\beta_2)$ such that \eqref{eq:matched-M} holds.

\end{proof}

By Theorem \ref{thm-existence-self-similar1} there exist a unique constant
$R_o=R_{\mathrm o}(\widehat\beta,1)\in(0,\infty]$
and a unique maximal positive radially symmetric solution
\[
f^o_{\widehat\beta}\in C^1([0,R_o))\cap C^2((0,R_o))
\]
of \eqref{f-elliptic-eqn} in $B_{R_o}$ with $\alpha=\frac{2\widehat\beta+1}{1-m}$ and $\beta=\widehat{\beta}$,
satisfying \eqref{f-fr-origin-value} with $\eta_0=1$.
By Theorem \ref{thm-existence-self-similar2}  there exists a unique constant $R_g=R_g(\widehat{\beta},1)\in(0,\infty]$ and a unique maximal positive solution $g_{\widehat\beta}$ of the inverted equation \eqref{eq-inversion} in $[0,R_g)$ with $\widetilde{\alpha}$, $\widetilde{\beta}$ given by \eqref{q-and-tilde-def-intro} with $\alpha=\frac{2\widehat\beta+1}{1-m}$ and $\beta=\widehat{\beta}$, satisfying $g_{\widehat{\beta}}(0)=1$.
Moreover the function $f^{\infty}_{\widehat\beta}$ given by 
\eqref{outer-profile-definition-intro}
is the unique maximal positive outer radial branch solution of \eqref{f-elliptic-eqn} in $\mathbb{R}^n\setminus B_{R_g^{-1}}$ satisfying \eqref{f-infty-behaviour2} with $\eta=1$.

Let  $W_P:=W^o_{\widehat\beta}$ and $W_Q:=W^{\infty}_{\widehat\beta}$ be given by \eqref{eq:W-origin-infinity}.   Let $P_{\widehat\beta}$ and $Q_{\widehat\beta}$ be defined by \eqref{P-beta-defn} and \eqref{Q-beta-defn} respectively with $\beta=\widehat\beta$. Then \eqref{P-defn} and \eqref{Q-defn} hold.

\begin{lem}\label{lem:W1-=s-t-star}
Let $t_P:=\tau_P(\widehat\beta)$ be given by \eqref{tau-P-definition}. Then $t_P\in\mathbb{R}$ and 
\begin{equation}\label{W1-=s-star}
W_P(t_P)=s_*.
\end{equation}

\end{lem}
\begin{proof}
By (iii) of Lemma \ref{lem:turning-finite},

\begin{equation*}
P_{\widehat\beta}(s)=\sqrt{2h(s_*)(s_*-s)}+o(\sqrt{s_*-s}).
 \end{equation*}
Then there exists a constant $0<\delta_3<s_*/2$ such that 
\begin{equation} \label{eq:t-star-lower-P}
 P_{\widehat\beta}(s)\ge\sqrt{h(s_*)(s_*-s)}\quad\forall s_*-\delta_3\le s\le s_*.
 \end{equation}
Let $s_1=s_*-\delta_3$. Then there exists a constant $t_1$ such that $W_P(t_1)=s_1$ and by \eqref{P-defn}, \eqref{eq:t-star-lower-P},

\begin{equation*}
\frac{ds}{dt}=P_{\widehat\beta}(s)\quad\forall s\in(0,M_P(\widehat\beta))\quad\Rightarrow\quad t_P:=t_1+\int_{s_1}^{s_*}\frac{d\xi}{P_{\widehat\beta}(\xi)}<\infty
\end{equation*}
with  $W_P(t_P)=s_*$ and \eqref{W1-=s-star} follows.
\end{proof}

By an argument similar to the proof of Lemma \ref{lem:W1-=s-t-star} but with 
 (iv) of Lemma \ref{lem:turning-finite} and \eqref{Q-beta-defn} replacing (iii) of Lemma \ref{lem:turning-finite} and \eqref{P-beta-defn} in the proof we get the following result.

\begin{lem}\label{lem:W2-=s-t-star}
Let $t_Q:=\tau_Q(\widehat\beta)$ be given by \eqref{tau-Q-definition}. Then $t_Q\in\mathbb{R}$ and 
\begin{equation}\label{W2-=s-star}
W_Q(t_Q)=s_*.
\end{equation}

\end{lem}

We now let 
\begin{equation}\label{hat-f-defn}
f_P(r)=e^{\frac{2t_P}{1-m}}f^o_{\widehat{\beta}}(e^{t_P}r)\quad\forall 0\le r\le 1
\end{equation}
and
\begin{equation}\label{hat-W1-defn}
\widehat{W}_P(t)=e^{\frac{2mt}{1-m}}f_P(e^t)^m\quad\forall t\le 0.
\end{equation}
Then $f_P$ is a  radially symmetric solution of \eqref{f-elliptic-eqn} in $B_1$ with $\alpha=\frac{2\widehat\beta+1}{1-m}$, $\beta=\widehat\beta$, which satisfies
\begin{equation}\label{hat-f-origin-value}
f_P(0)=e^{\frac{2t_P}{1-m}}\quad\mbox{ and }\quad f_{P,r}(0)=0
\end{equation}
and by \eqref{eq:W-origin-infinity}, \eqref{hat-f-defn} and \eqref{hat-W1-defn}
\begin{equation*}
\widehat{W}_P(t)=e^{\frac{2m(t+t_P)}{1-m}}f^o_{\widehat{\beta}}(e^{t+t_P})^m=W_P(t+t_P)\quad\forall t\le 0.
\end{equation*}
Hence
\begin{equation}\label{widehat-W1-at-0=0}
\widehat{W}_P(0)=W_P(t_P)=s_*,
\end{equation} 
\begin{equation}\label{widehat-W1'-positive}
\widehat{W}_P'(t)=W_P'(t+t_P)>0\quad\forall t<0
\end{equation}
and
\begin{equation}\label{widehat-W1'-0=0}
\widehat{W}_P'(0)=W_P'(t_P)=P_{\widehat{\beta}}(s_*)=0.
\end{equation}
Similarly let 
\begin{equation}\label{hat-f2-defn}
f_Q(r)=e^{\frac{2t_Q}{1-m}}f^{\infty}_{\widehat{\beta}}(e^{t_Q}r)\quad\forall r\ge 1
\end{equation}
and
\begin{equation}\label{hat-W2-defn}
\widehat{W}_Q(t)=e^{\frac{2mt}{1-m}}f_Q(e^t)^m\quad\forall t\ge 0.
\end{equation}
Then $f_Q$ is the unique maximal positive outer radial branch of \eqref{f-elliptic-eqn} in $\mathbb{R}^n\setminus B_1$ with $\beta=\widehat\beta$, $\alpha=\frac{2\widehat\beta+1}{1-m}$, which satisfies 
\begin{equation}\label{hat-f2-infty-behaviour2} 
\lim_{r\to\infty}r^{\frac{n-2}{m}}f_Q(r)=e^{-\frac{n-2-nm}{m(1-m)}t_Q}\lim_{r\to\infty}(e^{t_Q}r)^{\frac{n-2}{m}}f^{\infty}_{\widehat{\beta}}(e^{t_Q}r)=e^{-\frac{n-2-nm}{m(1-m)}t_Q}
\end{equation}
and by \eqref{eq:W-origin-infinity}, \eqref{hat-f2-defn} and \eqref{hat-W2-defn},
\begin{equation*}
\widehat{W}_Q(t)=e^{\frac{2m(t+t_Q)}{1-m}}f^{\infty}_{\widehat{\beta}}(e^{t+t_Q})^m=W_Q(t+t_Q)\quad\forall t\ge 0.
\end{equation*}
Hence 
\begin{equation}\label{widehat-W2-at-0=0}
\widehat{W}_Q(0)=W_Q(t_Q)=s_*,
\end{equation} 
\begin{equation}\label{widehat-W2'-negative}
\widehat{W}_Q'(t)=W_Q'(t+t_Q)<0\quad\forall t>0
\end{equation}
and
\begin{equation}\label{widehat-W2'-0=0}
\widehat{W}_Q'(0)=W_Q'(t_Q)=-Q_{\widehat{\beta}}(s_*)=0.
\end{equation}
Let
\begin{equation}\label{glued-W-definition}
 W(t)=
 \begin{cases}
  \widehat W_P(t),&t\leq0,\\
  \widehat W_Q(t),&t\geq0.
 \end{cases}
\end{equation}

\begin{lem}[Smooth gluing]\label{lem:smooth-gluing}
$W\in C^\infty(\R)$ is positive and satisfies \eqref{eq:lienard} on $\R$ with $\beta=\widehat\beta\in(\beta_-,\beta_+)$ and
\begin{equation}\label{glued-W-monotonicity}
 W'(t)>0\quad\forall t<0,
 \qquad W'(0)=0,
 \qquad W'(t)<0\quad\forall t>0.
\end{equation}
Moreover there exists a radially symmetric function $f\in C^{\infty}(\mathbb{R}^n)$ which satisfies \eqref{eq:f-from-W} and \eqref{f-elliptic-eqn} in $\mathbb{R}^n$ with $\beta=\widehat{\beta}$, $\alpha=\frac{2\widehat\beta+1}{1-m}$,
\begin{equation}\label{f-origin-value6}
f(0)=e^{\frac{2t_P}{1-m}}\quad\mbox{ and }\quad f_r(0)=0
\end{equation}
and 
\begin{equation}\label{f-infty-behaviour6} 
\lim_{r\to\infty}r^{\frac{n-2}{m}}f(r)=e^{-\frac{n-2-nm}{m(1-m)}t_Q}
\end{equation} 
for some constants $t_P$, $t_Q$.  In addition, $f_r(r)<0$ for every $r>0$.

\end{lem}
\begin{proof}
By \eqref{hat-W1-defn},\eqref{widehat-W1-at-0=0},  \eqref{hat-W2-defn} and \eqref{widehat-W2-at-0=0},
\begin{equation}\label{W12-at-zero=0}
\widehat{W}_P(0)=\widehat{W}_Q(0)=f_P(1)^m=f_Q(1)^m=s_*.
\end{equation} 
By \eqref{widehat-W1'-0=0} and \eqref{widehat-W2'-0=0},
\begin{equation}\label{widehat-W1-W2'-0=0}
\widehat{W}_P'(0)=\widehat{W}_Q'(0)=0.
\end{equation}
By Lemma \ref{lem:lienard-reduction} both $\widehat{W}_P$ and $\widehat{W}_Q$ satisfy \eqref{eq:lienard} with \(\beta=\widehat\beta\).
Hence by putting $t=0$, \(\beta=\widehat\beta\) and $W=\widehat{W}_P, \widehat{W}_Q$, in \eqref{eq:lienard} we get,
\begin{equation}\label{widehat-W12''-equal-origin=0}
\widehat{W}_P''(0)=\widehat{W}_Q''(0)=-h(s_*).
\end{equation} 
By \eqref{glued-W-definition}, \eqref{W12-at-zero=0}, \eqref{widehat-W1-W2'-0=0} and  \eqref{widehat-W12''-equal-origin=0}, $W\in C^2(\R)$ is positive and satisfies \eqref{eq:lienard} on $\R$ with $\beta=\widehat\beta\in(\beta_-,\beta_+)$. By differentiating 
\eqref{eq:lienard} we can express $W'''$ as a function of $W, W'$ and $W''$. Hence $W\in C^3(\mathbb{R})$. By repeating the argument we get that $W\in C^{\infty}(\R)$.
By \eqref{glued-W-definition}, \eqref{widehat-W1'-positive}, \eqref{widehat-W2'-negative} and \eqref{widehat-W1-W2'-0=0} we get \eqref{glued-W-monotonicity}.

By \eqref{W12-at-zero=0} we can  let 
\begin{equation}\label{f-new-defn}
f(e^t)=\left\{\begin{aligned}
&f_P(e^t)\quad\forall t\le 0\\
&f_Q(e^t)\quad\forall t\ge 0
\end{aligned}\right.
\end{equation}
Then $f$ satisfies \eqref{eq:f-from-W}. By \eqref{hat-f-origin-value} and \eqref{hat-f2-infty-behaviour2}, we get  \eqref{f-origin-value6}
and \eqref{f-infty-behaviour6}. 
Since $W\in C^{\infty}(\mathbb{R})$ solves
\eqref{eq:lienard} with \(\beta=\widehat\beta\), $f\in C^{\infty}(\mathbb{R}\setminus\{0\})$. By Lemma \ref{lem:lienard-reduction}, $f$ satisfies \eqref{f-ode} with $\beta=\widehat{\beta}$ and $\alpha=\frac{2\widehat\beta+1}{1-m}$. Since $f=f_P$ satisfies \eqref{f-elliptic-eqn} in a neighbourhood of the origin with these parameters, $f\in C^{\infty}(\mathbb{R}^n)$ and satisfies \eqref{f-elliptic-eqn} throughout $\mathbb{R}^n$.

Finally, differentiating \eqref{eq:W-def} gives
\[
 \frac{W'(t)}{W(t)}=a+m\frac{r f_r(r)}{f(r)},
 \qquad r=e^t.
\]
For $t<0$, Lemma \ref{lem:linear-comparison} gives $0<W'(t)=P_{\widehat\beta}(W(t))<aW(t)$; for $t>0$, \eqref{glued-W-monotonicity} gives $W'(t)<0<aW(t)$; and at $t=0$ one has $W'(0)=0<aW(0)$.  Thus $r f_r(r)/f(r)<0$ for every $r>0$, proving $f_r(r)<0$.

\end{proof}

\begin{rmk}\label{rmk:Phi-normalization}
Let $f$ be given by Lemma \ref{lem:smooth-gluing}, \(\lambda=f(0)^{-\frac{1-m}{2}}\) and 
\begin{equation}
 \Phi(r):=\lambda^{\frac{2}{1-m}}f(\lambda r).
 \label{eq:Phi-normalization}
\end{equation}
Then
\begin{equation}
 \Phi(0)=1,
 \qquad \Phi_r(0)=0,
 \qquad \Phi_r(r)<0\quad\forall r>0,
 \label{eq:Phi-properties}
\end{equation}
and
\begin{equation}\label{eq:Phi-properties2}
\lim_{r\to\infty}r^{\frac{n-2}{m}}\Phi(r)=\lambda^{-\frac{n-2-nm}{m(1-m)}}\lim_{r\to\infty}(\lambda r)^{\frac{n-2}{m}}f(\lambda r)=f(0)^{\frac{n-2-nm}{2m}}e^{-\frac{n-2-nm}{m(1-m)}t_Q}=:C_*.
\end{equation}
The function \(\Phi\) satisfies \eqref{f-elliptic-eqn} in $\mathbb{R}^n$ with $\beta=\widehat{\beta}$, $\alpha=\frac{2\widehat\beta+1}{1-m}$.

\end{rmk}

We are now ready to prove Theorem \ref{thm-existence-self-similar3}.

\noindent $\text{\bf Proof of Theorem \ref{thm-existence-self-similar3}}$: Let \(\rho _1>0\) and \(\eta _0>0\). Let $\widehat{\beta}\in(\beta_-,\beta_+)$ and $\Phi$ be as in Remark \ref{rmk:Phi-normalization}.  Let
\[
 \widehat\alpha:=\frac{2\widehat\beta+1}{1-m}
\]
and
\begin{equation}
 \beta_*:=\rho _1\widehat\beta,
 \qquad
 \alpha_*:=\rho _1\frac{2\widehat\beta+1}{1-m}
           =\frac{2\beta_*+\rho _1}{1-m}.
 \label{eq:scaled-parameters}
\end{equation}
Then $\beta_*\in (\beta_-(\rho _1),\beta_+(\rho _1))$.
We define
\begin{equation}
 k:=\bigl(\rho _1\eta _0^{1-m}\bigr)^{1/2},
 \qquad
 f(x):=\eta _0\Phi(kx).
 \label{eq:final-profile}
\end{equation}
With \(y=kx\),
\begin{align*}
 \Delta\!\left(\frac{f^m}{m}\right)(x)
 &=\eta _0^m k^2
   \Delta\!\left(\frac{\Phi^m}{m}\right)(y)
 =\rho _1\eta _0
   \Delta\!\left(\frac{\Phi^m}{m}\right)(y),
\end{align*}
and
\[
 \alpha_*f+\beta_*x\cdot\nabla f
 =\rho _1\eta _0
 \left(\widehat\alpha\Phi+
       \widehat\beta y\cdot\nabla\Phi\right)(y).
\]
Thus \(f\)  satisfies \eqref{f-elliptic-eqn} in $\mathbb{R}^n$ with $\alpha=\alpha_*$ and $\beta=\beta_*$.  By \eqref{eq:Phi-properties},
\[
 f(0)=\eta _0,
 \qquad f_r(0)=0,
 \qquad f_r(r)<0\quad\forall r>0.
\]
Also by \eqref{eq:Phi-properties2},
\begin{align*}
 \lim_{r\to\infty}r^{\frac{n-2}{m}}f(r)
 &=C_*\eta _0 k^{-\frac{n-2}{m}}=C_*\rho _1^{-\frac{n-2}{2m}}
   \eta _0^{-\frac{n-2-nm}{2m}}
\end{align*}
and \eqref{f-infty-limit} follows. For uniqueness, let $\widetilde f$ be
another positive radial solution with the same $\beta_*$ and the same origin
data.  Lemma \ref{lem-local-existence} gives $\widetilde f=f$ near the
origin, and ordinary ODE uniqueness continues the equality on every compact
interval on which the profiles are positive.  Since both are positive on
$[0,\infty)$, they agree everywhere and  Theorem \ref{thm-existence-self-similar3} is proved.

\hfill\(\square\)

\end{document}